\documentclass[11pt]{article}
\usepackage{amsfonts}
\usepackage{graphics}
\usepackage{indentfirst}
\usepackage{color}
\usepackage{cite}
\usepackage{latexsym}
\usepackage[paper=a4paper, left=1.8cm, right=1.8cm, top=1.8cm, bottom=1.6cm, headheight=5.5pt, footskip=0.8cm, footnotesep=0.8cm, centering, includefoot]{geometry}
\usepackage{amsmath}
\allowdisplaybreaks
\usepackage{amssymb}
\usepackage[colorlinks, linkcolor=red]{hyperref}
\hypersetup{urlcolor=red, citecolor=red}
\usepackage[dvips]{epsfig}
\usepackage{amscd}

\newtheorem{theorem}{Theorem}[section]
\newtheorem{remark}{Remark}[section]

\newtheorem{lemma}{Lemma}[section]
\newtheorem{corollary}{Corollary}[section]
\newtheorem{proposition}{Proposition}[section]

\DeclareMathOperator{\curl}{curl}

\DeclareMathOperator{\divv}{div}

\makeatletter
\@addtoreset{equation}{section}
\makeatother
\makeatletter
\@addtoreset{equation}{section}
\makeatother

\title{Global well-posedness and large time behavior for compressible non-isothermal nematic liquid crystal flows with vacuum at infinity
\thanks{Y. Liu was supported by National Natural Science Foundation of China (No. 11901288), Scientific Research Foundation of Jilin Provincial Education Department (No. JJKH20210873KJ), Natural Science Foundation of Changchun Normal University, and Postdoctoral Science Foundation of
China (No. 2021M691219). X. Zhong was supported by National Natural Science Foundation of China (Nos. 11901474, 12071359), the Chongqing Talent Plan for Young Topnotch Talents (No. CQYC202005074), and the Innovation Support Program for Chongqing Overseas Returnees (No. cx2020083).
}
}

\author{Yang Liu$\,^{\rm 1,2}\,$\quad Xin Zhong$\,^{\rm 3}\,${\thanks{Corresponding author. E-mail address:
 liuyang0405@ccsfu.edu.cn (Y. Liu),  xzhong1014@amss.ac.cn (X. Zhong)}}
\date{}\\
\footnotesize $^{\rm 1}\,$
School of Mathematics, Jilin University, Changchun 130012, P. R. China\\
\footnotesize $^{\rm 2}\,$ College of Mathematics, Changchun Normal
University, Changchun 130032, P. R. China\\
\footnotesize $^{\rm 3}\,$ School of Mathematics and Statistics, Southwest University, Chongqing 400715, P. R. China}

\date{ }

\begin{document}
\maketitle

\begin{abstract}
We investigate the Cauchy problem of three-dimensional compressible non-isothermal nematic liquid crystal flows in $\mathbb{R}^3$. We derive the global existence and uniqueness of strong solutions with both interior and far field vacuum states provided that the initial data are of small total energy. This improves our previous work \cite{LZ21} in the sense that the initial data allow possibly large oscillations. The key analysis is based on delicate energy estimates and the structural characteristic of the system under consideration. Moreover, we also show the algebraic decay estimates of the solution. The results could be viewed as an extension of the studies in Li-Xu-Zhang \cite{LXZ18} for the isothermal case to the non-isothermal situation.
\end{abstract}



\section{Introduction and main results}

The motion of compressible nematic liquid crystal flows
in the three-dimensional (3D) space can be written as the following simplified non-isothermal Ericksen-Leslie model (see \cite{DL19,HP16,HP18,HP19}):
\begin{align}\label{a1}
\begin{cases}
 \rho_t+\divv(\rho u)=0,\\
\rho u_t+\rho u\cdot\nabla u-\mu\Delta u-(\lambda+\mu)\nabla \divv u+\nabla p=-\nabla d\cdot\Delta d,\\
c_v\rho(\theta_t+u\cdot\nabla\theta)+p\divv u-\kappa\Delta\theta
=\lambda(\divv u)^2+\frac{\mu}{2}|\nabla u+(\nabla u)^{tr}|^2
+|\Delta d+|\nabla d|^2d|^2,\\
d_t+u\cdot\nabla d=\Delta d+|\nabla d|^2d, \\
|d|=1,
\end{cases}
\end{align}
where the unknowns $\rho$, $u=(u^1,u^2,u^3)$, $\theta$, and $d\in\Bbb S^2$ are
the density, the velocity, the absolute temperature, and the macroscopic average of the nematic liquid crystal orientation field, respectively. $p=R\rho\theta$, with positive constant $R$, is the pressure.
The constant viscosity coefficients $\mu$ and $\lambda$ satisfy the physical restrictions
\begin{align}\label{a0}
\mu>0, \ 2\mu+3\lambda\ge 0.
\end{align}
The positive constants $c_v$ and $\kappa$ are the heat capacity and the ratio of the heat conductivity coefficient
over the heat capacity, respectively.

We consider the Cauchy problem of \eqref{a1} with the initial condition
\begin{align}\label{a2}
(\rho, u, \theta, d)|_{t=0}=(\rho_0, u_0, \theta_0, d_0), \quad x\in\Bbb R^3,
\end{align}
and the far field behavior
\begin{align}\label{a3}
(\rho,u,\theta,d)(x, t)\rightarrow (0, 0, 0,\mathbf{1}), ~{\rm as}~|x|\rightarrow\infty, ~t>0,
\end{align}
where $\mathbf{1}$ is a given unit vector and $|d_0|=1$.

In 2011, Feireisl-Rocca-Schimperna \cite{FRS2011} introduced a
non-isothermal model for incompressible nematic liquid crystals in the spirit of the simplified version of the Leslie-Ericksen model proposed by Lin and Liu
\cite{LL95}. Later on, Feireisl-Fr{\'e}mond-Rocca-Schimperna \cite{FFRS2012} proposed a new approach to the modeling of non-isothermal liquid crystals and investigated global-in-time weak solutions for a Ginzburg-Landau approximation system to relax the nonlinear constraint $|d|=1$. Since then much attention has been paid to various non-isothermal nematic liquid crystals models. Li and Xin \cite{LX16} proved the global existence of weak solutions on $\mathbb{T}^2$. Yet the uniqueness and regularity of such weak solutions remain open. Bian and Xiao \cite{BX17} showed the global strong solution to the non-isothermal nematic liquid crystal flows with temperature-dependent viscosity provided that the initial data is a sufficiently small perturbation around the trivial equilibrium state. De Anna and Liu \cite{DL19} studied the general non-isothermal Ericksen-Leslie system and obtained the global existence of strong solutions in Besov spaces provided the Besov norm of the initial data is sufficiently small. Du-Li-Wang \cite{DLW20} concerned the simplified non-isothermal Ericksen-Leslie system of nematic
liquid crystal flow on a three-dimensional smooth bounded domain. Existence of
a global weak solution is proven under the homogeneous Dirichlet boundary condition on the velocity field, the Neumann condition on the director field, and the non-flux condition on the heat flux, and for suitable initial data.
It is worth mentioning that some mathematical analysis concerning the hydrodynamic flow of non-isothermal nematics in the $Q$-tensor formalism were addressed
by Feireisl-Rocca-Schimperna-Zarnescu \cite{FRSZ2014,FRSZ2015}. For more background and theoretical results about nematic liquid crystals both in the
isothermal and non-isothermal situation, we refer the reader to the monograph \cite{SV12} and survey articles \cite{LW14,HP18}.

Recently there has been considerable interest in the mathematical study
for the compressible non-isothermal nematic liquid crystals. In 2016, Hieber and Pr{\"u}ss \cite{HP16} derived a general compressible non-isothermal Ericksen-Leslie model for the flow of nematic liquid crystals and analyzed the simplified non-isothermal model mathematically. Later, they \cite{HP19} showed the local well-posedness with in $L^q$ setting and the global solution for initial data close to equilibria points without structural assumptions on the Leslie coefficients in bounded domains. At the same time, De Anna and Liu \cite{DL19} proposed a non-isothermal model describing the evolution of a nematic liquid-crystal material under the action of thermal effects in the whole space.
On the other hand, for the simplified compressible non-isothermal model \eqref{a1},
Fan-Li-Nakamura \cite{FLN18} proved the local well-posedness of strong solutions to the initial boundary value problem under the compatibility conditions \eqref{qqw1} and \eqref{qqw} for the initial data. Then Liu and Zhong \cite{LZ21} extended it to be a global one for the problem \eqref{a1}--\eqref{a3} provided that $\|\rho_0\|_{L^\infty}+\|\nabla d_0\|_{L^3}$ is suitably small and $3\mu>\lambda$. While, for the absence of vacuum states, Guo-Xi-Xie \cite{GXX17} studied global existence and long time behavior of smooth solutions to \eqref{a1}--\eqref{a3} provided that the initial data are close to the constant equilibrium state in $H^3$-framework. When the term $|\nabla d|^2$ in $\eqref{a1}_3$ and $\eqref{a1}_4$ is replaced by the Ginzburg-Landau approximation term $\frac{1-|d|^2}{\varepsilon}$, Guo-Xie-Xi \cite{GXX16} established the global existence of weak solutions. The global well-posedness of strong solutions to the one-dimensional problem of \eqref{a1} can be found in \cite{M20,TS20}.
There are also very interesting mathematical results concerning the global existence of (weak, strong or classical) solutions of compressible isothermal nematic liquid crystal flows, please refer to \cite{W16,JJW2013,JJW14,LXZ18,LLW15}.

In this paper, we further investigate the global existence and uniqueness of strong solutions to the problem \eqref{a1}--\eqref{a3}. Our result reveals that the small initial energy is sufficient to guarantee the global well-posedness. Moreover, we also show that the solutions converge algebraically to zero in $L^2$ as time goes to infinity. This result could be viewed as the first one on time-decay rates of solutions to the Cauchy problem of multi-dimensional compressible non-isothermal nematic liquid crystals with vacuum. The main novelty of the present article consists in allowing possibly large oscillations of the initial data as well as obtaining algebraic decay estimates of the solution.

Before formulating our main result, we first explain the notations and conventions used throughout this paper.
 For simplicity, in what follows, we denote
\begin{align*}
\int_{\Bbb R^3}fdx=\int fdx.
\end{align*}
For $1\le p\le \infty$ and integer $k\ge 0$, the standard homogeneous and inhomogeneous
Sobolev spaces as follows:
\begin{align*}
\begin{cases}
L^p=L^p(\Bbb R^3),~ ~ W^{k, p}=L^p\cap D^{k, p},~~ H^k=W^{k, 2}, \\
  D^{k, p}=\{u\in L_{loc}^1(\Bbb R^3): \|\nabla^ku\|_{L^p}<\infty\}, ~D^k=D^{k, 2},\\
  D_0^1=\{u\in L^6(\Bbb R^3): \|\nabla u\|_{L^2}<\infty\}.
\end{cases}
\end{align*}

We begin with the global existence of strong solutions to the problem \eqref{a1}--\eqref{a3} with initial data being subjected to the
small energy assumption.
\begin{theorem}\label{thm1}
 For $q\in (3, 6)$ and positive constant $\bar\rho$, assume that the initial data $(\rho_0, u_0, \theta_0\ge 0, d_0)$ satisfies
\begin{align}
& \rho_0\in L^1\cap H^1\cap W^{1, q}, ~\big(\rho_0|u_0|^2,\rho_0\theta_0,|\nabla d_0|^2\big)\in L^1,\label{qq1} \\
& \big(\nabla u_0, \nabla\theta_0\big)\in H^1,\ \nabla d_0\in H^2,\ |d_0|=1,\\ \label{1.6}
& 0\le \inf\rho_0(x)\le \sup\rho_0(x)\le \bar{\rho},
\end{align}
and the compatibility conditions
\begin{align}
& \mu\Delta u_0+(\mu+\lambda)\nabla\divv u_0-\nabla(R\rho_0\theta_0)-\nabla d_0\cdot\Delta d_0=\sqrt{\rho_0}g_1,\label{qqw1}\\
& \kappa\Delta\theta_0+\lambda(\divv u_0)^2+\frac{\mu}{2}|\nabla u_0+(\nabla u_0)^{tr}|^2+|\Delta d_0+|\nabla d_0|^2d_0|^2=\sqrt{\rho_0}g_2,\label{qqw}
\end{align}
for some $g_1, g_2\in L^2$.
There exists a small positive
constant $\varepsilon_0$ depending only on  $\mu$, $\lambda$, $c_v$, $R$, $\kappa$, $\bar{\rho}$, and the initial data such that if
\begin{align}\label{1.9}
\int\Big(\frac12\rho_0|u_0|^2+c_v\rho_0\theta_0+
\frac12|\nabla d_0|^2\Big)dx\triangleq\Bbb E_0\le \varepsilon_0,
\end{align}
then the problem \eqref{a1}--\eqref{a3} has a unique global strong solution $(\rho, u, \theta, d)$ satisfying
\begin{align}\label{trr}
\begin{cases}
\rho \in L^\infty\big([0,\infty);L^1\cap L^\infty\big)\cap C\big([0, \infty); L^1\cap H^1\cap W^{1, q}\big),\\
(u, \theta)\in C\big([0, \infty); D_0^1\cap D^2\big)\cap L_{loc}^2\big([0, \infty); D^{2, q}\big),\\
(u_t, \theta_t)\in L_{loc}^2\big([0, \infty); D_0^1\big), \ \big(\sqrt{\rho}u_t, \sqrt{\rho}\theta_t\big)\in L_{loc}^\infty\big([0, \infty); L^2\big),\\
\nabla d\in L^\infty\big([0, \infty); H^2\big)\cap C\big([0, \infty); H^2\big)\cap L^2\big([0, \infty); H^3\big), \\
d_t\in C([0, \infty); H^1)\cap L^2\big([0, \infty); H^2\big).
\end{cases}
\end{align}
\end{theorem}

\begin{remark}
Although our Theorem \ref{thm1} has small energy, yet the oscillations of the initial data could be arbitrarily large. Thus we improve our previous work \cite{LZ21} where the initial density must be small pointwisely. Moreover, we remove the additional restriction $3\mu>\lambda$ used in \cite{LZ21}.
\end{remark}

\begin{remark}
Compared with the 3D isothermal case \cite{LXZ18}, there is no need to assume that
$u_0,\nabla d_0\in \dot{H}^{\beta}$ with $\beta\in(\frac12,1]$ for the global existence of strong solutions.
\end{remark}

The proof of Theorem \ref{thm1} is based on delicate energy estimates.
The key point in our analysis is to derive the time-independent upper bound of the density. Compared with the 3D isothermal situation \cite{LXZ18}, as pointed out in \cite{LZ21}, the basic energy estimate provides us
useless information on dissipation estimates of $u$, $\theta$, and $d$. To overcome this difficulty, motivated by \cite{L21}, we introduce a small time $t_0$ and recover the crucial dissipation estimate of the form
$\int_0^T\big(\|\nabla u\|_{L^2}^2+\|\nabla \theta\|_{L^2}^2+\|\nabla^2d\|_{L^2}^2\big)dt$ in terms of a power function of $\Bbb E_0$ (see Lemma \ref{l31}) rather than $\|\rho_0\|_{L^\infty}+\|\nabla d_0\|_{L^3}$ used in \cite{LZ21}. Moreover, we point out that the method used here to derive the $L^\infty(0,T;L^2)$ estimates of $\nabla u$ and $\nabla^2d$ is different from that of in \cite{LZ21}. In order to obtain this control, the key observation is to treat the temperature simultaneously by introducing a cut-off function $\sigma(t)\triangleq\min\{t,t_0\}$ (see Lemma \ref{l33}) instead of $E\triangleq\frac{|u|}{2}+c_v\theta$, the sum of specific kinetic energy and specific internal energy. Furthermore, unlike \cite{LZ21}, it is impossible to get the
smallness of $\|\nabla d\|_{L^3}$ which plays a critical role in dealing with the coupling terms $u\cdot\nabla d$ and $\nabla d\cdot\Delta d$.
Next, since the Zlotnik inequality which was used in \cite{LXZ18}
for isothermal flow to get the upper bound of the density does not work here, we adopt some ideas introduced by \cite{D1997,WZ17} to establish the time-independent upper bound of the density. Finally, having the time-independent and time-dependent
\textit{a priori} estimates at hand, together with the local well-posedness theory, we succeed in showing Theorem \ref{thm1} by an argument of contradiction.

Next we derive the algebraic decay estimates of the global solution obtained in Theorem \ref{thm1}.
\begin{theorem}\label{thm2}
Let $(\rho, u, \theta, d)$ be the global solution obtained in Theorem \ref{thm1},
then there exists a constant $C$ depending only on $\mu$, $\lambda$, $c_v$, $\kappa$, $R$, $\bar{\rho}$, and the initial data such that, for large time $t$,
\begin{align}\label{1.11}
&\|\sqrt{\rho}u\|_{L^2}+\|\sqrt{\rho}\theta\|_{L^2}
+\|\nabla u\|_{L^2}+\|\nabla\theta\|_{L^2}+\|\sqrt{\rho}\dot{u}\|_{L^2}\le Ct^{-\frac14},\\ \label{1.12}
& \|\sqrt{\rho}\dot{\theta}\|_{L^2}+\|\nabla^2\theta\|_{L^2}
+\|\nabla^2d\|_{L^2}+\|\nabla d_t\|_{L^2}+\|\nabla^3d\|_{L^2}\le Ct^{-\frac12},
\end{align}
provided that
$\Bbb E_0\le \varepsilon_0'$ for some constant $\varepsilon_0'\in (0, \varepsilon_0]$.
\end{theorem}

\begin{remark}
Let $F$ and $w$ (see \eqref{f3}) be the effective viscous flux and the vorticity, respectively. Then it follows from \eqref{w5}, Sobolev's inequality, \eqref{trr}, \eqref{1.11}, and \eqref{1.12} that, for large time $t$,
\begin{align}\label{1.13}
\|\nabla F\|_{L^2}+\|\nabla w\|_{L^2}
& \le C\big(\|\rho\dot{u}\|_{L^2}+\||\nabla d||\nabla^2d|\|_{L^2}\big) \notag \\
& \leq C\|\rho\|_{L^\infty}^{\frac12}\|\sqrt{\rho}\dot{u}\|_{L^2}
+C\|\nabla d\|_{L^\infty}\|\nabla^2d\|_{L^2} \notag \\
& \leq C\|\rho\|_{L^\infty}^{\frac12}\|\sqrt{\rho}\dot{u}\|_{L^2}
+C\|\nabla d\|_{H^2}\|\nabla^2d\|_{L^2}
\notag \\
& \leq C\|\sqrt{\rho}\dot{u}\|_{L^2}+C\|\nabla^2d\|_{L^2} \notag \\
& \leq Ct^{-\frac14}.
\end{align}
This combined with \eqref{w7} and \eqref{1.11} yields that, for large time $t$,
\begin{align}\label{1.14}
\|\nabla u\|_{L^6}
\le C\big(\|\rho\dot{u}\|_{L^2}+\|\nabla\theta\|_{L^2}
+\||\nabla d||\nabla^2d|\|_{L^2}\big)\leq Ct^{-\frac14}.
\end{align}
Furthermore, we derive from H\"older's inequality, Sobolev's inequality, \eqref{1.14}, \eqref{1.11}, and \eqref{1.12} that, for any $q\in[2,6]$ and large time $t$,
\begin{align*}
& \|\nabla u\|_{L^q} \leq \|\nabla u\|_{L^2}^{\frac{6-q}{2q}}
\|\nabla u\|_{L^6}^{\frac{3q-6}{2q}}\leq Ct^{-\frac14},\\
& \|\nabla\theta\|_{L^q} \leq \|\nabla\theta\|_{L^2}^{\frac{6-q}{2q}}
\|\nabla\theta\|_{L^6}^{\frac{3q-6}{2q}}
\leq C\|\nabla\theta\|_{L^2}^{\frac{6-q}{2q}}
\|\nabla^2\theta\|_{L^2}^{\frac{3q-6}{2q}}\leq Ct^{-\frac{5q-6}{8q}}.
\end{align*}
\end{remark}

\begin{remark}
It is worth mentioning that Guo-Xi-Xie \cite{GXX17} established the decay rates for second order spatial derivatives of the density and the velocity provided that the initial perturbation is small in $H^3$-norm and bounded in $L^q$-norm ($1\leq q<\frac65$). However, we should point out that, for the initial density allowing vacuum states, whether the second derivatives of the velocity field decay or not remains unknown. The main difficulty lies in deriving time-independent bounds on the gradient of density (see \eqref{t91}).
\end{remark}

\begin{remark}
Compared with decay rates of the isothermal model \cite[Theorem 1.3]{W16}, our decay rates of the solution is slower due to the complicated structure of the system \eqref{a1}. Consequently, the mathematical studies between the isothermal model and the non-isothermal case of the 3D Cauchy problem appear highly differences.
\end{remark}

We divide the proof of Theorem \ref{thm2} into six steps. As mentioned above, due to the structural characteristic of the system \eqref{a1}, the basic energy estimate provides useless dissipation effect on the solution. However, for small initial energy (see \eqref{x4}), we can derive
\begin{align}\label{1.16}
\frac{d}{dt}\big(\|\sqrt{\rho}u\|_{L^2}^2+\|\sqrt{\rho}\theta\|_{L^2}^2\big)
+\|\nabla u\|_{L^2}^2+\|\nabla \theta\|_{L^2}^2
\le C\|\nabla^2d\|_{L^2}\|\nabla^3d\|_{L^2}^2.
\end{align}
This reveals that the time-decay rate of the second derivatives of the orientation field plays a crucial role in the initial layer analysis. In step 1, we obtain that $\|\nabla^2d\|_{L^2}$ decays with the rate of $t^{-\frac12}$ (see \eqref{dd1}) by analyzing \eqref{a1}$_4$. In step 2, with the help of \eqref{1.6} and the result of step 1, we deduce that $\|\sqrt{\rho}u\|_{L^2}$ and $\|\sqrt{\rho}\theta\|_{L^2}$ decay with the rate of $t^{-\frac14}$ (see \eqref{5.5}). In step 3, owing to very slowly time-decay estimates on the initial layer analysis, we fail to show the decay rates of $\|\nabla u\|_{L^2}$ and $\|\nabla \theta\|_{L^2}$ directly and should require the initial energy to be much more small (see \eqref{5.15}).  In steps 4 and 5, we devote the decay rates to $\|\nabla d_t\|_{L^2}$ and $\|\sqrt{\rho}\dot{\theta}\|_{L^2}$ by delicate energy estimates and the results of steps 1--3.
Finally, in step 6, the decay estimates on the third derivatives of the orientation field and the temperature can be obtained by applying $L^2$-estimates for the elliptic system and uniformly bounds we have just derived.

This paper is organized as follows. In Section 2, we collect some known facts and elementary inequalities which will be used later. Section 3 is devoted to the global \textit{a priori} estimates. The proof of Theorem \ref{thm1} will be done in Section 4. Finally, we give the proof of Theorem \ref{thm2} in Section 5.

\section{Preliminarily}
In this section, we recall some known facts and elementary inequalities which will be used later.

We begin with the following local existence of a unique strong solution to the problem \eqref{a1}--\eqref{a3}, whose proof can be performed
by using similar strategies (with a minor modification) as those in \cite{FLN18}.
\begin{lemma}\label{l21}
Under the conditions \eqref{qq1}--\eqref{qqw}, there is a small $T_*>0$ such that the Cauchy problem \eqref{a1}--\eqref{a3} admits a unique solution $(\rho, u, \theta, d)$
over $\Bbb R^3\times[0, T_*]$ satisfying for some constant $M_*>1$ depending on $T_*$ and the initial data
\begin{align}\label{w1}
&\sup_{0\le t\le T_*}\big(\|\nabla u\|_{H^1}^2+\|\nabla\theta\|_{H^1}^2+\|\nabla d\|_{H^2}^2
+\|\sqrt{\rho}\dot{u}\|_{L^2}^2+\|\sqrt{\rho}\dot{\theta}\|_{L^2}^2+\|\nabla d_t\|_{L^2}^2\big)\nonumber\\
&+\int_0^{T_*}\big(\|\nabla\dot{u}\|_{L^2}^2+\|\nabla^2d_t\|_{L^2}^2
+\|\nabla\dot{\theta}\|_{L^2}^2\big)dt\le M_*.
\end{align}
\end{lemma}

Next, the following well-known Gagliardo-Nirenberg inequality (see \cite[Chapter 10, Theorem 1.1]{D2016}) will be
used frequently later.
\begin{lemma}\label{l22}
Assume that $f\in D^{1, m}\cap L^r$ with $m,r\geq1$, then there exists a constant $C$ depending only on $q$, $m$, and $r$ such that
\begin{align}
\|f\|_{L^q}\le C\|\nabla f\|_{L^m}^\vartheta\|f\|_{L^r}^{1-\vartheta},
\end{align}
where $\vartheta=\big(\frac{1}{r}-\frac{1}{q}\big)/
\big(\frac{1}{r}-\frac{1}{m}+\frac13\big)$ and the admissible range of $q$ is the following:
\begin{itemize}
\item if $m<3$, then $q$ is between $r$ and $\frac{3m}{3-m}$;
\item if $m=3$, then $q\in [r, \infty)$;
\item if $m>3$, then $q\in [r, \infty]$.
\end{itemize}
\end{lemma}

Next,  the material derivatives of $f$, the effective viscous flux $F$,
and the vorticity $w$ are defined as follows,
\begin{align}\label{f3}
\dot{f}\triangleq u_t+u\cdot\nabla f, \ F\triangleq(2\mu+\lambda)\divv u-p, \ w\triangleq\curl u.
\end{align}
Then one derives the following elliptic equations from $\eqref{a1}_2$
by taking the operators $\divv$ and $\curl$, respectively.
\begin{align}\label{s22}
\Delta F=\divv(\rho\dot{u})-\divv(\nabla d\odot\nabla d),\ \mu\Delta w=\curl(\rho\dot{u}+\divv(\nabla d\odot\nabla d)).
\end{align}

We now state some elementary $L^q$-estimates of the elliptic system \eqref{s22}, whose proof can be done by similar methods as those in \cite{L21}.
\begin{lemma}
Let $(\rho, u, \theta, d)$ be a smooth solution of \eqref{a1}--\eqref{a3}. Then for any $q\in[2,6]$, there exists a positive constant $C$ depending only on $\mu$, $\lambda$, and $q$ such that
\begin{align}
&\|\nabla F\|_{L^q}+\|\nabla w\|_{L^q}\le C\big(\|\rho\dot{u}\|_{L^q}+
\||\nabla d||\nabla^2 d|\|_{L^q}\big),\label{w5}\\
&\|\nabla u\|_{L^6}\le C\big(\|\rho\dot{u}\|_{L^2}+\|\nabla\theta\|_{L^2}+\||\nabla d||\nabla^2d|\|_{L^2}\big).\label{w7}
\end{align}
\end{lemma}

Finally, the following Beal-Kato-Majda type inequality (see \cite[Lemma 2.3]{HLX20112}) will be used to estimate $\|\nabla u\|_{L^\infty}$.
\begin{lemma}\label{l24}
Suppose that $q\in (3, \infty)$, then there exists a constant $C$ depending only on $q$ such that for any $\nabla u\in L^2\cap W^{1, q}$,
\begin{align}
\|\nabla u\|_{L^\infty}\le C\big(\|\divv u\|_{L^\infty}
+\|\curl u\|_{L^\infty}\big)\log\big(e+\|\nabla^2u\|_{L^q}\big)+C\big(1+\|\nabla u\|_{L^2}\big).
\end{align}
\end{lemma}

\section{A priori estimates}

In this section, we will establish some necessary {\it a priori} bounds for
smooth solutions to the Cauchy problem \eqref{a1}--\eqref{a3} to extend the local strong
solutions guaranteed by Lemma \ref{l21}. Thus, let $T>0$ be a fixed time and
$(\rho, u, \theta, d)$ be the smooth solution to \eqref{a1}--\eqref{a3} on $\Bbb R^3\times(0, T]$ with initial data
$(\rho_0, u_0, \theta_0, d_0)$ satisfying \eqref{qq1}--\eqref{qqw}. In what follows, $C, C_i$, and $c_i\ (i=1, 2, \cdots)$ denote generic positive
constants which rely only on $R$, $c_v$, $\mu$, $\lambda$, $\kappa$, $\bar{\rho}$, $M_*$, and the initial values, and $C(\alpha)$ is used to emphasize the dependence of $C$ on $\alpha$.

First of all, if we multiply $\eqref{a1}_1$ by a cut-off function and then use a standard limit procedure, we can obtain that (see also \cite[Lemma 3.1]{WZ17})
\begin{align}\label{lz3.1}
\|\rho(t)\|_{L^1}=\|\rho_0\|_{L^1}, \quad 0\leq t\leq T.
\end{align}
Applying standard maximum principle (see \cite[p. 43]{F2004}) to \eqref{a1} along with $\rho_0,\theta_0\geq0$ indicates that
\begin{equation}\label{lz3.2}
\inf_{\mathbb{R}^3\times[0,T]}\rho(x,t)\geq0,\ \ \inf_{\mathbb{R}^3\times[0,T]}\theta(x,t)\geq0.
\end{equation}
Multiplying $\eqref{a1}_2$ by $u$ and $\eqref{a1}_4$ by $\Delta d+|\nabla d|^2d$, respectively, and integrating the resulting equations over $\mathbb{R}^3$,
integrating $\eqref{a1}_3$ over $\mathbb{R}^3$, and summing them together, we obtain from integration by parts and $\eqref{a1}_5$ that
\begin{equation*}
\frac{d}{dt}\int\Big(\frac12\rho|u|^2+c_v\rho\theta+\frac12|\nabla d|^2\Big)dx
=0.
\end{equation*}
This implies that
\begin{equation}\label{lz3.3}
\int\Big(\frac12\rho|u|^2+c_v\rho\theta+\frac12|\nabla d|^2\Big)dx
=\int\Big(\frac12\rho_0|u_0|^2+c_v\rho_0\theta_0+\frac12|\nabla d_0|^2\Big)dx,
\quad 0\leq t\leq T.
\end{equation}

\subsection{Time-independent estimates}
Denote
\begin{align*}
& A_1(T)\triangleq\sup_{0\le t\le T}\big(\|\nabla u\|_{L^2}^2+\|\nabla^2d\|_{L^2}^2+\|\nabla\theta\|_{L^2}^2\big),\\
& A_2(T) \triangleq \int_0^T\big(\|\sqrt{\rho}\dot{u}\|_{L^2}^2+\|\nabla d_t\|_{L^2}^2
+\|\nabla^3d\|_{L^2}^2+\|\nabla\theta\|_{L^2}^2+\|\nabla^2\theta\|_{L^2}^2\big)dt,\\
& A_3(T)\triangleq\int_0^T\|\nabla^2d\|_{L^2}^2dt.
\end{align*}

\begin{proposition}\label{p1}
There exists some positive constant $\varepsilon_0$ depending only on
$\mu$, $\lambda$, $c_v$, $R$, $\kappa$, $\bar{\rho}$, and the initial data such that if $(\rho, u, \theta, d)$ is the smooth solution of \eqref{a1}--\eqref{a3} on $\mathbb{R}^3\times(0, T]$ satisfying
\begin{align}\label{ee1}
\sup_{0\le t\le T}\|\rho\|_{L^\infty}\le 2\bar{\rho}, \ A_1(T)\le 2 M_*,
\ A_2(T)\le 2\Bbb E_0^\frac23,\ A_3(T)\le 2\Bbb E_0^\frac13,
\end{align}
then the following holds true
\begin{align}\label{ee2}
\sup_{0\le t\le T}\|\rho\|_{L^\infty}\le \frac32\bar{\rho}, \ A_1(T)\le \frac32M_*,
\ A_2(T)\le \Bbb E_0^\frac23,\ A_3(T)\le \Bbb E_0^\frac13,
\end{align}
provided that
\begin{align*}
\Bbb E_0\le \varepsilon_0.
\end{align*}
\end{proposition}

Before proving Proposition \ref{p1}, we establish some necessary \textit{a priori} estimates, see Lemmas \ref{l31}--\ref{l34} below.
\begin{lemma}\label{l31}
Let $(\rho, u, \theta, d)$ be the smooth solution of \eqref{a1}--\eqref{a3} satisfying \eqref{ee1}, then we have
\begin{align}\label{f5}
&\sup_{0\le t\le T}\big(\|\sqrt{\rho}u\|_{L^2}^2+\|\sqrt{\rho}\theta\|_{L^2}^2
+\|\nabla d\|_{L^2}^2\big)+\int_0^T\big(\|\nabla u\|_{L^2}^2+\|\nabla\theta\|_{L^2}^2
+\|\nabla^2d\|_{L^2}^2\big)dt\nonumber\\
&\le C\Big(t_0+\Bbb E_0^\frac{4}{5}\Big)+C\Bbb E_0^\frac{2}{3}\sup_{t_0\le t\le T}\big(\|\nabla u\|_{L^2}^2+\|\nabla\theta\|_{L^2}^2
+\|\nabla^2d\|_{L^2}^2\big),
\end{align}
where $t_0\in (0, T_*)$ is to be determined and $T_*$ is taken from Lemma \ref{l21}.
\end{lemma}
{\it Proof.}
1. Multiplying $\eqref{a1}_2$ by $u$ and integrating by parts over $\Bbb R^3$, we get
\begin{align}\label{f4}
&\frac{1}{2}\frac{d}{dt}\int\rho|u|^2dx+\int\big[\mu|\nabla u|^2
+(\mu+\lambda)(\divv u)^2\big]dx\nonumber\\
&=-\int (u\cdot\nabla)d\cdot\Delta ddx+\int R\rho\theta\divv udx.
\end{align}
Multiplying $\eqref{a1}_4$ by $-\Delta d$ and integrating by parts lead to
\begin{align}\label{zf4}
\frac{1}{2}\frac{d}{dt}\int|\nabla d|^2dx+\int|\nabla^2 d|^2dx=\int(u\cdot\nabla)d\cdot\Delta ddx+\int|\nabla d|^4dx,
\end{align}
where we have used $\Delta d\cdot d=-|\nabla d|^2$ due to $|d|=1$.
Then we derive from \eqref{f4}, \eqref{zf4}, H\"older's inequality, Sobolev's inequality, Gagliardo-Nirenberg inequality, \eqref{lz3.1}, and \eqref{ee1} that
\begin{align*}
&\frac{1}{2}\frac{d}{dt}\int\big(\rho|u|^2+|\nabla d|^2\big)dx
+\int\big[\mu|\nabla u|^2+(\mu+\lambda)(\divv u)^2
+|\nabla^2d|^2\big]dx\nonumber\\
&=\int|\nabla d|^4dx+\int R\rho\theta\divv udx\nonumber\\
&\le (\mu+\lambda)\|\divv u\|_{L^2}^2
+C\int\rho^2\theta^2dx+\int|\nabla d|^4dx\nonumber\\
&\le (\mu+\lambda)\|\divv u\|_{L^2}^2
+C(\bar{\rho})\|\rho\|_{L^\frac32}\|\theta\|_{L^6}^2+\|\nabla d\|_{L^\infty}^2\|\nabla d\|_{L^2}^2\nonumber\\
&\le (\mu+\lambda)\|\divv u\|_{L^2}^2+C\|\rho\|_{L^1}^\frac23\|\rho\|_{L^\infty}^\frac13
\|\nabla\theta\|_{L^2}^2+C\|\nabla^2d\|_{L^2}\|\nabla^3 d\|_{L^2}\|\nabla d\|_{L^2}^2\nonumber\\
&\le (\mu+\lambda)\|\divv u\|_{L^2}^2+c_1\|\nabla\theta\|_{L^2}^2
+C\big(\|\nabla^2d\|_{L^2}^2+\|\nabla^3 d\|_{L^2}^2\big)\|\nabla d\|_{L^2}^2,
\end{align*}
which yields
\begin{align}\label{f6}
& \frac{d}{dt}\big(\|\sqrt{\rho}u\|_{L^2}^2+\|\nabla d\|_{L^2}^2\big)
+2\mu\|\nabla u\|_{L^2}^2
+2\|\nabla^2d\|_{L^2}^2 \notag \\
& \leq 2c_1\|\nabla\theta\|_{L^2}^2
+C\big(\|\nabla^2d\|_{L^2}^2+\|\nabla^3 d\|_{L^2}^2\big)\|\nabla d\|_{L^2}^2.
\end{align}

2. Multiplying $\eqref{a1}_3$ by $\theta$ together with integration by parts and \eqref{a1}$_5$ gives rise to
\begin{align}\label{f7}
&\frac{c_v}{2}\frac{d}{dt}\int\rho\theta^2dx+\kappa\int|\nabla\theta|^2dx\nonumber\\
&=\int\theta\Big(\lambda(\divv u)^2+\frac{\mu}{2}|\nabla u+(\nabla u)^{tr}|^2-R\rho\theta\divv u\Big)dx\nonumber\\
&\quad+\int\big(\Delta d+|\nabla d|^2d\big)\cdot\big(\Delta d+|\nabla d|^2d\big)\theta dx\nonumber\\
&\le \frac{\mu\kappa}{6c_1}\int|\nabla u|^2dx+C\int\theta^2\big(|\nabla u|^2+\rho\theta^2\big)dx+C\int|\nabla d||\nabla^3 d|\theta dx\nonumber\\
&\quad+C\int|\nabla\theta||\nabla d||\nabla^2d|dx+C\int|\nabla d|^4\theta dx+C\int|\nabla d|^2|\nabla^2d|\theta dx\nonumber\\
&\triangleq\frac{\mu\kappa}{6c_1}\|\nabla u\|_{L^2}^2+\sum_{i=1}^5I_i.
\end{align}
By virtue of H\"older's inequality, Sobolev's inequality, and Gagliardo-Nirenberg inequality, we infer from \eqref{ee1} that
\begin{align*}
I_1&\le C\|\theta\|_{L^\infty}^2\big(\|\nabla u\|_{L^2}^2+\|\sqrt{\rho}\theta\|_{L^2}^2\big)\nonumber\\
&\le C\|\theta\|_{L^\infty}^2\big(\|\nabla u\|_{L^2}^2+\|\nabla\theta\|_{L^2}^2\big),\\
I_2+I_3&\le C\|\theta\|_{L^6}\|\nabla^3d\|_{L^2}\|\nabla d\|_{L^3}
+C\|\nabla\theta\|_{L^2}\|\nabla d\|_{L^3}\|\nabla^2 d\|_{L^6}\nonumber\\
&\le C\|\nabla\theta\|_{L^2}\|\nabla^3 d\|_{L^2}\|\nabla d\|_{L^2}^\frac{1}{2}\|\nabla^2d\|_{L^2}^\frac{1}{2}\nonumber\\
&\le \frac{\kappa}{2}\|\nabla\theta\|_{L^2}^2+C\|\nabla^3 d\|_{L^2}^2\|\nabla d\|_{L^2}^2+C\|\nabla^3d\|_{L^2}^2\|\nabla^2d\|_{L^2}^2,\\
I_4 &\le C\|\theta\|_{L^\infty}\|\nabla d\|_{L^2}\|\nabla d\|_{L^6}^3\nonumber\\
&\le C(M_*)\|\theta\|_{L^\infty}\|\nabla d\|_{L^2}\|\nabla^2d\|_{L^2}\nonumber\\
&\le \frac{\kappa}{6c_1}\|\nabla^2d\|_{L^2}^2+C\|\theta\|_{L^\infty}^2\|\nabla d\|_{L^2}^2,\\
I_5&\le C\|\theta\|_{L^\infty}\|\nabla d\|_{L^2}\|\nabla d\|_{L^3}
\|\nabla^2d\|_{L^6}\nonumber\\
&\le C\|\theta\|_{L^\infty}^2\|\nabla d\|_{L^2}^2+C\|\nabla d\|_{L^3}^2\|\nabla^2d\|_{L^6}^2\nonumber\\
&\le C\|\theta\|_{L^\infty}^2\|\nabla d\|_{L^2}^2+C\big(\|\nabla d\|_{L^2}^2+\|\nabla^2d\|_{L^2}^2\big)\|\nabla^3d\|_{L^2}^2.
\end{align*}
Substituting the above inequalities into \eqref{f7}, we arrive at
\begin{align}\label{f8}
c_v\frac{d}{dt}\|\sqrt{\rho}\theta\|_{L^2}^2+\kappa\|\nabla\theta\|_{L^2}^2
&\le\frac{\mu\kappa}{3c_1}\|\nabla u\|_{L^2}^2+\frac{\kappa}{3c_1}\|\nabla^2d\|_{L^2}^2
+C\|\theta\|_{L^\infty}^2\big(\|\nabla u\|_{L^2}^2+\|\nabla\theta\|_{L^2}^2\big)\nonumber\\
&\quad+C\|\theta\|_{L^\infty}^2\|\nabla d\|_{L^2}^2+C\|\nabla^3 d\|_{L^2}^2\|\nabla d\|_{L^2}^2
+C\|\nabla^3d\|_{L^2}^2\|\nabla^2d\|_{L^2}^2,
\end{align}
Adding \eqref{f8} multiplied by $\frac{3c_1}{\kappa}$ to \eqref{f6} implies that
\begin{align}\label{f9}
&\frac{d}{dt}\big(\|\sqrt{\rho}u\|_{L^2}^2+\|\nabla d\|_{L^2}^2
+\|\sqrt{\rho}\theta\|_{L^2}^2\big)+\|\nabla u\|_{L^2}^2+\|\nabla\theta\|_{L^2}^2+\|\nabla^2d\|_{L^2}^2\nonumber\\
&\le C\|\theta\|_{L^\infty}^2\|\nabla d\|_{L^2}^2+C\|\theta\|_{L^\infty}^2\big(\|\nabla u\|_{L^2}^2+\|\nabla\theta\|_{L^2}^2\big)
+C\|\nabla^2d\|_{L^2}^2\|\nabla d\|_{L^2}^2\nonumber\\
&\quad+C\|\nabla^3 d\|_{L^2}^2\|\nabla d\|_{L^2}^2+C\|\nabla^3d\|_{L^2}^2\|\nabla^2d\|_{L^2}^2.
\end{align}

3. Applying the Gagliardo-Nirenberg inequality
\begin{align}\label{zx}
\|\theta\|_{L^\infty}\le C\|\nabla \theta\|_{L^2}^\frac12\|\nabla^2 \theta\|_{L^2}^\frac12,
\end{align}
we derive from \eqref{ee1} that
\begin{align}\label{g12}
\int_0^T\|\theta\|_{L^\infty}^2dt\le C\int_0^T\big(\|\nabla\theta\|_{L^2}^2+\|\nabla^2\theta\|_{L^2}^2\big)dt\le C\Bbb E_0^\frac23.
\end{align}
By H{\"o}lder's inequality, \eqref{1.6}, Sobolev's inequality (see \cite[Theorem]{T1976}), and \eqref{1.9}, we have
\begin{align}\label{3.11}
\|\sqrt{\rho_0}\theta_0\|_{L^2}^2
& \le \|\rho_0\theta_0^\frac45\|_{L^{\frac54}}\|\theta_0^\frac65\|_{L^5} \notag \\
& \le \|\rho_0\|_{L^\infty}^{\frac15}\|\rho_0\theta_0\|_{L^1}^\frac45
\|\theta_0\|_{L^6}^\frac65 \notag \\
& \le \bar{\rho}^{\frac15}c_v^{-\frac45}\|c_v\rho_0\theta_0\|_{L^1}^\frac45
\frac{1}{3^{\frac35}}\Big(\frac2\pi\Big)^{\frac45}
\|\nabla\theta_0\|_{L^2}^\frac65 \notag \\
& \le \bar{\rho}^{\frac15}\frac{1}{3^{\frac35}}\Big(\frac{2}{c_v\pi}\Big)^{\frac45}
\|\nabla\theta_0\|_{L^2}^\frac65\Bbb E_0^\frac45.
\end{align}
Integrating \eqref{f9} with respect to $t$ gives
\begin{align}\label{zx2}
&\sup_{0\le t\le T}\big(\|\sqrt{\rho}u\|_{L^2}^2+\|\sqrt{\rho}\theta\|_{L^2}^2
+\|\nabla d\|_{L^2}^2\big)+\int_0^T\big(\|\nabla u\|_{L^2}^2+\|\nabla\theta\|_{L^2}^2
+\|\nabla^2d\|_{L^2}^2\big)dt\nonumber\\
&\le 2\Bbb E_0^2+\bar{\rho}^{\frac15}\frac{1}{3^{\frac35}}\Big(\frac{2}{c_v\pi}\Big)^{\frac45}
\|\nabla\theta_0\|_{L^2}^\frac65\Bbb E_0^\frac{4}{5} \notag \\
& \quad +\sup_{0\le t\le t_0}\big(\|\nabla u\|_{L^2}^2+\|\nabla\theta\|_{L^2}^2
+\|\nabla^2d\|_{L^2}^2\big)\big(\|\theta\|_{L^\infty}^2+\|\nabla^3d\|_{L^2}^2\big)\nonumber\\
&\quad+C\sup_{t_0\le t\le T}\big(\|\nabla u\|_{L^2}^2+\|\nabla\theta\|_{L^2}^2
+\|\nabla^2d\|_{L^2}^2\big)\int_{t_0}^T\big(\|\theta\|_{L^\infty}^2+\|\nabla^3d\|_{L^2}^2\big)dt\nonumber\\
&\quad+C\sup_{0\le t\le T}\|\nabla d\|_{L^2}^2\int_0^T\big(\|\theta\|_{L^\infty}^2+\|\nabla^3d\|_{L^2}^2
+\|\nabla^2d\|_{L^2}^2\big)dt\nonumber\\
&\le C\Big(t_0+\Bbb E_0^\frac45\Big)+C\Bbb E_0^\frac{2}{3}\sup_{t_0\le t\le T}\big(\|\nabla u\|_{L^2}^2+\|\nabla\theta\|_{L^2}^2
+\|\nabla^2d\|_{L^2}^2\big)\nonumber\\
&\quad+C_1\Bbb E_0^\frac23\sup_{0\le t\le T}\|\nabla d\|_{L^2}^2+C_1\Bbb E_0^\frac13\sup_{0\le t\le T}\|\nabla d\|_{L^2}^2,
\end{align}
provided that
\begin{align*}
\Bbb E_0\leq \frac{\bar\rho}{2}\frac{1}{3^3}\Big(\frac{1}{c_v\pi}\Big)^4
\|\nabla\theta_0\|_{L^2}^6
=\frac{\bar\rho\|\nabla\theta_0\|_{L^2}^6}{54c_v^4\pi^4}.
\end{align*}
Hence, \eqref{f5} follows from \eqref{zx2}
provided that
$$\Bbb E_0\le \varepsilon_1\triangleq\min
\left\{\frac{\bar\rho\|\nabla\theta_0\|_{L^2}^6}{54c_v^4\pi^4}, \frac{1}{(2C_1)^\frac32}\right\}.$$
This finishes the proof of Lemma \ref{l31}.
\hfill $\Box$

\begin{lemma}\label{l32}
Let $(\rho, u, \theta, d)$ be the smooth solution of \eqref{a1}--\eqref{a3} satisfying \eqref{ee1}, then it holds that
\begin{align}\label{f12}
&\frac{d}{dt}\int\Big[\frac{\mu}{2}|\nabla u|^2+\frac{\kappa\delta}{2}|\nabla\theta|^2+|\nabla^2d|^2
+\Big(\frac{2c_2c_3}{\mu}+1\Big)\delta^\frac{3}{2}\rho|\dot{u}|^2+2c_2\delta^\frac32|\nabla d_t|^2\Big]dx\nonumber\\
&\quad+\int\Big(\rho|\dot{u}|^2+\frac14|\nabla d_t|^2+\frac14|\nabla^3d|^2+\frac{c_v\delta}{4}\rho|\dot{\theta}|^2
+\frac{\mu c_2\delta^\frac32}{2}|\nabla \dot{u}|^2+2c_2\delta^\frac32|d_{tt}|^2+c_2\delta^\frac32|\nabla^2d_t|^2\Big)dx\nonumber\\
&\le \frac{d}{dt}\int\Big(R\rho\theta{\rm div}u+\big(\nabla d\odot\nabla d\big):\nabla u
-\frac12\divv u|\nabla d|^2+\lambda\delta\theta(\divv u)^2
+\frac{\mu\delta}{2}\theta|\nabla u+(\nabla u)^{tr}|^2\Big)dx\nonumber\\
&\quad+\frac{d}{dt}\int\big(\delta\theta|\Delta d+|\nabla d|^2d|^2\big)dx
+\delta^{-\frac32}\big(\|\theta\nabla u\|_{L^2}^2+\|\nabla u\|_{L^2}^2+\|\nabla u\|_{L^2}^4\|\nabla\theta\|_{L^2}^2\big)
\nonumber\\
&\quad+6\delta^\frac32\|\nabla u\|_{L^4}^4+\delta^\frac52\|\nabla^3d\|_{L^2}^4
+C\|\nabla^2d\|_{L^2}^4+C\|\nabla u\|_{L^2}^4\|\nabla^2d\|_{L^2}^2
+C\|\nabla d\|_{L^2}\|\nabla^2d\|_{L^2}^3\nonumber\\[3pt]
&\quad+C\delta\|\theta\|_{L^\infty}^2\big(\|\nabla^2d\|_{L^2}^2+\|\nabla d\|_{L^2}\|\nabla^2d\|_{L^2}^3\big),
\end{align}
where $\delta$ is a small positive constant and to be determined.
\end{lemma}
{\it Proof.}
1. Multiplying $\eqref{a1}_2$ by $\dot{u}$ and integrating the resulting equality over $\Bbb R^3$
lead to
\begin{align}\label{f13}
\int\rho|\dot{u}|^2dx&=-\int\dot{u}\cdot\nabla pdx+\mu\int\Delta u\cdot\dot{u}dx+(\mu+\lambda)\int\nabla\divv u\cdot\dot{u}dx\nonumber\\
&\quad+\frac12\int\dot{u}\cdot\nabla|\nabla d|^2dx-\int\dot{u}\cdot\divv(\nabla d\odot\nabla d)dx\triangleq\sum_{i=1}^5J_i.
\end{align}
Since $\dot{\theta}=\theta_t+u\cdot\nabla\theta$,
we have $\rho\dot{\theta}=\rho\theta_t+\rho u\cdot\nabla\theta=(\rho\theta)_t+\divv(\rho u\theta)$ due to \eqref{a1}$_1$. Thus, we get that, for any $\delta>0$,
\begin{align}\label{f14}
J_1&=\int(\divv u)_tpdx-\int(u\cdot\nabla u)\cdot\nabla pdx\nonumber\\
&=\frac{d}{dt}\int R\rho\theta\divv udx-R\int\big[(\rho\theta)_t+\divv(\rho\theta u)\big]\divv udx+\int R\rho\theta\partial_ju^k\partial_ku^jdx\nonumber\\
&\le \frac{d}{dt}\int R\rho\theta{\rm div}udx+\delta^\frac32\big(\|\sqrt{\rho}\dot{\theta}\|_{L^2}^2+\|\theta\nabla u\|_{L^2}^2\big)
+C\delta^{-\frac32}\|\nabla u\|_{L^2}^2.
\end{align}
Here and in what follows, we use the Einstein convention that the repeated indices denote the summation.
Integration by parts implies that
\begin{align}\label{3.20}
J_2&=-\frac{\mu}{2}\frac{d}{dt}\|\nabla u\|_{L^2}^2-\mu\int\partial_iu^j\partial_i(u^k\partial_ku^j)dx\nonumber\\
&\le -\frac{\mu}{2}\frac{d}{dt}\|\nabla u\|_{L^2}^2-\mu\int\Big(\partial_iu^j\partial_iu^k
\partial_ku^j-\frac{1}{2}|\nabla u|^2\divv u\Big)dx\nonumber\\
&\le -\frac{\mu}{2}\frac{d}{dt}\|\nabla u\|_{L^2}^2+C\|\nabla u\|_{L^3}^3.
\end{align}
Similarly to \eqref{3.20}, we have
\begin{align}
J_3&\le -\frac{\mu+\lambda}{2}\frac{d}{dt}\|{\rm div}u\|_{L^2}^2+C\|\nabla u\|_{L^3}^3.
\end{align}
We deduce from $\eqref{a1}_4$, Sobolev's inequality, and Gagliardo-Nirenberg inequality that
\begin{align}
J_4&=-\frac{1}{2}\int\divv u_t|\nabla d|^2dx-\frac{1}{2}\int\divv (u\cdot\nabla u)|\nabla d|^2dx\nonumber\\
&=-\frac12\frac{d}{dt}\int\divv u|\nabla d|^2dx+\int\divv u\nabla d:\nabla\Delta ddx+\int\divv u\nabla d:\nabla(|\nabla d|^2d)dx\nonumber\\
&\quad-\int\divv u\partial_id^j\partial_iu\cdot\nabla d^jdx-\frac12\int\partial_iu\cdot\nabla u^i|\nabla d|^2dx
+\frac{1}{2}\int|\nabla d|^2(\divv u)^2dx\nonumber\\
&\le -\frac12\frac{d}{dt}\int\divv u|\nabla d|^2dx
+\frac{1}{16}\|\nabla\Delta d\|_{L^2}^2+C\||\nabla u|\nabla d\|_{L^2}^2
+C\||\nabla^2d||\nabla d|\|_{L^2}^2+C\|\nabla d\|_{L^6}^6\nonumber\\
&\le -\frac12\frac{d}{dt}\int\divv u|\nabla d|^2dx
+\frac{1}{16}\|\nabla\Delta d\|_{L^2}^2+C\delta^\frac{3}{2}\|\nabla u\|_{L^4}^4+C\|\nabla d\|_{L^4}^4\nonumber\\
&\quad+C\|\nabla d\|_{L^\infty}^2\|\nabla^2d\|_{L^2}^2+C\|\nabla^2d\|_{L^2}^6\nonumber\\
&\le -\frac12\frac{d}{dt}\int\divv u|\nabla d|^2dx
+\frac{1}{16}\|\nabla\Delta d\|_{L^2}^2+C\delta^\frac{3}{2}\|\nabla u\|_{L^4}^4+C\|\nabla d\|_{L^2}\|\nabla^2d\|_{L^2}^3\nonumber\\
&\quad+C\|\nabla^2 d\|_{L^2}\|\nabla^3d\|_{L^2}\|\nabla^2d\|_{L^2}^2+C\|\nabla^2d\|_{L^2}^6\nonumber\\
&\le-\frac12\frac{d}{dt}\int\divv u|\nabla d|^2dx+\frac{1}{16}\|\nabla^3d\|_{L^2}^2
+C\delta^\frac32\|\nabla u\|_{L^4}^4+C\|\nabla d\|_{L^2}\|\nabla^2d\|_{L^2}^3+C\|\nabla^2d\|_{L^2}^6.
\end{align}
Similarly, we arrive at
\begin{align}\label{f18}
J_5\le  \frac{d}{dt}\int(\nabla d\odot\nabla d):\nabla udx+\frac{1}{16}\|\nabla^3d\|_{L^2}^2
+C\delta^\frac{3}{2}\|\nabla u\|_{L^4}^4+C\|\nabla d\|_{L^2}^4\|\nabla^2d\|_{L^2}^2+C\|\nabla^2d\|_{L^2}^6.
\end{align}
Inserting \eqref{f14}--\eqref{f18} into \eqref{f13} yields that
\begin{align}\label{f19}
&\frac{1}{2}\frac{d}{dt}\int\big[\mu|\nabla u|^2+(\mu+\lambda)
(\divv u)^2\big]dx+\int\rho|\dot{u}|^2dx\nonumber\\
&\le \frac{d}{dt}\int\Big(R\rho\theta\divv u+(\nabla d\odot\nabla d):\nabla u
-\frac12\divv u|\nabla d|^2\Big)dx+\frac18\|\nabla^3d\|_{L^2}^2\nonumber\\
&\quad+\delta^\frac32C\big(\|\sqrt{\rho}\dot{\theta}\|_{L^2}^2+\|\theta\nabla u\|_{L^2}^2+\|\nabla u\|_{L^4}^4\big)
+C\delta^{-\frac32}\|\nabla u\|_{L^2}^2+C\|\nabla^2d\|_{L^2}^6
+C\|\nabla d\|_{L^2}\|\nabla^2d\|_{L^2}^3.
\end{align}

2. It follows from $\eqref{a1}_4$, H\"older's and Gagliardo-Nirenberg inequalities that
\begin{align}\label{f20}
&\frac{d}{dt}\int|\nabla^2 d|^2dx+\int\big(|\nabla d_t|^2+|\nabla^3d|^2\big)dx\nonumber\\
&=\int|\nabla d_t-\nabla\Delta d|^2dx\nonumber\\
&=\int|-\nabla(u\cdot\nabla d)+\nabla\big(|\nabla d|^2d\big)|^2dx\nonumber\\
&\le C\||\nabla u||\nabla d|\|_{L^2}^2+C\||u||\nabla^2d|\|_{L^2}^2+C\|\nabla d\|_{L^6}^6+C\||\nabla d||\nabla^2d|\|_{L^2}^2\nonumber\\
&\le C\|\nabla u\|_{L^2}^2\|\nabla d\|_{L^\infty}^2+C\|u\|_{L^6}^2\|\nabla^2d\|_{L^3}^2
+C\|\nabla^2d\|_{L^2}^6+C\|\nabla d\|_{L^\infty}^2\|\nabla^2d\|_{L^2}^2\nonumber\\
&\le C\|\nabla u\|_{L^2}^2\|\nabla^2 d\|_{L^2}\|\nabla^3d\|_{L^2}
+C\|\nabla^2d\|_{L^2}^3\|\nabla^3d\|_{L^2}+C\|\nabla^2d\|_{L^2}^2\nonumber\\
&\le \frac{1}{8}\|\nabla^3d\|_{L^2}^2
+C\|\nabla^2d\|_{L^2}^6+C\|\nabla u\|_{L^2}^4\|\nabla^2d\|_{L^2}^2.
\end{align}

3. Multiplying $\eqref{a1}_3$ by $\dot{\theta}$ and integrating the resulting equation by parts yields that
\begin{align}\label{f21}
\frac{\kappa}{2}\frac{d}{dt}\int|\nabla\theta|^2dx
+c_v\int\rho|\dot{\theta}|^2dx
&=\kappa\int u\cdot\nabla\theta\Delta\theta dx
-\int\rho\theta\divv u\dot{\theta}dx+\lambda\int(\divv u)^2\dot{\theta}dx\nonumber\\
&\quad+\frac{\mu}{2}\int|\nabla u+(\nabla u)^{tr}|^2\dot{\theta}dx+\int|\Delta d+|\nabla d|^2d|^2\dot{\theta}dx\triangleq\sum_{i=1}^5N_i.
\end{align}
Let us deal with the right-hand side terms in \eqref{f21}.  In view of the standard estimate of elliptic system, one obtains
\begin{align}\label{f22}
\|\nabla^2\theta\|_{L^2}&\le C\big(\|\rho\dot{\theta}\|_{L^2}+\|\rho\theta{\rm div}u\|_{L^2}+\|\nabla u\|_{L^4}^2
+\|\Delta d+|\nabla d|^2d\|_{L^4}^2\big)\nonumber\\
&\le C\|\sqrt{\rho}\dot{\theta}\|_{L^2}+C\|\theta\nabla u\|_{L^2}+\|\nabla u\|_{L^4}^2
+C\|\Delta d+|\nabla d|^2d\|_{L^4}^2.
\end{align}
Thus, we obtain from Sobolev's inequality, Gagliardo-Nirenberg inequality, and \eqref{f22} that
\begin{align}
N_1&=\kappa\int\Big(\frac12|\nabla\theta|^2{\rm div}u-\partial_i\theta\partial_iu^k\partial_k\theta\Big)dx\nonumber\\
&\le C\|\nabla u\|_{L^2}\|\nabla\theta\|_{L^2}^\frac{1}{2}\|\nabla^2\theta\|_{L^2}^\frac{3}{2}\nonumber\\
&\le \delta^\frac12\|\nabla^2\theta\|_{L^2}^2
+C\delta^{-\frac32}\|\nabla u\|_{L^2}^4\|\nabla\theta\|_{L^2}^2\nonumber\\
&\le \delta^\frac12\big(\|\sqrt{\rho}\dot{\theta}\|_{L^2}^2+C\|\theta\nabla u\|_{L^2}^2+\|\nabla u\|_{L^4}^4
+\|\Delta d+|\nabla d|^2d\|_{L^4}^4\big)+C\delta^{-\frac32}\|\nabla u\|_{L^2}^4\|\nabla\theta\|_{L^2}^2\nonumber\\
&\le \delta^\frac{1}{2}\big(\|\sqrt{\rho}\dot{\theta}\|_{L^2}^2+C\|\theta\nabla u\|_{L^2}^2+\|\nabla u\|_{L^4}^4\big)+
C\delta^{-\frac32}\|\nabla u\|_{L^2}^4\|\nabla\theta\|_{L^2}^2\nonumber\\
&\quad+\delta^\frac12\|\nabla^2d\|_{L^2}\|\nabla^3d\|_{L^2}^3+\delta^\frac12
\|\nabla d\|_{L^\infty}^2\|\nabla d\|_{L^6}^6\nonumber\\
&\le \delta^\frac12\big(\|\sqrt{\rho}\dot{\theta}\|_{L^2}^2+C\|\theta\nabla u\|_{L^2}^2+\|\nabla u\|_{L^4}^4\big)+
C\delta^{-\frac32}\|\nabla u\|_{L^2}^4\|\nabla\theta\|_{L^2}^2\nonumber\\
&\quad+\delta^\frac32\|\nabla^3d\|_{L^2}^4+\delta^\frac32\|\nabla^3d\|_{L^2}^2
+C\|\nabla^2d\|_{L^2}^4+C\|\nabla^2d\|_{L^2}^{14}.
\end{align}
It is easy to see that
\begin{align}
N_2\le \delta^\frac12\|\sqrt{\rho}\dot{\theta}\|_{L^2}^2+C\delta^{-\frac12}\|\theta\nabla u\|_{L^2}^2.
\end{align}
Integration by parts together with \eqref{f3} leads to
\begin{align}\label{f25}
N_3&=\lambda\int(\divv u)^2\theta_tdx
+\lambda\int(\divv u)^2(u\cdot\nabla\theta)dx\nonumber\\
&=\lambda\frac{d}{dt}\int(\divv u)^2\theta dx-2\lambda\int\theta\divv u
\divv \dot{u}dx
+2\lambda\int\theta\divv u\divv (u\cdot\nabla u)dx
+\lambda\int(\divv u)^2(u\cdot\nabla\theta)dx\nonumber\\
&=\lambda\frac{d}{dt}\int(\divv u)^2\theta dx-2\lambda\int\theta\divv u
\divv \dot{u}dx
+2\lambda\int\theta\divv u\partial_iu^j\partial_ju^idx
+\lambda\int u\cdot\nabla(\theta(\divv u)^2)dx\nonumber\\
&=\lambda\frac{d}{dt}\int(\divv u)^2\theta dx-2\lambda\int\theta\divv u\divv \dot{u}dx
+2\lambda\int\theta\divv u\partial_iu^j\partial_ju^idx
-\lambda\int\theta(\divv u)^3dx\nonumber\\
&\le \lambda\frac{d}{dt}\int(\divv u)^2\theta dx
+\delta\big(\|\nabla u\|_{L^4}^4+\|\nabla\dot{u}\|_{L^2}^2\big)+C\delta^{-1}\|\theta\nabla u\|_{L^2}^2.
\end{align}
Similarly to \eqref{f25}, one has
\begin{align}
N_4&\le \frac{\mu}{2}\frac{d}{dt}\int\theta|\nabla u+(\nabla u)^{tr}|^2
+\delta\big(\|\nabla u\|_{L^4}^4+\|\nabla\dot{u}\|_{L^2}^2\big)+C\delta^{-1}\|\theta\nabla u\|_{L^2}^2.
\end{align}
For the term $N_5$, we deduce from Gagliardo-Nirenberg, H\"older's inequality, and Young's inequalities that
\begin{align}
N_5&=\int|\Delta d+|\nabla d|^2d|^2\theta_tdx+\int|\Delta d+|\nabla d|^2d|^2(u\cdot\nabla\theta)dx\nonumber\\
&=\frac{d}{dt}\int|\Delta d+|\nabla d|^2d|^2\theta dx-\int\big(|\Delta d+|\nabla d|^2d|^2\big)_t\theta dx
+\int|\Delta d+|\nabla d|^2d|^2(u\cdot\nabla\theta)dx\nonumber\\
&\le \frac{d}{dt}\int|\Delta d+|\nabla d|^2d|^2\theta dx
+C\|\theta\|_{L^\infty}\big(\|\nabla^2 d\|_{L^2}+\|\nabla d\|_{L^4}^2\big)\nonumber\\[3pt]
&\quad\times\big(\|\nabla^2d_t\|_{L^2}+\|d_t\|_{L^6}\|\nabla d\|_{L^2}\|\nabla d\|_{L^3}+\|\nabla d\|_{L^3}\|\nabla d_t\|_{L^6}\big)\nonumber\\
&\quad+C\big(\|\nabla^2d\|_{L^6}^2+\||\nabla d|^2\|_{L^6}^2\big)\|u\|_{L^6}\|\nabla\theta\|_{L^2}\nonumber\\
&\le \frac{d}{dt}\int|\Delta d+|\nabla d|^2d|^2\theta dx
+C\|\theta\|_{L^\infty}^2\big(\|\nabla^2d\|_{L^2}^2+\|\nabla d\|_{L^4}^4\big)
+\delta^\frac12\|\nabla^2 d_t\|_{L^2}^2\nonumber\\
&\quad+\delta^\frac12\|\nabla d_t\|_{L^2}^2\|\nabla d\|_{L^2}^3\|\nabla^2d\|_{L^2}
+\delta^\frac12C\big(\|\nabla d\|_{L^2}^2+\|\nabla^2d\|_{L^2}^2\big)\|\nabla^2d_t\|_{L^2}^2\nonumber\\
&\quad+C\|\nabla u\|_{L^2}\|\nabla\theta\|_{L^2}\big(\|\nabla^3d\|_{L^2}^2+C\|\nabla d\|_{L^\infty}^2\|\nabla^2d\|_{L^2}^2\big)\nonumber\\
&\le \frac{d}{dt}\int|\Delta d+|\nabla d|^2d|^2\theta dx
+C\|\theta\|_{L^\infty}^2\big(\|\nabla^2d\|_{L^2}^2+\|\nabla d\|_{L^2}\|\nabla^2d\|_{L^2}^3\big)
+\delta^\frac12\|\nabla^2 d_t\|_{L^2}^2\nonumber\\
&\quad+\delta^\frac12\|\nabla d_t\|_{L^2}^2
+\delta^\frac12C\big(\|\nabla d\|_{L^2}^2+\|\nabla^2d\|_{L^2}^2\big)\|\nabla^2d_t\|_{L^2}^2+C\|\nabla^3d\|_{L^2}^2
+C\|\nabla^2 d\|_{L^2}^3\|\nabla^3 d\|_{L^2}\nonumber\\
&\le \frac{d}{dt}\int|\Delta d+|\nabla d|^2d|^2\theta dx
+C\|\theta\|_{L^\infty}^2\big(\|\nabla^2d\|_{L^2}^2+\|\nabla d\|_{L^2}\|\nabla^2d\|_{L^2}^3\big)
\nonumber\\
&\quad+\delta^\frac12\|\nabla d_t\|_{L^2}^2
+\delta^\frac12C\|\nabla^2d_t\|_{L^2}^2+C\|\nabla^3d\|_{L^2}^2+C\|\nabla^2 d\|_{L^2}^6.
\end{align}
Inserting the above inequalities back into \eqref{f21}, we get
\begin{align}\label{f28}
&\frac{\kappa}{2}\frac{d}{dt}\int|\nabla\theta|^2dx
+\big(c_v-2\delta^\frac12\big)\int\rho|\dot{\theta}|^2dx\nonumber\\
&\le \frac{d}{dt}\int\theta\Big(\lambda(\divv u)^2+\frac{\mu}{2}|\nabla u+(\nabla u)^{tr}|^2+
|\Delta d+|\nabla d|^2d|^2\Big)dx\nonumber\\
&\quad+\delta^{-\frac32}\big(\|\theta\nabla u\|_{L^2}^2+\|\nabla u\|_{L^2}^4\|\nabla\theta\|_{L^2}^2\big)
+3\delta^\frac12\|\nabla u\|_{L^4}^4+2\delta\|\nabla\dot{u}\|_{L^2}^2\nonumber\\[3pt]
&\quad+C_2\|\nabla^3d\|_{L^2}^2+\delta^\frac32\|\nabla^3d\|_{L^2}^4
+C\|\nabla^2d\|_{L^2}^4
+\delta^\frac{1}{2}\|\nabla d_t\|_{L^2}^2+\delta^\frac12C\|\nabla^2d_t\|_{L^2}^2\nonumber\\[3pt]
&\quad
+C\|\theta\|_{L^\infty}^2\big(\|\nabla^2d\|_{L^2}^2+\|\nabla d\|_{L^2}\|\nabla^2d\|_{L^2}^3\big).
\end{align}
Thus, combining \eqref{f19}, \eqref{f20}, and \eqref{f28} multiplied by $\delta$ altogether, we derive that, for all $0<\delta\le \min\Big\{1, \frac{1}{\sqrt[3]{4}},
\frac{1}{4C_2}, \frac{c_v^2}{16}\Big\}$,
\begin{align}\label{f29}
&\frac{d}{dt}\int\Big(\frac{\mu}{2}|\nabla u|^2+\frac{\kappa\delta}{2}|\nabla\theta|^2+|\nabla^2d|^2\Big)dx
+\int\Big(\rho|\dot{u}|^2+\frac12|\nabla d_t|^2+\frac12|\nabla^3d|^2+\frac{c_v\delta}{2}\rho|\dot{\theta}|^2
\Big)dx\nonumber\\
&\le \frac{d}{dt}\int\Big(R\rho\theta\divv u+(\nabla d\odot\nabla d):\nabla u
-\frac{1}{2}\divv u|\nabla d|^2+\lambda\theta(\divv u)^2+\frac{\mu}{2}\theta|\nabla u+(\nabla u)^{tr}|^2\Big)dx\nonumber\\
&\quad+\frac{d}{dt}\int\big(\theta|\Delta d+|\nabla d|^2d|^2\big)dx
+\delta^{-\frac32}\big(\|\theta\nabla u\|_{L^2}^2+\|\nabla u\|_{L^2}^2+\|\nabla u\|_{L^2}^4\|\nabla\theta\|_{L^2}^2\big)
+2\delta^2\|\nabla\dot{u}\|_{L^2}^2\nonumber\\
&\quad+6\delta^\frac32\|\nabla u\|_{L^4}^4+C\|\nabla^2d\|_{L^2}^4
+C\|\nabla u\|_{L^2}^4\|\nabla^2d\|_{L^2}^2+C\|\nabla d\|_{L^2}\|\nabla^2d\|_{L^2}^3+\delta^\frac{3}{2}c_2\|\nabla^2d_t\|_{L^2}^2\nonumber\\[3pt]
&\quad+\delta^\frac52\|\nabla^3d\|_{L^2}^4
+C\delta\|\theta\|_{L^\infty}^2\big(\|\nabla^2d\|_{L^2}^2+\|\nabla d\|_{L^2}\|\nabla^2d\|_{L^2}^3\big).
\end{align}

4. Operating $\partial_t+\divv(u\cdot)$ to the $j$-th component of $\eqref{a1}_2$ and
multiplying the resulting equation by $\dot{u}^j$, one gets by some calculations that
\begin{align}\label{f30}
&\frac{1}{2}\frac{d}{dt}\int\rho|\dot{u}|^2dx\nonumber\\
&=\mu\int\dot{u}^j(\partial_t\Delta u^j+{\rm div}(u\Delta u^j))dx
+(\mu+\lambda)\int\dot{u}^j(\partial_t\partial_j(\divv u)+{\rm div}(u\partial_j({\rm div}u)))dx\nonumber\\
&\quad-\int\dot{u}^j(\partial_jp_t+{\rm div}(u\partial_jp))dx
-\int\dot{u}^j\big(\partial_t(\nabla d\cdot\Delta d_j)+{\rm div}(u\cdot\nabla d\cdot\Delta d_j)\big)dx\triangleq\sum_{i=1}^4K_i.
\end{align}
Integration by parts leads to
\begin{align}\label{3.36}
K_1&=-\mu\int\big(\partial_i\dot{u}^j\partial_t\partial_iu^j+\Delta u^ju\cdot\nabla\dot{u}^j\big)dx\nonumber\\
&=-\mu\int\big(|\nabla\dot{u}|^2-\partial_i\dot{u}^ju^k\partial_k\partial_iu^j-\partial_i\dot{u}^j\partial_iu^k\partial_ku^j
+\Delta u^ju\cdot\nabla\dot{u}^j\big)dx\nonumber\\
&=-\mu\int\big(|\nabla\dot{u}|^2+\partial_i\dot{u}^j\partial_ku^k\partial_iu^j-\partial_i\dot{u}^j\partial_iu^k\partial_ku^j
-\partial_iu^j\partial_iu^k\partial_k\dot{u}^j\big)dx\nonumber\\
&\le -\frac{3\mu}{4}\|\nabla\dot{u}\|_{L^2}^2+C\|\nabla u\|_{L^4}^4.
\end{align}
Similarly, one has
\begin{align}
K_2&\le -\frac{\mu+\lambda}{2}\|{\rm div}\dot{u}\|_{L^2}^2+C\|\nabla u\|_{L^4}^4.
\end{align}
It follows from integration by parts, $\eqref{a1}_1$, and Young's inequality that
\begin{align}
K_3&=\int\big(\partial_j\dot{u}^jp_t+\partial_j pu\cdot\nabla\dot{u}^j\big)dx\nonumber\\
&=\int\partial_j\dot{u}^jp_tdx-\int p\partial_j\big(u\cdot\nabla\dot{u}^j\big)dx\nonumber\\
&=\int\partial_j\dot{u}^j\big[(\rho\theta)_t+\divv(\rho u\theta)\big]dx
-\int\rho\theta\partial_ju\cdot\nabla\dot{u}^jdx\nonumber\\
&=\int\rho\dot{\theta}\partial_j\dot{u}^jdx-\int\rho\theta\partial_ju\cdot\nabla\dot{u}^jdx\nonumber\\
&\le \frac{\mu}{4}\|\nabla\dot{u}\|_{L^2}^2+C\|\theta\nabla u\|_{L^2}^2+C\|\sqrt{\rho}\dot{\theta}\|_{L^2}^2.
\end{align}
Integrating by parts and applying H\"older's inequality and Gagliardo-Nirenberg
inequality, we arrive at
\begin{align}\label{3.39}
K_4&\le C\int|\nabla\dot{u}|\big(|\nabla d||\nabla d_t|+|u||\nabla d||\nabla^2d|\big)dx\nonumber\\
&\le \frac{\mu}{8}\|\nabla\dot{u}\|_{L^2}^2+C\|\nabla d\|_{L^6}^2\|\nabla d_t\|_{L^3}^2
+C\|u\|_{L^6}^2\|\nabla d\|_{L^6}^2\|\nabla^2d\|_{L^6}^2\nonumber\\
&\le \frac{\mu}{8}\|\nabla\dot{u}\|_{L^2}^2+C\|\nabla d\|_{L^6}^2\|\nabla d_t\|_{L^2}\|\nabla^2d_t\|_{L^2}+C\|\nabla u\|_{L^2}^2\|\nabla^2d\|_{L^2}^2\|\nabla^3d\|_{L^2}^2\nonumber\\
&\le \frac{\mu}{8}\|\nabla\dot{u}\|_{L^2}^2+\delta\|\nabla^2d_t\|_{L^2}^2
+C\|\nabla^2d\|_{L^2}^4\|\nabla d_t\|_{L^2}^2+C\|\nabla^3d\|_{L^2}^2.
\end{align}
Substituting \eqref{3.36}--\eqref{3.39} into \eqref{f30}, we have
\begin{align}\label{f35}
&\frac{1}{2}\frac{d}{dt}\int\rho|\dot{u}|^2dx+\frac{\mu}{2}\int|\nabla\dot{u}|^2dx\nonumber\\
&\le \delta\|\nabla^2 d_t\|_{L^2}^2+C\|\theta\nabla u\|_{L^2}^2+C\|\sqrt{\rho}\dot{\theta}\|_{L^2}^2
+C\|\nabla^2d\|_{L^2}^4\|\nabla d_t\|_{L^2}^2+C\|\nabla^3d\|_{L^2}^2.
\end{align}

5. Differentiation $\eqref{a1}_4$ with respect to $t$ leads to
\begin{align}
d_{tt}-\Delta d_t=(-u\cdot\nabla d+|\nabla d|^2d)_t
=-u_t\cdot\nabla d-u\cdot\nabla d_t+|\nabla d|^2d_t+(|\nabla d|^2)_td,
\end{align}
which combined with \eqref{a1}$_5$, H{\"o}lder's inequality, Sobolev's inequality, Gagliardo-Nirenberg inequality, and \eqref{ee1} implies that
\begin{align}\label{rrf}
&\frac{d}{dt}\int|\nabla d_t|^2dx+\int\big(|d_{tt}|^2+|\nabla^2d_t|^2\big)dx\nonumber\\
&=\int|d_{tt}-\Delta d_t|^2dx\nonumber\\
&\le C\int|u_t|^2|\nabla d|^2dx+C\int|u|^2|\nabla d_t|^2dx+C\int|\nabla d|^4|d_t|^2dx+C\int|\nabla d|^2|\nabla d_t|^2dx\nonumber\\
&\le C\int|\dot{u}|^2|\nabla d|^2dx+C\int|u|^2|\nabla u|^2|\nabla d|^2dx
+C\int|u|^2|\nabla d_t|^2dx\nonumber\\
&\quad+C\int|\nabla d|^4|d_t|^2dx+C\int|\nabla d|^2|\nabla d_t|^2dx\nonumber\\
&\le C\|\dot{u}\|_{L^6}^2\|\nabla d\|_{L^3}^2+C\|u\|_{L^6}^2\|\nabla u\|_{L^4}^2\|\nabla d\|_{L^{12}}^2
+C\|u\|_{L^6}^2\|\nabla d_t\|_{L^2}\|\nabla d_t\|_{L^6}\nonumber\\
&\quad+C\|\nabla d\|_{L^6}^4\|d_t\|_{L^6}^2+C\|\nabla d\|_{L^6}^2\|\nabla d_t\|_{L^2}\|\nabla d_t\|_{L^6}\nonumber\\
&\le C\|\nabla\dot{u}\|_{L^2}^2\|\nabla d\|_{L^2}\|\nabla^2d\|_{L^2}
+C\|\nabla u\|_{L^2}^2\|\nabla u\|_{L^4}^2\|\nabla d\|_{L^3}\|\nabla^3d\|_{L^2}\nonumber\\
&\quad+C\|\nabla u\|_{L^2}^2\|\nabla d_t\|_{L^2}\|\nabla^2d_t\|_{L^2}+C\|\nabla^2d\|_{L^2}^4\|\nabla d_t\|_{L^2}^2
+C\|\nabla^2d\|_{L^2}^2\|\nabla d_t\|_{L^2}\|\nabla^2d_t\|_{L^2}\nonumber\\
&\le \frac12\|\nabla^2d_t\|_{L^2}^2+C\|\nabla\dot{u}\|_{L^2}^2+C\|\nabla u\|_{L^4}^4+C\|\nabla^3d\|_{L^2}^2
+C\big(\|\nabla u\|_{L^2}^4+\|\nabla^2d\|_{L^2}^4\big)\|\nabla d_t\|_{L^2}^2.
\end{align}
This indicates that
\begin{align}\label{f37}
&\frac{d}{dt}\int|\nabla d_t|^2dx+\int\big(|d_{tt}|^2+|\nabla^2d_t|^2\big)dx\nonumber\\
&\le c_3\|\nabla\dot{u}\|_{L^2}^2+C\|\nabla u\|_{L^4}^4+C\|\nabla^3d\|_{L^2}^2
+C\big(\|\nabla u\|_{L^2}^4+\|\nabla^2d\|_{L^2}^4\big)\|\nabla d_t\|_{L^2}^2.
\end{align}
Adding \eqref{f35} multiplying by $(\frac{2c_3}{\mu}+1)$ to \eqref{f37}, one has
\begin{align}\label{f38}
&\frac{d}{dt}\int\Big[\Big(\frac{c_3}{\mu}+\frac12\Big)\rho|\dot{u}|^2+|\nabla d_t|^2\Big]dx
+\int\Big(\frac{\mu}{2}|\nabla \dot{u}|^2+|d_{tt}|^2+|\nabla^2d_t|^2\Big)dx\nonumber\\
&\le C\|\nabla u\|_{L^4}^4+C\|\nabla^3d\|_{L^2}^2
+C\big(\|\nabla u\|_{L^2}^4+\|\nabla^2d\|_{L^2}^4\big)\|\nabla d_t\|_{L^2}^2
+C\|\theta\nabla u\|_{L^2}^2+C\|\sqrt{\rho}\dot{\theta}\|_{L^2}^2\nonumber\\
&\le c_4\big(\|\nabla u\|_{L^4}^4+\|\nabla^3d\|_{L^2}^2+\|\nabla d_t\|_{L^2}^2+\|\sqrt{\rho}\dot{\theta}\|_{L^2}^2
+\|\theta\nabla u\|_{L^2}^2\big).
\end{align}
Thus, combining \eqref{f29} and \eqref{f38} multiplied by $2c_2\delta^\frac32$ together, we conclude the desired \eqref{f12} as long as $\delta$ is chosen suitably small such that
\begin{align*}
0<\delta\le\min\left\{\frac{1}{\sqrt[3]{4}}, \frac{1}{2C_2}, \frac{c_v^2}{16}, \frac{c_v^2}{(8c_2c_4)^2}, \frac{1}{(8c_2c_4)^\frac23}, \frac{\mu}{2(2c_2+\mu)}\right\}.
\end{align*}
The proof of Lemma \ref{l32} is complete.
\hfill $\Box$

Next we define
\begin{align}\label{f39}
\sigma(t)\triangleq\min\{t_0, t\},
\end{align}
where $t_0\in (0, T_*)$ is the same as that in Lemma \ref{l31}.
\begin{lemma}\label{l33}
Let $\sigma(t)$ be as in \eqref{f39} and let $(\rho, u, \theta, d)$ be the smooth solution of \eqref{a1}--\eqref{a3} satisfying \eqref{ee1}, then it holds that
\begin{align}\label{f42}
&\sup_{0\le t\le T}\big[\sigma\big(\|\nabla u\|_{L^2}^2
+\|\nabla\theta\|_{L^2}^2+\|\nabla^2d\|_{L^2}^2+\|\sqrt{\rho}\dot{u}\|_{L^2}^2
+\|\nabla d_t\|_{L^2}^2\big)\big]\nonumber\\
&\quad+\int_0^T\sigma\big(\|\sqrt{\rho}\dot{u}\|_{L^2}^2
+\|\sqrt{\rho}\dot{\theta}\|_{L^2}^2
+\|\nabla\dot{u}\|_{L^2}^2+\|d_{tt}\|_{L^2}^2+\|\nabla d_t\|_{L^2}^2+\|\nabla^3d\|_{L^2}^2+\|\nabla^2d_t\|_{L^2}^2\big)dt\nonumber\\
&\le C\Big(t_0+\Bbb E_0^\frac{4}{5}\Big).
\end{align}
\end{lemma}
{\it Proof.}
1. Multiplying \eqref{f12} by $\sigma(t)$ and integrating the resultant with respect to $t$, we have
\begin{align}\label{f43}
&\sup_{0\le t\le T}\big[\sigma\big(\|\nabla u\|_{L^2}^2
+\delta\|\nabla\theta\|_{L^2}^2+\|\nabla^2d\|_{L^2}^2+\delta^\frac{3}{2}\|\sqrt{\rho}\dot{u}\|_{L^2}^2
+\delta^\frac{3}{2}\|\nabla d_t\|_{L^2}^2\big)\big]\nonumber\\
&\quad+\int_0^T\sigma\big(\|\sqrt{\rho}\dot{u}\|_{L^2}^2
+\|\nabla d_t\|_{L^2}^2+\|\nabla^3d\|_{L^2}^2
+\delta\|\sqrt{\rho}\dot{\theta}\|_{L^2}^2
+\delta^\frac{3}{2}\|\nabla\dot{u}\|_{L^2}^2\big)dt\nonumber\\
&\quad+\delta^\frac{3}{2}\int_0^T\sigma\big(\|d_{tt}\|_{L^2}^2
+\|\nabla^2d_t\|_{L^2}^2\big)dt\nonumber\\
&\le C\int_0^{t_0}\big(\|\nabla u\|_{L^2}^2
+\delta\|\nabla\theta\|_{L^2}^2+\|\nabla^2d\|_{L^2}^2+\delta^\frac{3}{2}\|\sqrt{\rho}\dot{u}\|_{L^2}^2
+\delta^\frac{3}{2}\|\nabla d_t\|_{L^2}^2\big)dt\nonumber\\
&\quad+C\int_0^{t_0}\int\big(\rho\theta|\divv u|+\theta|\nabla u|^2+|\nabla u||\nabla d|^2+\theta|\Delta d+|\nabla d|^2d|^2\big)dxdt\nonumber\\
&\quad+C\sigma\int\big(\rho\theta|\nabla u|+|\nabla u||\nabla d|^2+\delta\theta|\nabla u|^2+\delta\theta|\Delta d+|\nabla d|^2d|^2\big)dx\nonumber\\
&\quad+C\int_0^T\sigma\big(\|\theta\nabla u\|_{L^2}^2+\|\nabla u\|_{L^2}^2
+\|\nabla u\|_{L^2}^4\|\nabla\theta\|_{L^2}^2\big)dt
+C\delta^\frac52\int_0^T\sigma\|\nabla^3d\|_{L^2}^4dt\nonumber\\
&\quad+C\delta\int_0^T\sigma\|\theta\|_{L^\infty}^2\big(\|\nabla^2d\|_{L^2}^2+\|\nabla d\|_{L^2}\|\nabla^2d\|_{L^2}^3\big)dt
+C\delta^\frac32\int_0^T\sigma\|\nabla u\|_{L^4}^4dt\nonumber\\
&\quad+C\int_0^T\sigma\big(\|\nabla^2d\|_{L^2}^2+\|\nabla u\|_{L^2}^4
+\|\nabla d\|_{L^2}^2\big)\|\nabla^2d\|_{L^2}^2dt.
\end{align}
It follows from \eqref{w1} that
\begin{align}\label{f44}
&C\int_0^{t_0}\big(\|\nabla u\|_{L^2}^2
+\delta\|\nabla\theta\|_{L^2}^2+\|\nabla^2d\|_{L^2}^2+\delta^\frac{3}{2}\|\sqrt{\rho}\dot{u}\|_{L^2}^2
+\delta^\frac{3}{2}\|\nabla d_t\|_{L^2}^2\big)dt\nonumber\\
&\le Ct_0\sup_{0\le t\le t_0}\big(\|\nabla u\|_{L^2}^2
+\delta\|\nabla\theta\|_{L^2}^2+\|\nabla^2d\|_{L^2}^2+\delta^\frac{3}{2}\|\sqrt{\rho}\dot{u}\|_{L^2}^2
+\delta^\frac{3}{2}\|\nabla d_t\|_{L^2}^2\big)\nonumber\\
&\le C(M_*)t_0.
\end{align}
Applying \eqref{w1}, \eqref{ee1}, H\"older's inequality, and Gagliardo-Nirenberg inequality, we have
\begin{align}\label{f45}
&C\int_0^{t_0}\int\big(\rho\theta|{\rm div}u|+\theta|\nabla u|^2+|\nabla u||\nabla d|^2+\theta|\Delta d+|\nabla d|^2d|^2\big)dxdt\nonumber\\
&\le C(\bar{\rho})\int_0^{t_0}\|\nabla u\|_{L^2}\big(\|\sqrt{\rho}\theta\|_{L^2}
+\|\nabla d\|_{L^4}^2\big)dt+C\int_0^{t_0}\|\theta\|_{L^\infty}\big(\|\nabla u\|_{L^2}^2+\|\Delta d+|\nabla d|^2d\|_{L^2}^2\big)dt\nonumber\\
&\le C\int_0^{t_0}\|\nabla u\|_{L^2}\Big(\|\sqrt{\rho}\theta\|_{L^2}
+\|\nabla d\|_{L^2}^\frac12\|\nabla^2d\|_{L^2}^\frac32\Big)dt\nonumber\\
&\quad+C\int_0^{t_0}\|\nabla \theta\|_{L^2}^\frac12\|\nabla^2\theta\|_{L^2}^\frac12\big(
\|\nabla u\|_{L^2}^2+\|\nabla^2 d\|_{L^2}^2
+\|\nabla d\|_{L^4}^4\big)dt\nonumber\\
&\le C\int_0^{t_0}\big(\|\nabla u\|_{L^2}^2+\|\nabla\theta\|_{L^2}^2+\|\nabla^2d\|_{L^2}^2
+\|\nabla d\|_{L^2}^2\|\nabla u\|_{L^2}^4\big)dt\nonumber\\
&\quad+\int_0^{t_0}\big(\|\nabla\theta\|_{L^2}^2+\|\nabla^2\theta\|_{L^2}^2
+\|\nabla u\|_{L^2}^4+\|\nabla^2d\|_{L^2}^4+\|\nabla d\|_{L^2}\|\nabla^2d\|_{L^2}^3\big)dt\nonumber\\
&\le C(\bar{\rho}, M_*)t_0.
\end{align}
By Cauchy-Schwarz inequality, \eqref{lz3.3}, and \eqref{ee1}, we get that
\begin{align}\label{f46}
C\sigma\int|\nabla u||\nabla d|^2dx&\le C\sigma\|\nabla u\|_{L^2}\|\nabla d\|_{L^4}^2\nonumber\\
&\le \frac14\sigma\|\nabla u\|_{L^2}^2+C\sigma\|\nabla d\|_{L^2}\|\nabla^2d\|_{L^2}^3\nonumber\\
&\le \frac14\sigma\|\nabla u\|_{L^2}^2+\delta^{-1}\sigma\|\nabla^2d\|_{L^2}^4
+C\delta\sigma\|\nabla d\|_{L^2}^2\|\nabla^2 d\|_{L^2}^2\nonumber\\
&\le \frac14\sigma\|\nabla u\|_{L^2}^2+C(M_*)\Big(\Bbb E_0^\frac23
+\delta\Big)\sup_{0\le t\le T}\big(\sigma\|\nabla^2 d\|_{L^2}^2\big).
\end{align}
Owing to \eqref{w7}, Gagliardo-Nirenberg inequality, and \eqref{ee1}, we have
\begin{align}\label{f47}
\|\nabla u\|_{L^6}&\le C\big(\|\rho\dot{u}\|_{L^2}+\|\nabla\theta\|_{L^2}+\||\nabla d||\nabla^2d|\|_{L^2}\big)\nonumber\\
&\le C\big(\|\sqrt{\rho}\dot{u}\|_{L^2}+\|\nabla\theta\|_{L^2}+\|\nabla d\|_{L^\infty}\|\nabla^2d\|_{L^2}\big)\nonumber\\
&\le C\Big(\|\sqrt{\rho}\dot{u}\|_{L^2}+\|\nabla\theta\|_{L^2}
+\|\nabla^2d\|_{L^2}^\frac32\|\nabla^3d\|_{L^2}^\frac12\Big)\nonumber\\
&\le C\big(\|\sqrt{\rho}\dot{u}\|_{L^2}+\|\nabla\theta\|_{L^2}
+\|\nabla^3d\|_{L^2}+\|\nabla^2d\|_{L^2}^3\big).
\end{align}
Applying $\nabla$ to $\eqref{a1}_4$ gives that
\begin{align}\label{3.52}
\nabla d_t-\nabla\Delta d=-\nabla(u\cdot\nabla d)+\nabla(|\nabla d|^2d).
\end{align}
Then it follows from the standard $L^2$-estimate for elliptic equations, Sobolev's inequality, and Gagliardo-Nirenberg inequality that
\begin{align*}
\|\nabla^3d\|_{L^2}&\le C\|\nabla d_t\|_{L^2}+C\||u||\nabla^2 d|\|_{L^2}+C\||\nabla u||\nabla d|\|_{L^2}+C\|\nabla d\|_{L^6}^3+C\||\nabla d||\nabla^2d|\|_{L^2}\nonumber\\
&\le C\|\nabla d_t\|_{L^2}+C\|u\|_{L^6}\|\nabla^2d\|_{L^3}+C\|\nabla d\|_{L^\infty}\|\nabla u\|_{L^2}
+C\|\nabla^2d\|_{L^2}^3+C\|\nabla d\|_{L^\infty}\|\nabla^2d\|_{L^2}\nonumber\\
&\le C\|\nabla d_t\|_{L^2}+C(\|\nabla u\|_{L^2}+\|\nabla^2d\|_{L^2})\|\nabla^2d\|_{L^2}^\frac12
\|\nabla^3d\|_{L^2}^\frac12
+C\|\nabla^2d\|_{L^2}^3\nonumber\\
&\le \frac{1}{2}\|\nabla^3d\|_{L^2}+C\|\nabla d_t\|_{L^2}
+C\|\nabla u\|_{L^2}^2\|\nabla^2d\|_{L^2}+C\|\nabla^2d\|_{L^2}^3.
\end{align*}
Hence, we arrive at
\begin{align}\label{f49}
\|\nabla^3d\|_{L^2}\le C\|\nabla d_t\|_{L^2}
+C\|\nabla u\|_{L^2}^2\|\nabla^2d\|_{L^2}+C\|\nabla^2d\|_{L^2}^3.
\end{align}
Inserting \eqref{f49} into \eqref{f47}, one obtains that
\begin{align}\label{f50}
\|\nabla u\|_{L^6}\le C\|\sqrt{\rho}\dot{u}\|_{L^2}+C\|\nabla d_t\|_{L^2}+C\|\nabla\theta\|_{L^2}
+C\|\nabla^2d\|_{L^2}^3+C\|\nabla u\|_{L^2}^2\|\nabla^2d\|_{L^2},
\end{align}
which combined with H{\"o}lder's inequality and Sobolev's inequality leads to
\begin{align*}
& C\int\theta|\nabla u|^2dx \notag \\
&\le C\delta\|\theta\|_{L^6}\|\nabla u\|_{L^2}^\frac32\|\nabla u\|_{L^6}^\frac12\nonumber\\
&\le C\big(\|\sqrt{\rho}\dot{u}\|_{L^2}+\|\nabla d_t\|_{L^2}+\|\nabla\theta\|_{L^2}
+\|\nabla^2d\|_{L^2}^3+\|\nabla u\|_{L^2}^2\|\nabla^2d\|_{L^2}\big)^\frac12
\|\nabla u\|_{L^2}^\frac12\|\nabla\theta\|_{L^2}\|\nabla u\|_{L^2}\nonumber\\
&\le C\delta^{-\frac78}\|\nabla u\|_{L^2}^2\|\nabla\theta\|_{L^2}^2
+C\delta^{-\frac18}\|\nabla u\|_{L^2}^2+\delta^{\frac{15}{8}}\|\sqrt{\rho}\dot{u}\|_{L^2}^2
+\delta^{\frac{15}{8}}\|\nabla d_t\|_{L^2}^2\nonumber\\
&\quad+\delta^{\frac{15}{8}}\|\nabla\theta\|_{L^2}^2
+\delta^{\frac{15}{8}}(\|\nabla u\|_{L^2}^4+\|\nabla^2d\|_{L^2}^4)\|\nabla^2d\|_{L^2}^2.
\end{align*}
This along with \eqref{w1} yields that
\begin{align}\label{f51}
C\delta\sigma\int\theta|\nabla u|^2dx
&\le C\delta^\frac18\sup_{0\le t\le T}\big(\sigma\|\nabla u\|_{L^2}^2\big)
+\delta^\frac{15}{8}\sup_{0\le t\le T}\big(\delta\sigma\|\nabla\theta\|_{L^2}^2\big)\nonumber\\
&\quad+\delta^\frac{11}{8}\sup_{0\le t\le T}\big(\delta^\frac32\sigma\|\sqrt{\rho}\dot{u}\|_{L^2}^2\big)
+\delta^\frac{11}{8}\sup_{0\le t\le T}\big(\delta^\frac32\sigma\|\nabla d_t\|_{L^2}^2\big)\nonumber\\
&\quad+C(M_*)\delta^\frac{15}{8}\sup_{0\le t\le T}\big(\sigma\|\nabla^2d\|_{L^2}^2\big).
\end{align}
By H\"older's inequality, Sobolev's inequality, \eqref{f49}, and \eqref{w1}, we have
\begin{align*}
&C\int\theta|\Delta d+|\nabla d|^2d|^2dx\nonumber\\
&\le C\|\theta\|_{L^6}\|\Delta d+|\nabla d|^2d\|_{L^2}\|\Delta d+|\nabla d|^2d\|_{L^3}\nonumber\\
&\le C\|\nabla\theta\|_{L^2}\big(\|\nabla^2d\|_{L^2}+\|\nabla d\|_{L^4}^2\big)\big(\|\nabla^2d\|_{L^2}^\frac{1}{2}\|\nabla^3d\|_{L^2}^\frac{1}{2}
+\|\nabla d\|_{L^6}^2\big)\nonumber\\
&\le C\delta^{-\frac78}\|\nabla\theta\|_{L^2}^2\big(\|\nabla^2d\|_{L^2}^2+\|\nabla d\|_{L^2}^2\|\nabla^2d\|_{L^2}^3\big)
+\delta^\frac78\|\nabla^3d\|_{L^2}^2
+C\delta^\frac78\|\nabla^2d\|_{L^2}^2+\delta^\frac78\|\nabla^2d\|_{L^2}^4\nonumber\\
&\le C\delta^{-\frac78}\|\nabla^2d\|_{L^2}^2
+C\delta^\frac78\|\nabla d_t\|_{L^2}^2+C\delta^{-\frac78}\|\nabla d\|_{L^2}^2+C(M_*)\delta^\frac78\|\nabla^2d\|_{L^2}^2.
\end{align*}
This implies that
\begin{align*}
& C\delta\sigma\int\theta|\Delta d+|\nabla d|^2d|^2dx \notag \\
&\le C\sup_{0\le t\le T}\|\nabla d\|_{L^2}^2
+C\delta^\frac18\sup_{0\le t\le T}\big(\sigma\|\nabla^2 d\|_{L^2}^2\big)
+\delta^\frac38\sup_{0\le t\le T}\big(\delta^\frac32\sigma\|\nabla d_t\|_{L^2}^2\big),
\end{align*}
which together with \eqref{f46} and \eqref{f51} yields that
\begin{align}\label{f53}
&C\sigma\int\big(\rho\theta|\nabla u|+|\nabla u||\nabla d|^2+\delta\theta|\nabla u|^2+\delta\theta|\Delta d+|\nabla d|^2d|^2\big)dx\nonumber\\
&\le C\sup_{0\le t\le T}\big(\|\sqrt{\rho}\theta\|_{L^2}^2+\|\nabla d\|_{L^2}^2\big)
+\Big(\frac14+C\delta^\frac18\Big)\sup_{0\le t\le T}\big(\sigma\|\nabla u\|_{L^2}^2\big)\nonumber\\
&\quad+C\Big(\Bbb E_0^\frac23+\delta^\frac18\Big)\sup_{0\le t\le T}\big(\sigma\|\nabla^2 d\|_{L^2}^2\big)
+\delta^\frac38\sup_{0\le t\le T}(\delta^\frac32\sigma\|\nabla d_t\|_{L^2}^2)
+\delta^\frac{15}{8}\sup_{0\le t\le T}(\delta\sigma\|\nabla\theta\|_{L^2}^2)\nonumber\\
&\quad+\delta^\frac{11}{8}\sup_{0\le t\le T}\big(\delta^\frac32\sigma\|\sqrt{\rho}\dot{u}\|_{L^2}^2\big)\nonumber\\
&\le C\Big(t_0+\Bbb E_0^\frac45\Big)+\Big[\frac14+C_3\Big(\Bbb E_0^\frac23+\delta^\frac18\Big)\Big]
\sup_{0\le t\le T}\big(\sigma\|\nabla u\|_{L^2}^2\big)
+\delta^\frac{11}{8}\sup_{0\le t\le T}\big(\delta^\frac32\sigma\|\sqrt{\rho}\dot{u}\|_{L^2}^2\big)\nonumber\\
&\quad+C_3\Big(\Bbb E_0^\frac23+\delta^\frac18\Big)\sup_{0\le t\le T}\big(\sigma\|\nabla^2 d\|_{L^2}^2\big)
+C_3\Big(\Bbb E_0^\frac23\delta^{-1}+\delta^\frac{15}{8}\Big)\sup_{0\le t\le T}\big(\delta\sigma\|\nabla\theta\|_{L^2}^2\big)\nonumber\\
&\quad+\delta^\frac{3}{8}\sup_{0\le t\le T}\big(\delta^\frac32\sigma\|\nabla d_t\|_{L^2}^2\big).
\end{align}
Let
\begin{align*}
\delta=\delta_1\triangleq\min\left\{\frac{1}{\sqrt[3]{4}}, \frac{1}{2C_2}, \frac{c_v^2}{16}, \frac{c_v^2}{(8c_2c_4)^2}, \frac{1}{(8c_2c_4)^\frac23}, \frac{\mu}{2(2c_2+\mu)}, \frac{1}{(8C_3)^8}\right\},
\end{align*}
substituting \eqref{f44}, \eqref{f45}, and \eqref{f53} into \eqref{f43}, we conclude that
\begin{align}\label{f54}
&\sup_{0\le t\le T}\big[\sigma\big(\|\nabla u\|_{L^2}^2
+\|\nabla\theta\|_{L^2}^2+\|\nabla^2d\|_{L^2}^2+\|\sqrt{\rho}\dot{u}\|_{L^2}^2
+\|\nabla d_t\|_{L^2}^2\big)\big]\nonumber\\
&\quad+\int_0^T\sigma\big(\|\sqrt{\rho}\dot{u}\|_{L^2}^2+\|\nabla d_t\|_{L^2}^2
+\|\nabla^3d\|_{L^2}^2+\|\sqrt{\rho}\dot{\theta}\|_{L^2}^2
+\|\nabla\dot{u}\|_{L^2}^2+\|d_{tt}\|_{L^2}^2+\|\nabla^2d_t\|_{L^2}^2\big)dt\nonumber\\
&\le C\Big(t_0+\Bbb E_0^\frac45\Big)+C\int_0^T\sigma\big(\|\theta\nabla u\|_{L^2}^2+\|\nabla u\|_{L^2}^2
+\|\nabla u\|_{L^2}^4\|\nabla\theta\|_{L^2}^2\big)dt\nonumber\\
&\quad+C\delta\int_0^T\sigma\|\theta\|_{L^\infty}^2\big(\|\nabla^2d\|_{L^2}^2+\|\nabla d\|_{L^2}\|\nabla^2d\|_{L^2}^3\big)dt
+C\delta^\frac52\int_0^T\sigma\|\nabla^3d\|_{L^2}^4dt\nonumber\\
&\quad+C\delta^\frac32\int_0^T\sigma\|\nabla u\|_{L^4}^4dt+C\int_0^T\sigma\big(\|\nabla^2d\|_{L^2}^2+\|\nabla u\|_{L^2}^4
+\|\nabla d\|_{L^2}^2\big)\|\nabla^2d\|_{L^2}^2dt\nonumber\\
&\triangleq C\Big(t_0+\Bbb E_0^\frac45\Big)+\bar{R}_1
+\bar{R}_2+\bar{R}_3+\bar{R}_4+\bar{R}_5,
\end{align}
provided $$\Bbb E_0\le \varepsilon_2\triangleq\min\left\{\varepsilon_1,
\Big(\frac{\delta_1}{8C_3}\Big)^\frac32\right\}.$$

2. It follows from H\"older's, Young's and Gagliardo-Nirenberg inequalities, \eqref{ee1}, \eqref{g12}, and \eqref{f5} that
\begin{align*}
\bar{R}_1&\le \sup_{0\le t\le T}\big(\sigma\|\nabla u\|_{L^2}^2\big)\int_0^T\big(\|\nabla \theta\|_{L^2}^2
+\|\nabla^2\theta\|_{L^2}^2\big)dt
+C(M_*)\int_0^T\big(\|\nabla u\|_{L^2}^2+\|\nabla\theta\|_{L^2}^2\big)dt\nonumber\\
&\le C\Bbb E_0^\frac{2}{3}\sup_{t_0\le t\le T}\big(\|\nabla u\|_{L^2}^2+\|\nabla\theta\|_{L^2}^2
+\|\nabla^2d\|_{L^2}^2\big)+C\Big(t_0+\Bbb E_0^\frac45\Big)
+C\Bbb E_0^\frac{2}{3}\sup_{0\le t\le T}\big(\sigma\|\nabla u\|_{L^2}^2\big)\nonumber\\
&\le  C\Big(t_0+\Bbb E_0^\frac45\Big)+C\Bbb E_0^\frac23\sup_{0\le t\le T}\big[\sigma\big(\|\nabla u\|_{L^2}^2
+\|\nabla\theta\|_{L^2}^2+\|\nabla^2d\|_{L^2}^2\big)\big],\\
\bar{R}_2&\le C(M_*)\sup_{0\le t\le T}\big(\sigma\|\nabla^2d\|_{L^2}^2\big)\int_0^T\big(\|\nabla\theta\|_{L^2}^2
+\|\nabla^2\theta\|_{L^2}^2\big)dt\nonumber\\
\bar{R}_3&\le C\sup_{0\le t\le T}\sigma\|\nabla^3d\|_{L^2}^2\int_0^T\|\nabla^3d\|_{L^2}^2dt\nonumber\\
&\le C\Bbb E_0^\frac32\sup_{0\le t\le T}\big[\sigma\big(\|\nabla d_t\|_{L^2}^2+\|\nabla u\|_{L^2}^2+\|\nabla^2d\|_{L^2}^2\big)\big],\\
\bar{R}_4&\le C(M_*)\int_0^T\sigma\|\nabla u\|_{L^2}\big(\|\sqrt{\rho}\dot{u}\|_{L^2}+\|\nabla d_t\|_{L^2}+\|\nabla\theta\|_{L^2}
+\|\nabla^2d\|_{L^2}\big)^3dt\nonumber\\
&\le \sup_{0\le t\le T}\big(\sigma\|\nabla d_t\|_{L^2}^2\big)\int_0^T\|\nabla d_t\|_{L^2}^2dt
+\sup_{0\le t\le T}\big(\sigma\|\sqrt{\rho}\dot{u}\|_{L^2}^2\big)
\int_0^T\|\sqrt{\rho}\dot{u}\|_{L^2}^2dt\nonumber\\
&\quad+C\int_0^T\big(\|\nabla u\|_{L^2}^4
+\|\nabla\theta\|_{L^2}^4+\|\nabla^2d\|_{L^2}^4\big)dt
\nonumber\\
&\le C\Big(t_0+\Bbb E_0^\frac{4}{5}\Big)+C\Bbb E_0^\frac23\sup_{t_0\le t\le T}\big(\|\nabla u\|_{L^2}^2
+\|\nabla^2d\|_{L^2}^2+\|\nabla\theta\|_{L^2}^2\big)\nonumber\\
&\quad+2\Bbb E_0^\frac23\sup_{0\le t\le T}\big[\sigma\big(\|\sqrt{\rho}\dot{u}\|_{L^2}^2
+\|\nabla d_t\|_{L^2}^2\big)\big]\nonumber\\
&\le C\Bbb E_0^\frac{2}{3}\sup_{0\le t\le T}\big[\sigma\big(\|\nabla u\|_{L^2}^2
+\|\nabla^2d\|_{L^2}^2+\|\nabla\theta\|_{L^2}^2+\|\sqrt{\rho}\dot{u}\|_{L^2}^2
+\|\nabla d_t\|_{L^2}^2\big)\big]+C\Big(t_0+\Bbb E_0^\frac{4}{5}\Big),\\
\bar{R}_5&\le C\int_0^T\sigma\big(\|\nabla u\|_{L^2}^4+\|\nabla^2d\|_{L^2}^4+\|\nabla d\|_{L^2}^2\|\nabla^2d\|_{L^2}^2\big)dt\nonumber\\
&\le C(M_*)\int_0^T\big(\|\nabla u\|_{L^2}^2+\|\nabla^2d\|_{L^2}^2\big)dt\nonumber\\
&\le C\Big(t_0+\Bbb E_0^\frac45\Big)+C\Bbb E_0^\frac23\sup_{t_0\le t\le T}\big(\|\nabla u\|_{L^2}^2
+\|\nabla\theta\|_{L^2}^2+\|\nabla^2d\|_{L^2}^2\big)\nonumber\\
&\le C\Big(t_0+\Bbb E_0^\frac45\Big)+C\Bbb E_0^\frac23\sup_{0\le t\le T}\big[\sigma\big(\|\nabla u\|_{L^2}^2
+\|\nabla\theta\|_{L^2}^2+\|\nabla^2d\|_{L^2}^2\big)\big].
\end{align*}
Substituting the above estimates on $\bar{R}_1$--$\bar{R}_5$ into \eqref{f54}, one obtains that
\begin{align}
&\sup_{0\le t\le T}\big[\sigma\big(\|\nabla u\|_{L^2}^2
+\|\nabla\theta\|_{L^2}^2+\|\nabla^2d\|_{L^2}^2+\|\sqrt{\rho}\dot{u}\|_{L^2}^2
+\|\nabla d_t\|_{L^2}^2\big)\big]\nonumber\\
&\quad+\int_0^T\sigma\big(\|\sqrt{\rho}\dot{u}\|_{L^2}^2+\|\nabla d_t\|_{L^2}^2
+\|\nabla^3d\|_{L^2}^2+\|\sqrt{\rho}\dot{\theta}\|_{L^2}^2
+\|\nabla\dot{u}\|_{L^2}^2+\|d_{tt}\|_{L^2}^2+\|\nabla^2d_t\|_{L^2}^2\big)dt\nonumber\\
&\le C_4\Bbb E_0^\frac23\sup_{0\le t\le T}\big[\sigma\big(\|\nabla u\|_{L^2}^2
+\|\nabla^2d\|_{L^2}^2+\|\nabla\theta\|_{L^2}^2+\|\sqrt{\rho}\dot{u}\|_{L^2}^2
+\|\nabla d_t\|_{L^2}^2\big)\big]+C\Big(t_0+\Bbb E_0^\frac{4}{5}\Big).
\end{align}
This implies the desired \eqref{f42} provided that
\begin{align*}
\Bbb E_0\le \varepsilon_3\triangleq\min\left\{\varepsilon_2,\Big(\frac{1}{2C_4}\Big)
^\frac32\right\}.
\end{align*}
The proof of Lemma \ref{l33} is complete.
\hfill $\Box$

Motivated by \cite{D1997,WZ17}, we now show a time-independent upper bound for the density, which is the key to derive higher order estimates of the solution.
\begin{lemma}\label{l34}
Let $(\rho, u, \theta, d)$ be the smooth solution of \eqref{a1}--\eqref{a3} satisfying \eqref{ee1}, then it holds that
\begin{align}\label{3.49}
\sup_{0\le t\le T}\|\rho(t)\|_{L^\infty}\le \frac{3\bar{\rho}}{2}.
\end{align}
\end{lemma}
{\it Proof.}
1. For $(x, t)\in \Bbb R^3\times[0, T]$ and $\eta>0$, denote
\begin{align}
\rho^\eta(y,s)=\rho(y, s)+\eta\exp\left\{-\int_0^s\divv u(X(\tau; x, t), \tau)d\tau\right\}>0,
\end{align}
where
$X(s; x, t)$ is given by
\begin{align}
\left\{
\begin{array}{ll}
\displaystyle
\frac{d}{ds}X(s; x, t)=u(X(s; x, t), s), \quad 0\le s<t,\\
X(t; x, t)=x.
\end{array}
\right.
\end{align}
Then we have $\frac{d}{ds}(f(X(s; x, t), s)=(f_s+u\cdot\nabla f)(X(s; x, t), s)$, and so it follows from $\eqref{a1}_1$ that
 \begin{align}\label{tb1}
 \frac{d}{ds}\big(\log(\rho^\eta(X(s; x, t), s)\big)
 =-\divv u(X(s; x, t), s).
 \end{align}
Since $\divv u=\frac{F+p}{2\mu+\lambda}$ due to \eqref{f3}, we obtain from \eqref{tb1} that
\begin{align}
Y'(s)=g(s)+b'(s),
\end{align}
with
\begin{align*}
&Y(s)\triangleq\log(\rho^\eta(X(s; x, t), s)), \ g(s)\triangleq-\frac{p(X(s; x, t), s)}{2\mu+\lambda},\
b(s)\triangleq-\frac{1}{2\mu+\lambda}\int_0^sF(X(\tau; x, t), \tau)d\tau.
\end{align*}

2. Rewrite $\eqref{a1}_2$ as
\begin{align}\label{dxy}
\partial_t\big[\Delta^{-1}\divv(\rho u)\big]-(2\mu+\lambda)
\divv u+p=-\Delta^{-1}\divv\divv(\rho u\otimes u)
-\Delta^{-1}\divv\divv(M(d)),
\end{align}
where
\begin{align*}
M(d)\triangleq\nabla d\odot\nabla d-\frac{1}{2}|\nabla d|^2\mathbb{I}_3\ \ \text{with}\ \mathbb{I}_3\ \text{being}\ 3\times3\ \text{identical matrix}.
\end{align*}
Thus we get from \eqref{f3} and \eqref{dxy} that
\begin{align}
F(X(\tau; x, t), \tau)
&=-\big[(-\Delta)^{-1}\divv(\rho u)\big]_{\tau}-(-\Delta)^{-1}\divv\divv(\rho u\otimes u)
-(-\Delta)^{-1}\divv\divv(M(d))\nonumber\\
&=-\big[(-\Delta)^{-1}\divv(\rho u)\big]_{\tau}-u\cdot\nabla(-\Delta)^{-1}\divv(\rho u)
+u\cdot\nabla(-\Delta)^{-1}\divv(\rho u)\nonumber\\
&\quad-(-\Delta)^{-1}\divv\divv(\rho u\otimes u)
-(-\Delta)^{-1}\divv\divv(M(d))\nonumber\\
&=-\frac{d}{d\tau}\big[(-\Delta)^{-1}\divv(\rho u)\big]+[u^i, R_{ij}](\rho u^j)-(-\Delta)^{-1}\divv\divv(M(d)).
\end{align}
Here $[u^i, R_{ij}]\triangleq u^iR_{ij}-R_{ij}u^i$ with $R_{ij}=\partial_i(-\Delta)^{-1}\partial_j$ being the Riesz transform on $\Bbb R^3$ (we refer the reader to \cite[Chapter 3]{S70} for its definition and properties). This yields that
\begin{align}\label{xyx}
b(t)-b(0)
&=\frac{1}{2\mu+\lambda}\int_0^t\Big[\frac{d}{d\tau}\big[(-\Delta)^{-1}\divv (\rho u)\big]-[u^i, R_{ij}](\rho u^j)
+(-\Delta)^{-1}\divv\divv(M(d))\Big]d\tau\nonumber\\
&=\frac{1}{2\mu+\lambda}(-\Delta)^{-1}\divv (\rho u)-\frac{1}{2\mu+\lambda}(-\Delta)^{-1}\divv (\rho_0u_0)
-\frac{1}{2\mu+\lambda}\int_0^t[u^i, R_{ij}](\rho u^j)d\tau\nonumber\\
&\quad-\frac{1}{2\mu+\lambda}\int_0^t(-\Delta)^{-1}\divv \divv (M(d))d\tau\nonumber\\
&\le \frac{1}{2\mu+\lambda}\|(-\Delta)^{-1}\divv (\rho u)\|_{L^\infty}
+\frac{1}{2\mu+\lambda}\|(-\Delta)^{-1}\divv (\rho_0u_0)\|_{L^\infty}\nonumber\\
&\quad+\frac{1}{2\mu+\lambda}\int_0^t\|[u^i, R_{ij}](\rho u^j)\|_{L^\infty}d\tau+
\frac{1}{2\mu+\lambda}\int_0^t\|(-\Delta)^{-1}\divv \divv (M(d))\|_{L^\infty}d\tau\nonumber\\
&\triangleq Z_1+Z_2+Z_3+Z_4.
\end{align}
Applying Gagliardo-Nirenberg inequality, Sobolev's inequality, the $L^q$-estimates of the Riesz transform,
 and H\"older's inequality, we obtain from \eqref{ee1}, \eqref{f5}, and \eqref{f42} that
\begin{align}\label{3.60}
Z_1&\le \frac{C}{2\mu+\lambda}\|(-\Delta)^{-1}\divv (\rho u)\|_{L^6}^\frac13
\|\nabla(-\Delta)^{-1}\divv (\rho u)\|_{L^4}^\frac23\nonumber\\
&\le C\|\rho u\|_{L^2}^\frac12\|\rho u\|_{L^6}^\frac12
\le C(\bar{\rho})\|\sqrt{\rho}u\|_{L^2}^\frac12\|\nabla u\|_{L^2}^\frac12\nonumber\\
&\le C\Big(t_0+\Bbb E_0^\frac45\Big)^\frac14.
\end{align}
Similarly to $Z_1$, we have
\begin{align}\label{3.61}
Z_2 \le  C\Big(t_0+\Bbb E_0^\frac{4}{5}\Big)^\frac14.
\end{align}
For $Z_3$, we deduce from Gagliardo-Nirenberg inequality, Sobolev's inequality, \eqref{ee1}, \eqref{f50}, and \eqref{w1} that
\begin{align}\label{uuy}
Z_3 & \le \frac{C}{2\mu+\lambda}\int_0^t\|[u^i, R_{ij}](\rho u^j)\|_{L^{12}}^\frac12\|\nabla[u^i, R_{ij}](\rho u^j)\|_{L^4}^\frac12d\tau\nonumber\\
&\le C\int_0^t\big(\|\nabla u\|_{L^3}\|\rho u\|_{L^{12}}\big)^\frac12\big(\|\nabla u\|_{L^6}\|\rho u\|_{L^{12}}\big)^\frac12d\tau\nonumber\\
&\le C\bar{\rho}\int_0^t\|\nabla u\|_{L^2}^{\frac12}
\|\nabla u\|_{L^6}^{\frac12}\|u\|_{L^{12}}d\tau\nonumber\\
&\le C\int_0^t\|\nabla u\|_{L^2}^{\frac12}
\|\nabla u\|_{L^6}^{\frac12}\|u\|_{L^6}^{\frac12}\|\nabla u\|_{L^3}^{\frac12}d\tau\nonumber\\
&\le C\int_0^t\|\nabla u\|_{L^2}^{\frac54}
\|\nabla u\|_{L^6}^\frac34d\tau\nonumber\\
&\le C\int_0^t\big(\|\nabla u\|_{L^2}^2+\|\nabla u\|_{L^6}^2\big)d\tau\nonumber\\
&\le C\int_0^t\big(\|\nabla u\|_{L^2}^2
+\|\sqrt{\rho}\dot{u}\|_{L^2}^2+\|\nabla d_t\|_{L^2}^2+\|\nabla\theta\|_{L^2}^2
+\|\nabla^2d\|_{L^2}^6+\|\nabla u\|_{L^2}^4\|\nabla^2d\|_{L^2}^2\big)d\tau\nonumber\\
&\le C\left(\int_0^{t_0}+\int_{t_0}^t\right)\big(\|\nabla u\|_{L^2}^2
+\|\sqrt{\rho}\dot{u}\|_{L^2}^2+\|\nabla d_t\|_{L^2}^2+\|\nabla\theta\|_{L^2}^2
+\|\nabla^2d\|_{L^2}^6+\|\nabla u\|_{L^2}^4\|\nabla^2d\|_{L^2}^2\big)d\tau\nonumber\\
&\le C\int_0^t\sigma\big(\|\nabla u\|_{L^2}^2
+\|\sqrt{\rho}\dot{u}\|_{L^2}^2+\|\nabla d_t\|_{L^2}^2+\|\nabla\theta\|_{L^2}^2+\|\nabla^2d\|_{L^2}^2\big)d\tau
+C(M_*)t_0.
\end{align}
Combining \eqref{f5} and \eqref{f42}, one obtains that
\begin{align}\label{sp}
\sup_{0\le t\le T}\big(\|\sqrt{\rho}u\|_{L^2}^2+\|\sqrt{\rho}\theta\|_{L^2}^2+\|\nabla d\|_{L^2}^2\big)
+\int_0^T\big(\|\nabla u\|_{L^2}^2+\|\nabla\theta\|_{L^2}^2+\|\nabla^2d\|_{L^2}^2\big)dt\le C\Big(t_0+\Bbb E_0^\frac45\Big).
\end{align}
This together with \eqref{uuy} and \eqref{f42} yields
\begin{align}\label{f68}
Z_3\le C\Big(t_0+\Bbb E_0^\frac45\Big).
\end{align}
For $Z_4$, by \eqref{ee1}, \eqref{w1}, \eqref{f42}, and the following Gagliardo-Nirenberg inequality
\begin{align*}
\|\nabla d\|_{L^\infty}\le C\|\nabla^2 d\|_{L^2}^\frac12\|\nabla^3 d\|_{L^2}^\frac12,
\end{align*}
we have
\begin{align}\label{3.68}
Z_4&\le \frac{1}{2\mu+\lambda}\int_0^t\|(-\Delta)^{-1}\divv \divv (M(d))\|_{L^\infty}d\tau\nonumber\\
&\le C\left(\int_0^{t_0}+\int_{t_0}^t\right)\big(\|\nabla^2d\|_{L^2}^2+\|\nabla^3d\|_{L^2}^2\big)d\tau\nonumber\\
&\le C(M_*)t_0+\int_0^t\sigma\big(\|\nabla^2d\|_{L^2}^2+\|\nabla^3d\|_{L^2}^2\big)d\tau\nonumber\\
&\le C\Big(t_0+\Bbb E_0^\frac45\Big).
\end{align}
Substituting \eqref{3.60}, \eqref{3.61}, \eqref{f68}, and \eqref{3.68} into \eqref{xyx}, we find that
\begin{align}\label{3.75}
b(t)-b(0)&\le C_5\Big(t_0+\Bbb E_0^\frac45\Big)^\frac14\le \log\frac32,
\end{align}
provided that
\begin{align*}
\Bbb E_0\le \varepsilon_4\triangleq\min\left\{\varepsilon_3, \Big(\frac{\log\frac32}{2^\frac14C_5}\Big)^5\right\}.
\end{align*}
Consequently, integrating \eqref{tb1} with respect to $s$ over $[0, t]$, we get from \eqref{3.75} that
\begin{align*}
\log\rho^\eta(x, t)=\log[\rho_0(X(t; x, 0))+\eta]+\int_0^tg(\tau)d\tau+b(t)-b(0)
\le \log(\bar{\rho}+\eta)+\log\frac32,
\end{align*}
due to \eqref{1.6} and $g=-\frac{p}{2\mu+\lambda}=-\frac{R\rho\theta}{2\mu+\lambda}\leq0$.
Letting $\eta\rightarrow 0^+$, we obtain the desired \eqref{3.49}.
\hfill $\Box$

With Lemmas \ref{l31}--\ref{l34} at hand, we are now in a position to prove Proposition \ref{p1}.

\textbf{Proof of Proposition \ref{p1}.}
1. It follows from \eqref{w1} and \eqref{f42} that
\begin{align}
&\sup_{0\le t\le T}\big(\|\nabla u\|_{L^2}^2
+\|\nabla\theta\|_{L^2}^2+\|\nabla^2d\|_{L^2}^2+\|\sqrt{\rho}\dot{u}\|_{L^2}^2
+\|\nabla d_t\|_{L^2}^2\big)\nonumber\\
&\le \sup_{0\le t\le t_0}\big(\|\nabla u\|_{L^2}^2
+\|\nabla\theta\|_{L^2}^2+\|\nabla^2d\|_{L^2}^2+\|\sqrt{\rho}\dot{u}\|_{L^2}^2
+\|\nabla d_t\|_{L^2}^2\big)\nonumber\\
&\quad+\sup_{0\le t\le T}\sigma\big(\|\nabla u\|_{L^2}^2
+\|\nabla\theta\|_{L^2}^2+\|\nabla^2d\|_{L^2}^2+\|\sqrt{\rho}\dot{u}\|_{L^2}^2
+\|\nabla d_t\|_{L^2}^2\big)\nonumber\\
&\le M_*+C_6\Big(t_0+\Bbb E_0^\frac{4}{5}\Big).
\end{align}
Thus, we have
\begin{align}\label{f72}
\sup_{0\le t\le T}\big(\|\nabla u\|_{L^2}^2+\|\nabla\theta\|_{L^2}^2+\|\nabla^2d\|_{L^2}^2
+\|\sqrt{\rho}\dot{u}\|_{L^2}^2
+\|\nabla d_t\|_{L^2}^2\big)\le \frac{3M_*}{2},
\end{align}
provided
\begin{align*}
t_0+\Bbb E_0^\frac{4}{5}\le \frac{M_*}{2C_6}.
\end{align*}

2. From \eqref{ee1}, \eqref{f22}, and \eqref{f50} that
\begin{align*}
\|\nabla^2\theta\|_{L^2}^2&\le C\big(\|\sqrt{\rho}\dot{\theta}\|_{L^2}^2+\|\theta\nabla u\|_{L^2}^2+\|\nabla u\|_{L^4}^4
+\|\Delta d+|\nabla d|^2d\|_{L^4}^4\big)\nonumber\\
&\le C\big(\|\sqrt{\rho}\dot{\theta}\|_{L^2}^2+\|\theta\|_{L^\infty}^2\|\nabla u\|_{L^2}^2+\|\nabla u\|_{L^2}\|\nabla u\|_{L^6}^3
+\|\nabla^2d\|_{L^4}^4+\|\nabla d\|_{L^8}^8\big)\nonumber\\
&\le \frac{1}{2}\|\nabla^2\theta\|_{L^2}^2
+C\|\sqrt{\rho}\dot{\theta}\|_{L^2}^2+C\|\nabla\theta\|_{L^2}^2+C\|\nabla u\|_{L^2}^2
+C\|\sqrt{\rho}\dot{u}\|_{L^2}^2
+C\|\nabla^2 d\|_{L^2}^2+C\|\nabla^3d\|_{L^2}^4\nonumber\\
&\le  \frac{1}{2}\|\nabla^2\theta\|_{L^2}^2
+C\|\sqrt{\rho}\dot{\theta}\|_{L^2}^2+C\|\nabla\theta\|_{L^2}^2+C\|\nabla u\|_{L^2}^2
+C\|\sqrt{\rho}\dot{u}\|_{L^2}^2
+C\|\nabla^2 d\|_{L^2}^2+C\|\nabla d_t\|_{L^2}^4,
\end{align*}
that is
\begin{align}\label{3.74}
\|\nabla^2\theta\|_{L^2}^2&\le C\|\sqrt{\rho}\dot{\theta}\|_{L^2}^2+C\|\nabla\theta\|_{L^2}^2+C\|\nabla u\|_{L^2}^2
+C\|\sqrt{\rho}\dot{u}\|_{L^2}^2+C\|\nabla^2 d\|_{L^2}^2+C\|\nabla d_t\|_{L^2}^4.
\end{align}
Thus, using \eqref{w1}, \eqref{f5}, and \eqref{f42}, we infer
\begin{align}\label{f74}
&\int_0^T\big(\|\sqrt{\rho}\dot{u}\|_{L^2}^2+\|\nabla d_t\|_{L^2}^2
+\|\nabla^3d\|_{L^2}^2+\|\nabla\theta\|_{L^2}^2+\|\nabla^2\theta\|_{L^2}^2\big)dt\nonumber\\
&=\left(\int_0^{t_0}+\int_{t_0}^T\right)\big(\|\sqrt{\rho}\dot{u}\|_{L^2}^2+\|\nabla d_t\|_{L^2}^2
+\|\nabla^3d\|_{L^2}^2+\|\nabla\theta\|_{L^2}^2+\|\nabla^2\theta\|_{L^2}^2\big)dt\nonumber\\
&\le M_*t_0+\int_0^T\sigma\Big(\|\sqrt{\rho}\dot{u}\|_{L^2}^2+\|\nabla d_t\|_{L^2}^2
+\|\nabla^3d\|_{L^2}^2+\|\nabla\theta\|_{L^2}^2+\|\nabla^2\theta\|_{L^2}^2\Big)dt\nonumber\\
&\quad+\int_0^T\sigma\big(\|\sqrt{\rho}\dot{\theta}\|_{L^2}^2+\|\nabla^2d\|_{L^2}^2+\|\nabla u\|_{L^2}^2\big)dt
+C\sup_{0\le t\le T}(\|\nabla d_t\|_{L^2}^2)\int_0^T\sigma\|\nabla d_t\|_{L^2}^2dt\nonumber\\
&\le C_7\Big(t_0+\Bbb E_0^\frac{4}{5}\Big)\le \Bbb E_0^\frac{2}{3},
\end{align}
provided
\begin{align*}
t_0+\Bbb E_0^\frac{4}{5}\le \frac{\Bbb E_0^\frac23}{C_7}.
\end{align*}

3. In the view of \eqref{f5} and \eqref{f42}, we have
\begin{align}
\int_0^T\|\nabla^2d\|_{L^2}^2dt&\le C\Big(t_0+\Bbb E_0^\frac45\Big)+C\Bbb E_0^\frac23\sup_{t_0\le t\le T}\big(\|\nabla u\|_{L^2}^2+\|\nabla\theta\|_{L^2}^2
+\|\nabla^2d\|_{L^2}^2\big)\nonumber\\
&\le C\Big(t_0+\Bbb E_0^\frac45\Big)
+C\Bbb E_0^\frac23\sup_{0\le t\le T}\sigma\big(\|\nabla u\|_{L^2}^2
+\|\nabla\theta\|_{L^2}^2+\|\nabla^2d\|_{L^2}^2\big)\nonumber\\
&\le C_8\Big(t_0+\Bbb E_0^\frac45\Big)\le E_0^\frac13,
\end{align}
provided
\begin{align*}
t_0+\Bbb E_0^\frac45\le \frac{\Bbb E_0^\frac13}{C_8}.
\end{align*}

Finally, \eqref{ee2} is valid if we select
\begin{align}
\Bbb E_0\le \varepsilon_0\quad {\rm and}\quad t_0\le t_1\triangleq\min\left\{\frac{T_*}{2}, \varepsilon_0\right\}
\end{align}
with
\begin{align*}
\varepsilon_0\triangleq\min\left\{\frac{\bar\rho\|\nabla\theta_0\|_{L^2}^6}{54c_v^4\pi^4}, \frac{1}{(2C_1)^\frac32}, \Big(\frac{\delta_1}{8C_3}\Big)^\frac32, \Big(\frac{1}{2C_4}\Big)^\frac32,
\Big(\frac{\log\frac32}{2^\frac14C_5}\Big)^5,
\Big(\frac{M_*}{2C_6}\Big)^\frac54, \Big(\frac{1}{2C_7}\Big)^\frac{15}{2}, \Big(\frac{1}{2C_8}\Big)^\frac{15}{7}\right\}.
\end{align*}
The proof of Proposition \ref{p1} is complete.
\hfill $\Box$

\begin{lemma}\label{l35}
Let $(\rho, u, \theta, d)$ be the smooth solution of \eqref{a1}--\eqref{a3} satisfying \eqref{ee1}, then it holds that
\begin{align}
\sup_{0\le t\le T}\big(\|\sqrt{\rho}\dot{\theta}\|_{L^2}^2+\|\nabla^2\theta\|_{L^2}^2+\|\nabla^3d\|_{L^2}^2\big)
+\int_0^T\|\nabla\dot{\theta}\|_{L^2}^2dt\le C.
\end{align}
\end{lemma}
{\it Proof.}
1. Operating $\partial_t+\divv(u\cdot)$ to $\eqref{a1}_3$ gives rise to
\begin{align}\label{t8}
c_v\rho(\dot{\theta}+u\cdot\nabla\dot{\theta})
&=\kappa\Delta\dot{\theta}+\kappa[\divv u\Delta\theta-
\partial_i(\partial_i\cdot\nabla\theta)-\partial_iu\cdot\nabla\partial_i\theta]
-\rho\dot{\theta}\divv u-\rho\theta\divv\dot{u}\nonumber\\
&\quad+\Big(\lambda(\divv u)^2+\frac{\mu}{2}|\nabla u+(\nabla u)^{tr}|^2\Big)\divv u
+2\lambda\big(\divv\dot{u}-\partial_ku^l\partial_lu^k\big)\divv u\nonumber\\
&\quad+\rho\theta\partial_ku^l\partial_lu^k
+\mu\big(\partial_iu^j+\partial_ju^i\big)\big(\partial_i\dot{u}^j
+\partial_j\dot{u}^i-\partial_iu^k\partial_ku^j-\partial_ju^k\partial_ku^i\big)\nonumber\\
&\quad+\partial_t|\Delta d+|\nabla d|^2d|^2
+\divv\big(|\Delta d+|\nabla d|^2d|^2u\big).
\end{align}
Then multiplying \eqref{t8} by $\dot{\theta}$ and integration by parts lead to
\begin{align}\label{t9}
&\frac{c_v}{2}\frac{d}{dt}\int\rho|\dot{\theta}|^2dx
+\kappa\int|\nabla\dot{\theta}|^2dx\nonumber\\
&\le C\int|\nabla u|(|\nabla^2\theta||\dot{\theta}|+|\nabla\theta||\nabla\dot{\theta}|)dx
+C\int|\nabla u|^2|\dot{\theta}|(|\nabla u|+\rho\theta)dx\nonumber\\
&\quad+C\int\rho|\dot{\theta}|^2|\nabla u|dx+C\int\rho\theta|\nabla\dot{u}||\dot{\theta}|dx
+C\int|\nabla u||\nabla\dot{u}||\dot{\theta}|dx\nonumber\\
&\quad+\int\dot{\theta}\partial_t|\Delta d+|\nabla d|^2d|^2dx+\int\dot{\theta}
\divv(|\Delta d+|\nabla d|^2d|^2u)dx
\triangleq\sum_{i=1}^7K_i.
\end{align}
By virtue of H\"older's, Young's, and Gagliardo-Nirenberg inequalities,
we can bound each term $K_i$ as follows
\begin{align*}
K_1&\le C\|\nabla u\|_{L^3}\big(\|\nabla^2\theta\|_{L^2}\|\dot{\theta}\|_{L^6}
+\|\nabla\dot{\theta}\|_{L^2}\|\nabla\theta\|_{L^6}\big)
\le \frac{\kappa}{14}\|\nabla\dot{\theta}\|_{L^2}^2+C\|\nabla u\|_{L^3}^2\|\nabla^2\theta\|_{L^2}^2,\\
K_2&\le C\|\rho\theta\|_{L^3}\|\nabla u\|_{L^4}^2\|\dot{\theta}\|_{L^6}
+C\|\nabla u\|_{L^3}\|\nabla u\|_{L^4}^2\|\dot{\theta}\|_{L^6}\nonumber\\
&\le \frac{\kappa}{14}\|\nabla\dot{\theta}\|_{L^2}^2
+C\big(\|\rho\theta\|_{L^3}^2+\|\nabla u\|_{L^3}^2\big)\|\nabla u\|_{L^4}^4,\\
K_3&\le C(\bar{\rho})\|\sqrt{\rho}\dot{\theta}\|_{L^2}\|\dot{\theta}\|_{L^6}\|\nabla u\|_{L^3}\le \frac{\kappa}{14}\|\nabla\dot{\theta}\|_{L^2}^2+C\|\sqrt{\rho}\dot{\theta}\|_{L^2}^2\|\nabla u\|_{L^3}^2,\\
K_4&\le C\|\rho\theta\|_{L^3}\|\nabla\dot{u}\|_{L^2}\|\dot{\theta}\|_{L^6}
\le \frac{\kappa}{14}\|\nabla\dot{\theta}\|_{L^2}^2
+C\|\rho\theta\|_{L^3}^2\|\nabla\dot{u}\|_{L^2}^2,\\
K_5&\le C\|\nabla u\|_{L^3}\|\nabla\dot{u}\|_{L^2}\|\dot{\theta}\|_{L^6}
\le \frac{\kappa}{14}\|\nabla\dot{\theta}\|_{L^2}^2+C\|\nabla u\|_{L^3}^2\|\nabla\dot{u}\|_{L^2}^2,\\
K_6&\le C\int\big(|\nabla^2d_t||\nabla^2d|+|\nabla^2d||\nabla d||\nabla d_t|+|\nabla d|^2|d_t||\nabla^2d|
\big)|\dot{\theta}|dx\nonumber\\
&\quad+\int\big(|\nabla^2d_t||\nabla d|^2+|\nabla d|^3|\nabla d_t|+|\nabla d|^4|d_t|\big)|\dot{\theta}|dx\nonumber\\
&\le C\|\nabla^2d\|_{L^3}\|\nabla^2d_t\|_{L^2}\|\dot{\theta}\|_{L^6}
+C\|\nabla^2d\|_{L^2}\|\nabla d\|_{L^6}\|\nabla d_t\|_{L^6}\|\dot{\theta}\|_{L^6}\nonumber\\
&\quad+C\|d_t\|_{L^6}\|\dot{\theta}\|_{L^6}\|\nabla^2d\|_{L^3}\|\nabla d\|_{L^6}^2
+C\|\dot{\theta}\|_{L^6}\|\nabla^2d_t\|_{L^2}\|\nabla d\|_{L^6}^2\nonumber\\
&\quad+C\|\dot{\theta}\|_{L^6}\|\nabla d\|_{L^6}^2\|\nabla d_t\|_{L^3}
+C\|\dot{\theta}\|_{L^6}\|d_t\|_{L^6}\|\nabla d\|_{L^6}\|\nabla d\|_{L^6}^3\nonumber\\
&\le \frac{\kappa}{14}\|\nabla\dot{\theta}\|_{L^2}^2
+C\|\nabla^2d\|_{L^3}^2\|\nabla^2d_t\|_{L^2}^2+C\|\nabla^2d\|_{L^2}^4\|\nabla^2d_t\|_{L^2}^2\nonumber\\
&\quad+C\|\nabla^2d\|_{L^3}^2\|\nabla^2d\|_{L^2}^4\|\nabla d_t\|_{L^2}^2
+C\|\nabla^2d\|_{L^2}^4\|\nabla^2d_t\|_{L^2}^2\nonumber\\
&\quad+C\big(\|\nabla^2 d\|_{L^2}^4+\|\nabla^2d\|_{L^2}^8\big)\|\nabla d_t\|_{L^2}^2\nonumber\\
&\le \frac{\kappa}{14}\|\nabla\dot{\theta}\|_{L^2}^2
+C\big(\|\nabla^2d\|_{L^3}^2+\|\nabla^2d\|_{L^2}^2\big)\|\nabla^2d_t\|_{L^2}^2
+C\|\nabla^2d\|_{L^2}^2\|\nabla d_t\|_{L^2}^2\nonumber\\
K_7&\le C\|u\|_{L^6}\|\nabla\dot{\theta}\|_{L^2}\big(\|\nabla^2d\|_{L^6}^2+\|\nabla d\|_{L^{12}}^4\big)\nonumber\\
&\le \frac{\kappa}{14}\|\nabla\dot{\theta}\|_{L^2}^2
+C\|\nabla u\|_{L^2}^2\|\nabla^3d\|_{L^2}^4+C\|\nabla u\|_{L^2}^2\|\nabla d\|_{L^3}^4\|\nabla^2 d\|_{L^6}^4\nonumber\\
&\le \frac{\kappa}{14}\|\nabla\dot{\theta}\|_{L^2}^2
+C\|\nabla u\|_{L^2}^2\|\nabla^3d\|_{L^2}^4.
\end{align*}
Substituting the above estimates on  $K_i\ (i=1, 2, \ldots, 7)$ into \eqref{t9} yields that
\begin{align}\label{f80}
\frac{d}{dt}\|\sqrt{\rho}\dot{\theta}\|_{L^2}^2
+\frac{\kappa}{2}\|\nabla\dot{\theta}\|_{L^2}^2
&\le C\big(\|\sqrt{\rho}\theta\|_{L^3}^2+\|\nabla u\|_{L^3}^2+\|\nabla^3d\|_{L^2}^2+\|\nabla^2d\|_{L^2}^2\big)\nonumber\\
&\quad\times\big(
\|\sqrt{\rho}\dot{\theta}\|_{L^2}^2+\|\nabla^2\theta\|_{L^2}^2+\|\nabla d_t\|_{H^1}^2
+\|\nabla u\|_{L^4}^4+\|\nabla\dot{u}\|_{L^2}^2\big).
\end{align}

2. It follows from \eqref{f5}, \eqref{w1}, \eqref{f50}, and \eqref{f72} that
\begin{align}\label{f81}
&\sup_{0\le t\le T}\big(\|\sqrt{\rho}\theta\|_{L^3}^2+\|\nabla u\|_{L^3}^2+\|\nabla^3d\|_{L^2}^2+\|\nabla^2d\|_{L^2}^2\big)\nonumber\\
&\le C(\bar{\rho})\sup_{0\le t\le T}\big(\|\sqrt{\rho}\theta\|_{L^2}^2+
\|\nabla \theta\|_{L^2}^2+\|\nabla u\|_{L^2}^2+\|\nabla u\|_{L^6}^2+\|\nabla^3d\|_{L^2}^2+\|\nabla^2d\|_{L^2}^2\big)\nonumber\\
&\le C(\bar{\rho})\sup_{0\le t\le T}\big(\|\sqrt{\rho}\theta\|_{L^2}^2+
\|\nabla \theta\|_{L^2}^2+\|\nabla u\|_{L^2}^2+\|\nabla d_t\|_{L^2}^2+\|\sqrt{\rho}\dot{u}\|_{L^2}^2+\|\nabla^2d\|_{L^2}^2\big)\nonumber\\
&\le C(M_*, \bar{\rho}),
\end{align}
and from \eqref{w1}, \eqref{ee2}, \eqref{f42}, and \eqref{f74} that
\begin{align}\label{f82}
&\int_0^T\big(\|\sqrt{\rho}\dot{\theta}\|_{L^2}^2
+\|\nabla^2\theta\|_{L^2}^2+\|\nabla^2d_t\|_{L^2}^2+\|\nabla d_t\|_{L^2}^2
+\|\nabla u\|_{L^4}^4+\|\nabla\dot{u}\|_{L^2}^2\big)\nonumber\\
&\le \left(\int_0^{t_0}+\int_{t_0}^T\right)\big(\|\sqrt{\rho}\dot{\theta}\|_{L^2}^2
+\|\nabla^2\theta\|_{L^2}^2+\|\nabla d_t\|_{H^1}^2
+\|\nabla u\|_{L^4}^4+\|\nabla\dot{u}\|_{L^2}^2\big)\nonumber\\
&\le C(M_*),
\end{align}
where we have used the following fact
\begin{align}\label{fff}
\|\nabla u\|_{L^4}^4&\le C\|\nabla u\|_{L^2}\|\nabla u\|_{L^6}^3\nonumber\\
&\le C(M_*)\|\nabla u\|_{L^2}\big(\|\sqrt{\rho}\dot{u}\|_{L^2}+\|\nabla d_t\|_{L^2}+\|\nabla\theta\|_{L^2}
+\|\nabla^2d\|_{L^2}\big)^3\nonumber\\
&\le C\|\nabla u\|_{L^2}^4+C\|\sqrt{\rho}\dot{\theta}\|_{L^2}^4+C\|\nabla^2d\|_{L^2}^4
+C\|\nabla\theta\|_{L^2}^4+C\|\nabla d_t\|_{L^2}^4.
\end{align}
Combining \eqref{f81} and \eqref{f82}, and integrating \eqref{f80} over
$[0, T]$, we have
\begin{align}\label{f83}
\sup_{0\le t\le T}\|\sqrt{\rho}\dot{\theta}\|_{L^2}^2+\int_0^T\|\nabla\dot{\theta}\|_{L^2}^2dt\le C.
\end{align}
This together with \eqref{f83} and \eqref{f72} leads to
\begin{align}
\|\nabla^2\theta\|_{L^2}^2&\le C\|\sqrt{\rho}\dot{\theta}\|_{L^2}^2+C\|\nabla\theta\|_{L^2}^2+C\|\nabla u\|_{L^2}^2
+C\|\sqrt{\rho}\dot{u}\|_{L^2}^2+C\|\nabla^2 d\|_{L^2}^2+C\|\nabla d_t\|_{L^2}^4\le C.
\end{align}
Thus, the proof of Lemma \ref{l35} is complete.
\hfill $\Box$

As a straightforward consequence of Proposition \ref{p1} and Lemmas \ref{l31}--\ref{l35}, we arrive at
\begin{corollary}
Under the assumptions in Proposition \ref{p1}, then there exists a
positive constant $C$ depending only on $\mu$, $\lambda$, $c_v$, $\kappa$, $R$, and the initial data such that the following estimates hold true
\begin{align}\label{tt}
&\sup_{0\le t\le T}\big(\|\nabla u\|_{H^1}^2+\|\nabla\theta\|_{H^1}^2+\|\nabla d\|_{H^2}^2+\|\nabla d_t\|_{L^2}^2+\|\sqrt{\rho}\dot{u}\|_{L^2}^2
+\|\sqrt{\rho}\dot{\theta}\|_{L^2}^2\big)\nonumber\\
&\quad+\int_0^T\big(\|\nabla\dot{u}\|_{L^2}^2+\|\nabla\dot{\theta}\|_{L^2}^2
+\|d_{tt}\|_{L^2}^2+\|\nabla^2d_t\|_{L^2}^2
+\|\nabla u\|_{L^2}^2+\|\nabla\theta\|_{H^1}^2+\|\nabla^2d\|_{H^1}^2\big)dt\le C.
\end{align}
\end{corollary}

\subsection{Time-dependent estimates}

\begin{lemma}\label{l36}
Under the same assumption in Theorem \ref{thm1}, we have
\begin{align}
\sup_{0\le t\le T}\big(\|\nabla\rho\|_{L^2\cap L^q}+\|\nabla^2u\|_{L^2}\big)+\int_0^T\big(\|\nabla^2u\|_{L^6}^2+\|\nabla^2\theta\|_{L^6}^2
+\|\nabla^4d\|_{L^2}^2\big)dt\le C(T).
\end{align}
\end{lemma}
{\it Proof.}
1. Taking spatial derivative $\nabla$ on the mass equation \eqref{a1}$_1$
leads to
\begin{equation*}
\partial_{t}\nabla\rho+u\cdot\nabla^2\rho
+\nabla u\cdot\nabla\rho+\divv u\nabla\rho+\rho\nabla\divv u=0.
\end{equation*}
For $q\in (3, 6)$, multiplying the above equality by $q|\nabla\rho|^{q-2}\nabla\rho$ gives that
\begin{align*}
(|\nabla\rho|^q)_t+\divv(|\nabla\rho|^qu)+(q-1)|\nabla\rho|^q\divv u+q|\nabla\rho|^{q-2}(\nabla\rho)^{tr}\nabla u(\nabla\rho)+q\rho|\nabla\rho|^{q-2}\nabla\rho\cdot\nabla\divv u=0.
\end{align*}
Thus, integration by parts over $\mathbb{R}^3$ yields that
\begin{align}\label{0t86}
\frac{d}{dt}\|\nabla\rho\|_{L^q}
\le C\big(\|\nabla u\|_{L^\infty}+1\big)\|\nabla\rho\|_{L^q}+\|\nabla^2u\|_{L^q}.
\end{align}
Applying the standard $L^q$-estimate to \eqref{a1}$_2$ leads to
\begin{align}\label{t87}
\|\nabla^2u\|_{L^q}\le C(\|\rho\dot{u}\|_{L^q}+\|\nabla (\rho\theta)\|_{L^q}+\||\nabla d||\nabla^2 d|\|_{L^q})
\le C(1+\|\nabla\dot{u}\|_{L^2}+\|\nabla\rho\|_{L^q}),
\end{align}
due to \eqref{tt}. So we obtain from \eqref{0t86} and \eqref{t87} that
\begin{align}\label{t86}
\frac{d}{dt}\|\nabla\rho\|_{L^q}\le C\big(\|\nabla u\|_{L^\infty}+1\big)\|\nabla\rho\|_{L^q}+C\big(1+\|\nabla\dot{u}\|_{L^2}
+\|\nabla\rho\|_{L^q}\big).
\end{align}

2.  It follows from the Gargliardo-Nirenberg inequality and \eqref{w5} that
\begin{align}\label{t88}
\|\divv u\|_{L^\infty}+\|\curl u\|_{L^\infty}&\le C(\|F\|_{L^\infty}+\|\rho\theta\|_{L^\infty}+\|w\|_{L^\infty})\nonumber\\
&\le C\|F\|_{L^6}^\frac12\|\nabla F\|_{L^6}^\frac12+C\|w\|_{L^6}^\frac12\|\nabla w\|_{L^6}^\frac12+C\|\nabla\theta\|_{L^2}^\frac12\|\nabla^2\theta\|_{L^2}^\frac12\nonumber\\
&\le C(1+\|\nabla F\|_{L^2}+\|\nabla F\|_{L^6}+\|\nabla w\|_{L^2}+\|\nabla w\|_{L^6})\nonumber\\
&\le C(1+\|\rho\dot{u}\|_{L^2}+\||\nabla d||\nabla^2d|\|_{L^2}+\|\rho\dot{u}\|_{L^6}
+\||\nabla d||\nabla^2d|\|_{L^6})\nonumber\\
&\le C(\bar{\rho})(1+\|\sqrt{\rho}\dot{u}\|_{L^2}+\|\nabla\dot{u}\|_{L^2}+\|\nabla^3d\|_{L^2}^2)\nonumber\\
&\le C(\bar{\rho})(1+\|\nabla\dot{u}\|_{L^2}),
\end{align}
which together with Lemma \ref{l24}, \eqref{t87}, and \eqref{t88} yields that
\begin{align}\label{t89}
\|\nabla u\|_{L^\infty}&\le C\big(\|\divv u\|_{L^\infty}
+\|\curl u\|_{L^\infty}\big)\log\big(e+\|\nabla^2u\|_{L^q}\big)
+C\|\nabla u\|_{L^2}+C\nonumber\\
&\le C\big(1+\|\nabla\dot{u}\|_{L^2}\big)\log\big(
e+\|\nabla\dot{u}\|_{L^2}+\|\nabla\rho\|_{L^q}\big)+C.
\end{align}
Substituting \eqref{t89} into \eqref{t86} yields that
\begin{align}\label{t90}
f'(t)\leq Cg(t)f(t)\ln f(t),
\end{align}
where
\begin{align*}
f(t)\triangleq e+\|\nabla\rho\|_{L^q},\ g(t)\triangleq 1+\|\nabla\dot{u}\|_{L^2}.
\end{align*}
It thus follows from \eqref{t90}, Gronwall's inequality, and \eqref{tt} that
\begin{align}\label{t91}
\sup_{0\le t\le T}\|\nabla\rho\|_{L^q}\le C(T).
\end{align}
It should be noted that \eqref{0t86} also holds true for $q=2$, hence we then derive from that Gronwall's inequality, \eqref{t89}, \eqref{t91}, \eqref{t87},
and \eqref{tt} that
\begin{align}\label{zt91}
\sup_{0\le t\le T}\|\nabla\rho\|_{L^2}\le C(T).
\end{align}
This combined with \eqref{t91}, \eqref{f72}, \eqref{f81}, and \eqref{f82} implies that
\begin{align}\label{t92}
\sup_{0\le t\le T}\|\nabla^2u\|_{L^2}&\le C\sup_{0\le t\le T}\big(\|\sqrt{\rho}\dot{u}\|_{L^2}
+\|\nabla(\rho\theta)\|_{L^2}+\||\nabla d||\nabla^2d|\|_{L^2}\big)\nonumber\\
&\le C+C\sup_{0\le t\le T}\big(\|\nabla\rho\|_{L^2}+\|\nabla\theta\|_{L^2}+\|\nabla^3 d\|_{L^2}\big)
\le C,
\end{align}
and
\begin{align}
\int_0^T\|\nabla^2u\|_{L^6}^2dt\le C\int_0^T\big(1+\|\nabla\dot{u}\|_{L^2}^2+\|\nabla\rho\|_{L^6}^2\big)dt\le C.
\end{align}

3. Taking the operator $\nabla^2$ to $\eqref{a1}_4$, we get
\begin{align}
-\nabla^2\Delta d=\nabla^2\big(|\nabla d|^2d-u\cdot\nabla d-d_t\big).
\end{align}
Using the $L^2$-estimates of elliptic system, we derive that
\begin{align}
\|\nabla^4 d\|_{L^2}&\le C\big(\|\nabla^2 d_t\|_{L^2}+\|\nabla^2(u\cdot\nabla d)\|_{L^2}+\|\nabla^2(|\nabla d|^2d)\|_{L^2}\big)\nonumber\\
&\le C\|\nabla^2 d_t\|_{L^2}+C\||\nabla^2u||\nabla d|\|_{L^2}+C\||\nabla u||\nabla^2d|\|_{L^2}
+C\||u||\nabla^3d|\|_{L^2}\nonumber\\
&\quad+C\||\nabla d|^2|\nabla^2d|\|_{L^2}+C\||\nabla^2d|^2\|_{L^2}+C\||\nabla d||\nabla^3d|\|_{L^2}\nonumber\\
&\le C\|\nabla^2 d_t\|_{L^2}+C\|\nabla d\|_{L^\infty}\|\nabla^2u\|_{L^2}+C\|\nabla u\|_{L^3}\|\nabla^2d\|_{L^6}
\nonumber\\
&\quad+C\|u\|_{L^6}\|\nabla^3d\|_{L^3}+C\|\nabla d\|_{L^6}^2\|\nabla^2d\|_{L^6}
+C\|\nabla^2 d\|_{L^\infty}\|\nabla^2d\|_{L^2}
+C\|\nabla d\|_{L^6}\|\nabla^3d\|_{L^3}\nonumber\\
&\le C\|\nabla^2d\|_{L^2}^\frac{1}{2}\|\nabla^3d\|_{L^2}^\frac{1}{2}\|\nabla^2u\|_{L^2}
+C\|\nabla u\|_{L^2}^\frac{1}{2}\|\nabla^2u\|_{L^2}^\frac{1}{2}\|\nabla^3d\|_{L^2}\nonumber\\
&\quad+C\big(\|\nabla u\|_{L^2}^2+\|\nabla^2 d\|_{L^2}^2\big)\|\nabla^3d\|_{L^2}
+C\|\nabla^2d_t\|_{L^2}+\frac{1}{2}\|\nabla^4d\|_{L^2},
\end{align}
which yields
\begin{align}
\|\nabla^4 d\|_{L^2}^2&\le C\|\nabla^2d\|_{L^2}\|\nabla^3d\|_{L^2}\|\nabla^2u\|_{L^2}^2
+C\|\nabla u\|_{L^2}\|\nabla^2 u\|_{L^2}\|\nabla^3d\|_{L^2}^2\nonumber\\
&\quad+C\big(\|\nabla u\|_{L^2}^4+\|\nabla^2 d\|_{L^2}^4\big)\|\nabla^3d\|_{L^2}^2
+C\|\nabla^2d_t\|_{L^2}^2,
\end{align}
from which, \eqref{t92}, and \eqref{f72}, one obtains that
\begin{align}\label{t97}
\int_0^T\|\nabla^4d\|_{L^2}^2dt&\le C\int_0^T(\|\nabla^3d\|_{L^2}^2+\|\nabla^2d\|_{L^2}^2+\|\nabla^2d_t\|_{L^2}^2)dt
\le C.
\end{align}
This together with \eqref{tt}, \eqref{f72}, \eqref{t92}, and \eqref{t97} yields
\begin{align}
\int_0^T\|\nabla^2\theta\|_{L^6}^2dt&\le C\int_0^T\big(\|\sqrt{\rho}\dot{\theta}\|_{L^6}^2
+\|\theta\nabla u\|_{L^6}^2+\|\nabla u\|_{L^{12}}^4
+\|\Delta d+|\nabla d|^2d\|_{L^{12}}^4\big)dt\nonumber\\
&\le C(\bar{\rho})\int_0^T\big(\|\nabla\dot{\theta}\|_{L^2}^2+\|\theta\|_{L^\infty}^2\|\nabla^2u\|_{L^2}^2
+\|\nabla u\|_{L^3}^2\|\nabla^2u\|_{L^6}^2\big)dt\nonumber\\
&\quad+C\int_0^T\big(\||\nabla^2d||\nabla^3d|\|_{L^2}^2+\||\nabla^2d||\nabla d|^3\|_{L^2}^2\big)dt\nonumber\\
&\le C\int_0^T\big(\|\nabla\dot{\theta}\|_{L^2}^2+\|\nabla^2\theta\|_{L^2}^2
+\|\nabla^2u\|_{L^6}^2+\|\nabla^4d\|_{L^2}^2+\|\nabla^3d\|_{L^2}^3\big)dt
\le C.
\end{align}
The proof of Lemma \ref{l36} is finished.
\hfill $\Box$

\section{Proof of Theorem \ref{thm1}}
In terms of Lemma \ref{l21}, the Cauchy problem \eqref{a1}--\eqref{a3} has a unique solution $(\rho, u, \theta, d)$ on $\Bbb R^3\times(0, T_*]$
for some small $T_*>0$. With the a priori estimates obtained in Section 3, we can extend the $(\rho, u, \theta, d)$ to all time. In fact, by the
continuity of density and Lemma \ref{l21}, we see that \eqref{ee1} and \eqref{ee2} are valid in $[0, T_1]$ with
$T_1\triangleq\min\left\{T_*,\frac{2\Bbb E_0^\frac23}{M_*}\right\}$.

Set
\begin{align}
T^*=\sup\{T~|~\eqref{ee1}~is~valid\}.
\end{align}
Then $T^*\ge T_1$. We next claim that
\begin{align}
T^*=\infty.
\end{align}
Suppose, by contradiction,
that $T^*<\infty$. Note that, all
the a {\it priori} estimates obtained in Section 3 are uniformly bounded for any $t<T^*$. Hence,
we define
\begin{align}
(\rho, u, \theta, d)(x, T^*)=\lim_{t\rightarrow T^*}(\rho, u, \theta, d)(x, t).
\end{align}
Furthermore, standard arguments yield that $\rho\dot{u}$, $\rho\dot{\theta}\in C([0, T^*); L^2)$, which implies
\begin{align}
(\rho\dot{u}, \rho\dot{\theta})(x, T^*)=\lim_{t\rightarrow T^*}(\rho\dot{u}, \rho\dot{\theta})(x, t)\in L^2.
\end{align}
Hence,
\begin{align*}
\big[-\mu\Delta u-(\mu+\lambda)\nabla{\rm div}u+\nabla(R\rho\theta)+\nabla d\cdot\Delta d\big]|_{t=T^*}&=\sqrt{\rho}(x, T^*)g_1(x),\\
\big[-\kappa\Delta \theta-\lambda(\divv u)^2-\frac{\mu}{2}|\nabla u+(\nabla u)^{tr}|^2-|\Delta d+|\nabla d|^2d|^2\big]|_{t=T^*}&=\sqrt{\rho}(x, T^*)g_2(x),
\end{align*}
with
\begin{align*}
g_1(x)=\left\{
\begin{array}{ll}
\displaystyle
\rho^{-\frac{1}{2}}(x, T^*)(\rho\dot{u})(x, T^*), &{\rm for}~ x\in\{x|\rho(x, T^*)>0\},\\[3pt]
0,  &{\rm for}~ x\in\{x|\rho(x, T^*)=0\},
\end{array}
\right.
\end{align*}
and
\begin{align*}
g_2(x)=\left\{
\begin{array}{ll}
\displaystyle
\rho^{-\frac{1}{2}}(x, T^*)(c_v\rho\dot{\theta}+\rho\theta{\rm div}u)(x, T^*), &{\rm for}~ x\in\{x|\rho(x, T^*)>0\},\\[3pt]
0,  &{\rm for}~ x\in\{x|\rho(x, T^*)=0\},
\end{array}
\right.
\end{align*}
satisfying $(g_1, g_2)\in L^2$ due to \eqref{tt}. Thus, $(\rho, u, \theta, d)$  satisfy compatibility conditions \eqref{qqw1} and \eqref{qqw}. Therefore, we can extend the local strong solution
beyond $T^*$ by taking $(\rho, u, \theta, d)(x, T^*)$ as the initial data and applying the Lemma \ref{l21}.
 This contradicts the assumption of $T^*$. This completes the proof of Theorem
\ref{thm1}. \hfill $\Box$

\section{Proof of Theorem \ref{thm2}}

We divide the proof of Theorem \ref{thm2} into several steps.

\textbf{Step 1. Estimate for $\|\nabla^2d(t)\|_{L^2}$}. We obtain from \eqref{f20} and \eqref{tt} that
\begin{align*}
\frac{d}{dt}\|\nabla^2 d\|_{L^2}^2+\|\nabla d_t\|_{L^2}^2+\|\nabla^3 d\|_{L^2}^2
\le C\left(\|\nabla^2 d\|_{L^2}^2+\|\nabla u\|_{L^2}^2\right)
\|\nabla^2 d\|_{L^2}^2,
\end{align*}
which leads to
\begin{align*}
\frac{d}{dt}\left(t\|\nabla^2 d\|_{L^2}^2\right)+t\big(\|\nabla d_t\|_{L^2}^2+\|\nabla^3 d\|_{L^2}^2\big)
&\le C\left(\|\nabla^2 d\|_{L^2}^2+\|\nabla u\|_{L^2}^2\right)
\big(t\|\nabla^2 d\|_{L^2}^2\big)+\|\nabla^2 d\|_{L^2}^2.
\end{align*}
This along with Gronwall's inequality and \eqref{tt} yields
\begin{align}\label{dd1}
\sup_{0\le t\le T}\big(t\|\nabla^2 d\|_{L^2}^2\big)+\int_0^Tt\left(\|\nabla d_t\|_{L^2}^2+\|\nabla ^3 d\|_{L^2}^2\right)dt\le C.
\end{align}

\textbf{Step 2. Estimate for $\|\sqrt{\rho}u(t)\|_{L^2}$ and $\|\sqrt{\rho}\theta(t)\|_{L^2}$}. It follows from \eqref{f4}, H{\"o}lder's inequality, Sobolev's inequality, \eqref{lz3.1}, and \eqref{3.49} that
\begin{align*}
&\frac{1}{2}\frac{d}{dt}\int\rho|u|^2dx+\int\big[\mu|\nabla u|^2+(\mu+\lambda)(\divv u)^2
\big]dx\nonumber\\
&=\int R\rho\theta\divv udx-\int (u\cdot\nabla)d\cdot\Delta ddx\nonumber\\
&\le  (\mu+\lambda)\|\divv u\|_{L^2}^2+C\|\rho^{\frac23}\|_{L^{\frac32}}\|\rho^{\frac43}\|_{L^\infty}
\|\theta\|_{L^6}^2
+C\|u\|_{L^6}\|\nabla d\|_{L^3}\|\nabla^2d\|_{L^2}\nonumber\\
&\le (\mu+\lambda)\|\divv u\|_{L^2}^2
+C\|\rho\|_{L^1}^\frac23\|\rho\|_{L^\infty}^\frac43
\|\nabla\theta\|_{L^2}^2+C\|\nabla u\|_{L^2}
\|\nabla d\|_{L^2}^{\frac12}\|\nabla^2d\|_{L^2}^{\frac32}\nonumber\\
&\le (\mu+\lambda)\|\divv u\|_{L^2}^2+\frac{\mu}{2}\|\nabla u\|_{L^2}^2
+C\|\nabla d\|_{L^2}\|\nabla^2d\|_{L^2}^3+C\|\nabla\theta\|_{L^2}^2,
\end{align*}
which together with \eqref{tt} implies that
\begin{align}\label{x1}
\frac{d}{dt}\|\sqrt{\rho}u\|_{L^2}^2+\mu\|\nabla u\|_{L^2}^2
\le C\|\nabla^2d\|_{L^2}\|\nabla d\|_{H^2}^2+c_6\|\nabla\theta\|_{L^2}^2.
\end{align}
In view of \eqref{f7}, we infer from H{\"o}lder's inequality, Gagliardo-Nirenberg inequality, and \eqref{tt} that
\begin{align*}
&\frac{c_v}{2}\frac{d}{dt}\int\rho\theta^2dx+\kappa\int|\nabla\theta|^2dx\nonumber\\
&=\int\theta\Big(\lambda(\divv u)^2+\frac{\mu}{2}|\nabla u+(\nabla u)^{tr}|^2-R\rho\theta\divv u\Big)dx
+\int|\Delta d+|\nabla d|^2d|^2\theta dx\nonumber\\
&\le \frac{\kappa\mu}{8c_6}\int|\nabla u|^2dx+C\int\theta^2\big(|\nabla u|^2+\rho\theta^2\big)dx+C\int|\nabla d||\nabla^3 d|\theta dx\nonumber\\
&\quad+C\int|\nabla\theta||\nabla d||\nabla^2d|dx+C\int|\nabla d|^4\theta dx+C\int|\nabla d|^2|\nabla^2d|\theta dx\nonumber\\
&\le \frac{\kappa\mu}{8c_6}\|\nabla u\|_{L^2}^2+C\|\theta\|_{L^\infty}^2\big(\|\nabla u\|_{L^2}^2+\|\nabla\theta\|_{L^2}^2\big)
+C\|\theta\|_{L^6}\|\nabla^3d\|_{L^2}\|\nabla d\|_{L^3}\nonumber\\
&\quad +C\|\nabla\theta\|_{L^2}\|\nabla d\|_{L^3}\|\nabla^2d\|_{L^6}
+C\|\nabla d\|_{L^6}^3\|\nabla d\|_{L^3}\|\theta\|_{L^6}+C\|\nabla d\|_{L^3}^2\|\theta\|_{L^6}\|\nabla^2d\|_{L^6}\nonumber\\
&\le \frac{\kappa\mu}{8c_6}\|\nabla u\|_{L^2}^2+C\|\nabla\theta\|_{L^2}\|\nabla^2\theta\|_{L^2}\big(\|\nabla u\|_{L^2}^2+\|\nabla\theta\|_{L^2}^2\big)+C\|\nabla\theta\|_{L^2}\|\nabla d\|_{L^2}^{\frac12}
\|\nabla^2d\|_{L^2}^{\frac12}\|\nabla^3d\|_{L^2}
\nonumber\\
&\quad +C\|\nabla\theta\|_{L^2}\|\nabla d\|_{L^2}^{\frac12}
\|\nabla^2d\|_{L^2}^{\frac72}
+C\|\nabla\theta\|_{L^2}\|\nabla d\|_{L^2}
\|\nabla^2d\|_{L^2}\|\nabla^3d\|_{L^2} \nonumber\\
&\le \frac{\kappa\mu}{8c_6}\|\nabla u\|_{L^2}^2+C\Big(t_0+\Bbb E_0^\frac45\Big)^\frac12\big(\|\nabla u\|_{L^2}^2+\|\nabla\theta\|_{L^2}^2\big)
+\frac{\kappa}{2}\|\nabla\theta\|_{L^2}^2
+C\|\nabla^2d\|_{L^2}\|\nabla d\|_{H^2}^2,
\end{align*}
where we have used, for $t\ge t_0$,
\begin{align*}
\|\nabla\theta\|_{L^2}\le C\sup_{t_0\le t\le T}\big(\sigma\|\nabla\theta\|_{L^2}^2\big)^\frac12
\le C\Big(t_0+\Bbb E_0^\frac45\Big)^\frac12,
\end{align*}
due to \eqref{f42}. Hence, we derive that
\begin{align}\label{x2}
&c_v\frac{d}{dt}\|\sqrt{\rho}\theta\|_{L^2}^2+\kappa\|\nabla\theta\|_{L^2}^2\nonumber\\
&\le \frac{\kappa\mu}{4c_6}\|\nabla u\|_{L^2}^2+C\|\nabla^2d\|_{L^2}
\|\nabla^2d\|_{H^1}^2
+C_9\Big(t_0+\Bbb E_0^\frac45\Big)^\frac12\big(\|\nabla u\|_{L^2}^2+\|\nabla\theta\|_{L^2}^2\big).
\end{align}
Combining \eqref{x1} and \eqref{x2} multiplied by $\frac{2c_6}{\kappa}$, and choosing
$$\Bbb E_0\le \varepsilon_4\triangleq\min\left\{\varepsilon_0,
\Big(\frac{\mu}{2^\frac52C_9}\Big)^{10}, \Big(\frac{\kappa}{2^\frac52C_9}\Big)^{10}\right\},$$
we get that
\begin{align}\label{x4}
\frac{d}{dt}\big(\|\sqrt{\rho}u\|_{L^2}^2+\|\sqrt{\rho}\theta\|_{L^2}^2\big)
+\|\nabla u\|_{L^2}^2+\|\nabla\theta\|_{L^2}^2
\le C\|\nabla^2d\|_{L^2}\|\nabla^2d\|_{H^1}^2.
\end{align}
Multiplying \eqref{x4} by $t^\frac12$ and integrating the resulting inequality over $[0, T]$, we obtain from \eqref{tt} and \eqref{dd1} that
\begin{align}\label{5.5}
&\sup_{0\le t\le T}
\big[t^\frac12\big(\|\sqrt{\rho}u\|_{L^2}^2+\|\sqrt{\rho}\theta\|_{L^2}^2\big)\big]
+\int_0^Tt^\frac12\big(\|\nabla u\|_{L^2}^2+\|\nabla\theta\|_{L^2}^2\big)dt\nonumber\\
&\le \sup_{0\le t\le T}\big(\|\sqrt{\rho}u\|_{L^2}^2+\|\sqrt{\rho}\theta\|_{L^2}^2\big)\int_0^{t_0}t^{-\frac12}dt
+C\int_{t_0}^T\big(\|\nabla u\|_{L^2}^2+\|\nabla\theta\|_{L^2}^2\big)dt\nonumber\\
&\quad+\sup_{0\le t\le T}\big(t\|\nabla^2d\|_{L^2}^2\big)^\frac12
\int_0^T\|\nabla^2d\|_{H^1}^2dt
\le C.
\end{align}

\textbf{Step 3. Estimate for $\|\nabla u(t)\|_{L^2}$},
\textbf{$\|\nabla\theta(t)\|_{L^2}$}, \textbf{and $\|\sqrt{\rho}\dot{u}(t)\|_{L^2}$}. Multiplying \eqref{f12} by $t^\frac12$, and integrating it by parts over $[0, T]$, we have
\begin{align}\label{k6}
& t^\frac12\int\Big[\frac{\mu}{2}|\nabla u|^2+\frac{\kappa\delta}{2}|\nabla\theta|^2+|\nabla^2d|^2
+\Big(\frac{2c_2c_3}{\mu}+1\Big)\delta^\frac{3}{2}\rho|\dot{u}|^2
+2c_2\delta^\frac32|\nabla d_t|^2\Big]dx\nonumber\\
&\quad+\int_0^T\int t^\frac12\Big(\rho|\dot{u}|^2+\frac14|\nabla d_t|^2+\frac{1}{4}|\nabla^3d|^2+\frac{c_v\delta}{4}\rho|\dot{\theta}|^2
+\frac{\mu c_2\delta^\frac{3}{2}}{2}|\nabla \dot{u}|^2+2c_2\delta^\frac{3}{2}|d_{tt}|^2+c_2\delta^\frac32|\nabla^2d_t|^2\Big)dxdt
\nonumber\\
&\le \int_0^T\int t^{-\frac12}\Big(\frac{\mu}{2}|\nabla u|^2+\frac{\kappa\delta}{2}|\nabla\theta|^2+
|\nabla^2d|^2+\Big(\frac{2c_2c_3}{\mu}+1\Big)\delta^\frac32\rho|\dot{u}|^2
+2c_2\delta^\frac32|\nabla d_t|^2\Big)dxdt\nonumber\\
&\quad+t^\frac12\int\Big(R\rho\theta{\rm div}u+(\nabla d\odot\nabla d):\nabla u
-\frac{1}{2}{\rm div}u|\nabla d|^2+\lambda\delta\theta({\rm div}u)^2+\frac{\mu\delta}{2}\theta|\nabla u+(\nabla u)^{tr}|^2\Big)dx\nonumber\\
&\quad+\int_0^T\int t^{-\frac12}\Big(R\rho\theta{\rm div}u+(\nabla d\odot\nabla d):\nabla u
-\frac12{\rm div}u|\nabla d|^2+\lambda\delta\theta({\rm div}u)^2+\frac{\mu\delta}{2}\theta|\nabla u+(\nabla u)^{tr}|^2\Big)dxdt\nonumber\\
&\quad+t^\frac12\int\big(\delta\theta|\Delta d+|\nabla d|^2d|^2\big)dx+\int_0^T\int t^{-\frac12}\big(\delta\theta|\Delta d+|\nabla d|^2d|^2\big)dxdt\nonumber\\
&\quad+\delta^{-\frac32}\int_0^Tt^\frac12\big(\|\theta\nabla u\|_{L^2}^2+\|\nabla u\|_{L^2}^2+\|\nabla u\|_{L^2}^4\|\nabla\theta\|_{L^2}^2\big)dt
\nonumber\\
&\quad+\int_0^Tt^\frac12\big(6\delta^\frac32\|\nabla u\|_{L^4}^4+\delta^\frac52\|\nabla^3d\|_{L^2}^4
+C\|\nabla^2d\|_{L^2}^4+C\|\nabla u\|_{L^2}^4\|\nabla^2d\|_{L^2}^2
+C\|\nabla d\|_{L^2}\|\nabla^2d\|_{L^2}^3\big)dt\nonumber\\[3pt]
&\quad+C\delta\int_0^Tt^\frac12\|\theta\|_{L^\infty}^2
\big(\|\nabla^2d\|_{L^2}^2+\|\nabla d\|_{L^2}\|\nabla^2d\|_{L^2}^3\big)dt.
\end{align}
It follows from \eqref{w1}, \eqref{f5}, and \eqref{tt} that
\begin{align}\label{k7}
&\int_0^T\int t^{-\frac12}\Big[\frac{\mu}{2}|\nabla u|^2+\frac{\kappa\delta}{2}|\nabla\theta|^2+
|\nabla^2d|^2+\Big(\frac{2c_2c_3}{\mu}+1\Big)\delta^\frac32\rho|\dot{u}|^2
+2c_2\delta^\frac32|\nabla d_t|^2\Big]dxdt\nonumber\\
&\le C\sup_{0\le t\le t_0}\big(\|\nabla u\|_{L^2}^2+\|\nabla\theta\|_{L^2}^2+\|\nabla^2d\|_{L^2}^2
+\|\sqrt{\rho}\dot{u}\|_{L^2}^2+\|\nabla d_t\|_{L^2}^2\big)\int_{0}^{t_0}t^{-\frac12}dt\nonumber\\
&\quad+C\int_{t_0}^T\big(\|\nabla u\|_{L^2}^2+\|\nabla\theta\|_{L^2}^2+\|\nabla^2d\|_{L^2}^2
+\|\sqrt{\rho}\dot{u}\|_{L^2}^2+\|\nabla d_t\|_{L^2}^2\big)dt\nonumber\\
&\le C.
\end{align}
Moreover, we deduce from \eqref{3.49}, \eqref{5.5}, and \eqref{tt} that
\begin{align}
&Ct^\frac12\int\big(\rho\theta|{\rm div}u|+\theta|\Delta d+|\nabla d|^2d|^2+|\nabla u||\nabla d|^2+\theta|\nabla u|^2\big)dx\nonumber\\
&\le \frac{1}{16}\sup_{0\le t\le T}\big(t^\frac12\|\nabla u\|_{L^2}^2\big)+C\sup_{0\le t\le T}\|\theta\|_{L^\infty}\sup_{0\le t\le T}\big(t\|\nabla^2d\|_{L^2}^2\big)^\frac12\nonumber\\
&\quad+C\sup_{0\le t\le T}\big(\|\nabla d\|_{L^2}\|\nabla^2d\|_{L^2}^2\big)\sup_{0\le t\le T}\big(t\|\nabla^2d\|_{L^2}^2\big)^\frac12+C\sup_{0\le t\le T}t^\frac12\|\sqrt{\rho}\theta\|_{L^2}^2\nonumber\\
&\quad+C\delta^\frac18\sup_{0\le t\le T}\big(t^\frac12\|\nabla u\|_{L^2}^2\big)
+\delta^\frac{15}{8}\sup_{0\le t\le T}\big(t^\frac12\delta\|\nabla\theta\|_{L^2}^2\big)
+\delta^\frac{11}{8}\sup_{0\le t\le T}\big(t^\frac12\delta^\frac32\|\sqrt{\rho}\dot{u}\|_{L^2}^2\big)\nonumber\\
&\quad+\delta^\frac{11}{8}\sup_{0\le t\le T}\big(t^\frac12\delta^\frac32\|\nabla d_t\|_{L^2}^2\big)
+\delta^\frac{15}{8}\sup_{0\le t\le T}\big(t^\frac12\|\nabla^2d\|_{L^2}^2\big)\nonumber\\
&\le \Big(\frac{1}{16}+C\delta^\frac18\Big)\sup_{0\le t\le T}\big(t^\frac12\|\nabla u\|_{L^2}^2\big)+\delta^\frac{15}{8}\sup_{0\le t\le T}\big(t^\frac12\delta\|\nabla\theta\|_{L^2}^2\big)
+\delta^\frac{11}{8}\sup_{0\le t\le T}\big(t^\frac12\delta^\frac32\|\sqrt{\rho}\dot{u}\|_{L^2}^2\big)\nonumber\\
&\quad+\delta^\frac{11}{8}\sup_{0\le t\le T}\big(t^\frac12\delta^\frac32\|\nabla d_t\|_{L^2}^2\big)
+\delta^\frac{15}{8}\sup_{0\le t\le T}\big(t^\frac12\|\nabla^2d\|_{L^2}^2\big)
+C.
\end{align}
We obtain from \eqref{w1}, H\"older's inequality, \eqref{lz3.1}, \eqref{3.49}, Gagliardo-Nirenberg inequality,
and \eqref{tt} that
\begin{align}
&\int_0^T\int t^{-\frac12}(\rho\theta|{\rm div}u|+\theta|\nabla u|^2+|\nabla u||\nabla d|^2+\theta|\Delta d+|\nabla d|^2d|^2)dxdt\nonumber\\
&\le\left(\int_0^{t_0}+\int_{t_0}^T\right)\int t^{-\frac12}\big(\rho\theta|{\rm div}u|+\theta|\nabla u|^2+|\nabla u||\nabla d|^2+\theta|\Delta d+|\nabla d|^2d|^2\big)dxdt\nonumber\\
&\le C\sup_{0\le t\le t_0}\big(\|\sqrt{\rho}\theta\|_{L^2}^2+\|\nabla\theta\|_{L^2}^2+\|\nabla^2\theta\|_{L^2}^2
+\|\nabla u\|_{L^2}^2+\|\nabla^2d\|_{L^2}^2+\|\nabla d\|_{L^2}^2\big)\int_0^{t_0}t^{-\frac12}dt\nonumber\\
&\quad+C\int_{t_0}^T\big(\|\sqrt{\rho}\theta\|_{L^2}\|\nabla u\|_{L^2}
+\|\theta\|_{L^\infty}\|\nabla u\|_{L^2}^2
+\|\nabla u\|_{L^2}\|\nabla d\|_{L^4}^2\big)dt\nonumber\\
&\quad+\int_{t_0}^T\|\theta\|_{L^\infty}\big(\|\nabla^2 d\|_{L^2}^2+\|\nabla d\|_{L^4}^4\big)dt\nonumber\\
&\le C\int_{t_0}^T\big(\|\nabla u\|_{L^2}^2+\|\nabla\theta\|_{L^2}^2\big)dt
+C\sup_{t_0\le t\le T}\|\theta\|_{L^\infty}\int_{t_0}^T\|\nabla u\|_{L^2}^2dt\nonumber\\
&\quad+C\sup_{t_0\le t\le T}\big(\|\nabla d\|_{L^2}^2\|\nabla u\|_{L^2}^2\big)
\int_{t_0}^T\|\nabla u\|_{L^2}^2dt+C\int_{t_0}^T\|\nabla^2d\|_{L^2}^2dt\nonumber\\
&\quad+C\sup_{t_0\le t\le T}\big(\|\nabla^2d\|_{L^2}^2\big)\int_{t_0}^T\big(\|\nabla\theta\|_{L^2}^2
+\|\nabla^2\theta\|_{L^2}^2\big)dt\nonumber\\
&\quad+C\sup_{t_0\le t\le T}\big(\|\nabla^2d\|_{L^2}\|\theta\|_{L^\infty}\|\nabla d\|_{L^2}\big)
\int_{t_0}^T\|\nabla^2d\|_{L^2}^2dt+C\nonumber\\
&\le C.
\end{align}
By Gagliardo-Nirenberg inequality, \eqref{tt}, and \eqref{dd1}, we get
\begin{align}
&\int_0^Tt^\frac12\|\theta\|_{L^\infty}^2\big(\|\nabla^2d\|_{L^2}^2+\|\nabla d\|_{L^2}\|\nabla^2d\|_{L^2}^3\big)dt\nonumber\\
&\le C\sup_{0\le t\le T}\big(t\|\nabla^2d\|_{L^2}^2\big)^\frac12
\int_{0}^T\big(\|\nabla\theta\|_{L^2}^2+\|\nabla^2\theta\|_{L^2}^2\big)dt\nonumber\\
&\quad+C\sup_{0\le t\le T}\big(t^\frac12\|\nabla^2d\|_{L^2}\big)\int_0^T\|\nabla^2d\|_{L^2}^2dt
\le C.
\end{align}
We infer from Gagliardo-Nirenberg inequality, \eqref{tt}, and \eqref{5.5} that
\begin{align}
&\int_0^Tt^\frac12\big(\|\theta\nabla u\|_{L^2}^2+\|\nabla u\|_{L^2}^2+\|\nabla u\|_{L^2}^4\|\nabla\theta\|_{L^2}^2\big)dx\nonumber\\
&\le C\sup_{0\le t\le T}\big(\|\theta\|_{L^\infty}+\|\nabla u\|_{L^2}^4+1\big)
\int_0^Tt^\frac12\big(\|\nabla u\|_{L^2}^2+\|\nabla\theta\|_{L^2}^2\big)dt
\le C.
\end{align}
From \eqref{w7}, \eqref{tt}, Gagliardo-Nirenberg inequality, and Young's inequality, we have
\begin{align}\label{5.17}
\|\nabla u\|_{L^4}^4&\le \|\nabla u\|_{L^2}\|\nabla u\|_{L^6}^3\nonumber\\
&\le C\|\nabla u\|_{L^2}\big(\|\sqrt{\rho}\dot{u}\|_{L^2}+\|\nabla\theta\|_{L^2}
+\|\nabla d\|_{L^\infty}\|\nabla^2d\|_{L^2}\big)^3\nonumber\\
&\le C\|\nabla u\|_{L^2}^4+C\|\sqrt{\rho}\dot{u}\|_{L^2}^4
+C\|\nabla\theta\|_{L^2}^4+C\|\nabla^2d\|_{L^2}^2\|\nabla^3d\|_{L^2}^2.
\end{align}
Hence, we obtain from \eqref{tt}, \eqref{5.5}, and Proposition \ref{p1} that
\begin{align}
\int_0^Tt^\frac12\|\nabla u\|_{L^4}^4dt
&\le \int_0^Tt^\frac12\|\nabla u\|_{L^2}\big(\|\sqrt{\rho}\dot{u}\|_{L^2}
+\|\nabla d_t\|_{L^2}+\|\nabla\theta\|_{L^2}+\|\nabla^2d\|_{L^2}\big)^3dt\nonumber\\
&\le C\sup_{0\le t\le T}\big[t^\frac12\big(\|\nabla d_t\|_{L^2}^2+\|\sqrt{\rho}\dot{u}\|_{L^2}^2\big)\big]
\int_0^T\big(\|\sqrt{\rho}\dot{u}\|_{L^2}^2+\|\nabla d_t\|_{L^2}^2\big)dt\nonumber\\
&\quad+C\int_0^Tt^\frac12(\|\nabla u\|_{L^2}^4+\|\nabla\theta\|_{L^2}^4
+\|\nabla^2d\|_{L^2}^4)ds\nonumber\\
&\le C\Bbb E_0^\frac23\sup_{0\le t\le T}\big[t^\frac12\big(\|\nabla d_t\|_{L^2}^2
+\|\sqrt{\rho}\dot{u}\|_{L^2}^2\big)\big]
+C\sup_{0\le t\le T}t\|\nabla^2d\|_{L^2}^2\int_0^T\|\nabla^2d\|_{L^2}^2ds\nonumber\\
&\quad+C\sup_{0\le t\le T}\big(\|\nabla u\|_{L^2}^2+\|\nabla\theta\|_{L^2}^2\big)
\int_0^Tt^\frac12\big(\|\nabla u\|_{L^2}^2+\|\nabla\theta\|_{L^2}^2\big)dt\nonumber\\
&\le C_{10}\Bbb E_0^\frac23\sup_{0\le t\le T}\big[t^\frac12\big(\|\nabla d_t\|_{L^2}^2
+\|\sqrt{\rho}\dot{u}\|_{L^2}^2\big)\big]+C.
\end{align}
From \eqref{f49}, \eqref{tt}, and Proposition \ref{p1}, we have
\begin{align}
\int_0^Tt^\frac12\|\nabla^3d\|_{L^2}^4dt
&\le C\sup_{0\le t\le T}\big[t^\frac12\big(\|\nabla d_t\|_{L^2}^2+\|\nabla u\|_{L^2}^2
+\|\nabla^2d\|_{L^2}^2\big)\big]\int_0^T\|\nabla^3d\|_{L^2}^2dt\nonumber\\
&\le C_{11}\Bbb E_0^\frac23\sup_{0\le t\le T}\big[t^\frac12\big(\|\nabla d_t\|_{L^2}^2+\|\nabla u\|_{L^2}^2
+\|\nabla^2d\|_{L^2}^2\big)\big].
\end{align}
Meanwhile, we infer from \eqref{tt}, \eqref{dd1}, and \eqref{5.5} that
\begin{align}\label{k14}
&C\int_0^Tt^\frac12\big(\|\nabla^2d\|_{L^2}^4+\|\nabla u\|_{L^2}^4\|\nabla^2d\|_{L^2}^2+\|\nabla d\|_{L^2}\|\nabla^2d\|_{L^2}^3\big)dt\nonumber\\
&\le C\sup_{0\le t\le T}\big(t^\frac12\|\nabla^2d\|_{L^2}\big)\int_{0}^T\|\nabla^2d\|_{L^2}^2dt
+C\sup_{0\le t\le T}\big(\|\nabla u\|_{L^2}^2\|\nabla^2d\|_{L^2}^2\big)
\int_0^Tt^\frac12\|\nabla u\|_{L^2}^2dt\nonumber\\
&\quad+C\sup_{0\le t\le T}\|\nabla d\|_{L^2}
\sup_{0\le t\le T}\big(t^\frac12\|\nabla^2d\|_{L^2}\big)
\int_0^T\|\nabla^2d\|_{L^2}^2dt\nonumber\\
&\le C.
\end{align}
Substituting \eqref{k7}--\eqref{k14} into \eqref{k6} and choosing $\delta$ suitably small, we arrive at
\begin{align}\label{5.15}
&\sup_{0\le t\le T}\big[t^\frac12\big(\|\nabla u\|_{L^2}^2
+\|\nabla\theta\|_{L^2}^2+\|\sqrt{\rho}\dot{u}\|_{L^2}^2+\|\nabla d_t\|_{L^2}^2\big)\big]\nonumber\\
&\quad+\int_0^Tt^\frac12\big(\|\sqrt{\rho}\dot{u}\|_{L^2}^2
+\|\nabla d_t\|_{L^2}^2+\|\nabla^3d\|_{L^2}^2+\|\sqrt{\rho}\dot{\theta}\|_{L^2}^2\big)dt\nonumber\\
&\quad+\int_0^Tt^\frac12\big(\|\nabla\dot{u}\|_{L^2}^2
+\|d_{tt}\|_{L^2}^2+\|\nabla^2d_t\|_{L^2}^2\big)dt\nonumber\\
&\le C,
\end{align}
provided that  $$\Bbb E_0\le \varepsilon_0'\triangleq\min\left\{\varepsilon_4, \Big(\frac{c_2}{8C_{10}}\Big)^\frac32, \Big(\frac{\mu}{2C_{11}}\Big)^\frac32\right\}.$$
Consequently, the desired \eqref{1.11} follows from \eqref{5.5} and \eqref{5.15}.

\textbf{Step 4. Estimate for $\|\nabla d_t(t)\|_{L^2}$}.
We derive from \eqref{rrf}, Sobolev's inequality, Gagliardo-Nirenberg inequality, Young's inequality, and \eqref{tt} that
\begin{align*}
&\frac{d}{dt}\|\nabla d_t\|_{L^2}^2
+\|d_{tt}\|_{L^2}^2+\|\nabla^2d_t\|_{L^2}^2\nonumber\\
&\le C\|\nabla\dot{u}\|_{L^2}^2\|\nabla d\|_{L^2}\|\nabla^2d\|_{L^2}
+C\|\nabla u\|_{L^2}^2\|\nabla u\|_{L^4}^2\|\nabla d\|_{L^3}\|\nabla^3d\|_{L^2}
\nonumber\\
&\quad+C\|\nabla u\|_{L^2}^2\|\nabla d_t\|_{L^2}\|\nabla^2d_t\|_{L^2}
+C\|\nabla^2d\|_{L^2}^4\|\nabla d_t\|_{L^2}^2
+C\|\nabla^2d\|_{L^2}^2\|\nabla d_t\|_{L^2}\|\nabla^2d_t\|_{L^2}\nonumber\\
&\le \frac12\|\nabla^2d_t\|_{L^2}^2
+C\big(\|\nabla\dot{u}\|_{L^2}^2+\|\nabla u\|_{L^2}^2\big)\|\nabla^2d\|_{L^2}
+C\|\nabla u\|_{L^4}^4
+C\big(\|\nabla^2d\|_{L^2}^2+\|\nabla u\|_{L^2}^2\big)\|\nabla d_t\|_{L^2}^2,
\end{align*}
which multiplied by $t$ together with \eqref{5.17} leads to
\begin{align}\label{5.16}
&\frac{d}{dt}\big(t\|\nabla d_t\|_{L^2}^2\big)
+t\|d_{tt}\|_{L^2}^2+t\|\nabla^2d_t\|_{L^2}^2\nonumber\\
&\le C\big(\|\nabla^2d\|_{L^2}^2+\|\nabla u\|_{L^2}^2\big)
\big(t\|\nabla d_t\|_{L^2}^2\big)
+Ct\big(\|\nabla\dot{u}\|_{L^2}^2+\|\nabla u\|_{L^2}^2\big)\|\nabla^2d\|_{L^2} \notag \\
& \quad +Ct\big(\|\nabla u\|_{L^2}^4+\|\sqrt{\rho}\dot{u}\|_{L^2}^4
+\|\nabla\theta\|_{L^2}^4\big)+Ct\|\nabla^2d\|_{L^2}^2\|\nabla^3d\|_{L^2}^2
+\|\nabla d_t\|_{L^2}^2.
\end{align}
We obtain from \eqref{3.52}, Sobolev's inequality, and \eqref{tt} that
\begin{align}\label{5.18}
\|\nabla d_t\|_{L^2}^2 & \le C\|\nabla^3d\|_{L^2}^2
+C\||u||\nabla^2d|\|_{L^2}^2+C\||\nabla u||\nabla d|\|_{L^2}^2
+C\||\nabla d|^3\|_{L^2}^2+C\||\nabla d||\nabla^2d|\|_{L^2}^2 \nonumber\\
&\le C\|\nabla^3d\|_{L^2}^2
+C\|u\|_{L^6}^2\|\nabla^2d\|_{L^3}^2+C\|\nabla d\|_{L^\infty}^2\|\nabla u\|_{L^2}^2
+C\|\nabla d\|_{L^6}^6+C\|\nabla d\|_{L^3}^2\|\nabla^2d\|_{L^6}^2 \nonumber\\
&\le C\|\nabla^3d\|_{L^2}^2
+C\|\nabla u\|_{L^2}^2\|\nabla^2d\|_{H^1}^2+C\|\nabla^2d\|_{L^2}
\|\nabla^3d\|_{L^2}\|\nabla u\|_{L^2}^2+C\|\nabla d\|_{H^1}^2\|\nabla^2d\|_{H^1}^2
\nonumber\\
&\le C\|\nabla^2 d\|_{H^1}^2+C\|\nabla u\|_{L^2}^2.
\end{align}
By virtue of \eqref{dd1}, \eqref{5.5}, \eqref{5.15}, and \eqref{tt}, we have
\begin{align*}
& \int_0^Tt\big(\|\nabla\dot{u}\|_{L^2}^2+\|\nabla u\|_{L^2}^2\big)\|\nabla^2d\|_{L^2}dt
\leq \sup_{0\le t\le T}\big(t^{\frac12}\|\nabla^2d\|_{L^2}\big)
\int_0^Tt^{\frac12}\big(\|\nabla\dot{u}\|_{L^2}^2+\|\nabla u\|_{L^2}^2\big)dt\leq C, \\
& \int_0^Tt\big(\|\nabla u\|_{L^2}^4+\|\sqrt{\rho}\dot{u}\|_{L^2}^4
+\|\nabla\theta\|_{L^2}^4\big)dt \\
& \leq \sup_{0\le t\le T}\big[t^{\frac12}\big(\|\nabla u\|_{L^2}^2+\|\sqrt{\rho}\dot{u}\|_{L^2}^2
+\|\nabla\theta\|_{L^2}^2\big)\big]
\int_0^Tt^{\frac12}\big(\|\nabla u\|_{L^2}^2+\|\sqrt{\rho}\dot{u}\|_{L^2}^2
+\|\nabla\theta\|_{L^2}^2\big)dt\leq C, \\
& \int_0^Tt\|\nabla^2d\|_{L^2}^2\|\nabla^3d\|_{L^2}^2dt
\leq \sup_{0\le t\le T}\big(t\|\nabla^2d\|_{L^2}^2\big)
\int_0^T\|\nabla^3d\|_{L^2}^2dt\leq C.
\end{align*}
These estimates combined with \eqref{5.16}, \eqref{5.18}, Gronwall's inequality, and \eqref{tt} yields that
\begin{align}\label{5.19}
\sup_{0\le t\le T}\big(t\|\nabla d_t\|_{L^2}^2\big)
+\int_0^Tt\big(\|d_{tt}\|_{L^2}^2
+\|\nabla^2d_t\|_{L^2}^2\big)dt\le C.
\end{align}

\textbf{Step 5. Estimate for $\|\sqrt{\rho}\dot{\theta}(t)\|_{L^2}$}. We obtain from \eqref{f80}, H\"older's inequality, Sobolev's inequality, \eqref{3.49}, \eqref{f47}, \eqref{5.17}, \eqref{3.74}, and \eqref{tt} that
\begin{align*}
& c_v\frac{d}{dt}\|\sqrt{\rho}\dot{\theta}\|_{L^2}^2+\kappa
\|\nabla\dot{\theta}\|_{L^2}^2\nonumber\\
&\le C\big(\|\sqrt{\rho}\theta\|_{L^2}\|\sqrt{\rho}\theta\|_{L^6}
+\|\nabla u\|_{L^2}\|\nabla u\|_{L^6}
+\|\nabla^2d\|_{L^2}^2+\|\nabla^3d\|_{L^2}^2\big)\nonumber\\
&\quad \times
\big(\|\sqrt{\rho}\dot{\theta}\|_{L^2}^2+\|\nabla^2\theta\|_{L^2}^2
+\|\nabla^2d_t\|_{L^2}^2+\|\nabla d_t\|_{L^2}^2
+\|\nabla u\|_{L^4}^4+\|\nabla\dot{u}\|_{L^2}^2\big)\nonumber\\
&\le C\big(\|\nabla\theta\|_{L^2}^2+\|\nabla u\|_{L^2}^2
+\|\sqrt{\rho}\dot{u}\|_{L^2}^2+\|\nabla^2d\|_{L^2}^2+\|\nabla^3d\|_{L^2}^2\big)\nonumber\\
&\quad\times\big(\|\sqrt{\rho}\dot{\theta}\|_{L^2}^2
+\|\nabla\theta\|_{L^2}^2+\|\nabla u\|_{L^2}^2
+\|\sqrt{\rho}\dot{u}\|_{L^2}^2+\|\nabla^2d\|_{L^2}^2
+\|\nabla^2d_t\|_{L^2}^2+\|\nabla d_t\|_{L^2}^2+\|\nabla\dot{u}\|_{L^2}^2\big),
\end{align*}
which multiplied by $t$ gives that
\begin{align}\label{5.20}
& c_v\frac{d}{dt}\big(t\|\sqrt{\rho}\dot{\theta}\|_{L^2}^2\big)
+\kappa t\|\nabla\dot{\theta}\|_{L^2}^2\nonumber\\
&\le C\big[t^{\frac12}\big(\|\nabla\theta\|_{L^2}^2+\|\nabla u\|_{L^2}^2
+\|\sqrt{\rho}\dot{u}\|_{L^2}^2\big)\big]\nonumber\\
&\quad\times t^{\frac12}\big(\|\sqrt{\rho}\dot{\theta}\|_{L^2}^2
+\|\nabla\theta\|_{L^2}^2+\|\nabla u\|_{L^2}^2
+\|\sqrt{\rho}\dot{u}\|_{L^2}^2+\|\nabla^2d_t\|_{L^2}^2+\|\nabla d_t\|_{L^2}^2+\|\nabla\dot{u}\|_{L^2}^2\big) \notag \\
& \quad +C\big(t\|\nabla^2d\|_{L^2}^2\big)
\big(\|\nabla\theta\|_{L^2}^2+\|\nabla u\|_{L^2}^2+\|\nabla\dot{\theta}\|_{L^2}^2
+\|\nabla\dot{u}\|_{L^2}^2+\|\nabla^2d\|_{H^1}^2+\|\nabla d_t\|_{H^1}^2\big)
 \notag \\
& \quad +C\|\nabla^2d\|_{H^1}^2
\big[t\big(\|\nabla d_t\|_{L^2}^2+\|\nabla^2d_t\|_{L^2}^2\big)\big]
+C\|\nabla\dot{\theta}\|_{L^2}^2 \notag \\
& \quad +C\big(\|\sqrt{\rho}\dot{\theta}\|_{L^2}^2
+\|\nabla\theta\|_{L^2}^2+\|\nabla u\|_{L^2}^2
+\|\sqrt{\rho}\dot{u}\|_{L^2}^2\big)\big(t\|\nabla^3d\|_{L^2}^2\big) \notag \\
& \quad
+C\big[t\big(\|\nabla d_t\|_{L^2}^2+\|\nabla^2d\|_{L^2}^2\big)\big]\|\nabla\dot{u}\|_{L^2}^2
\end{align}
where we have used \eqref{f49} and the following facts
\begin{align*}
\|\sqrt{\rho}\dot{u}\|_{L^2}^2
\leq \|\sqrt{\rho}\|_{L^3}^2\|\dot{u}\|_{L^6}^2
\leq C(\|\rho_0\|_{L^1},\bar\rho)\|\nabla\dot{u}\|_{L^2}^2,\
\|\sqrt{\rho}\dot{\theta}\|_{L^2}^2
\leq \|\sqrt{\rho}\|_{L^3}^2\|\dot{\theta}\|_{L^6}^2
\leq C(\|\rho_0\|_{L^1},\bar\rho)\|\nabla\dot{\theta}\|_{L^2}^2,
\end{align*}
due to H\"older's inequality, \eqref{lz3.1}, \eqref{3.49}, and Sobolev's inequality.
Consequently, integrating \eqref{5.20} over $[0,T]$, we deduce from \eqref{5.15}, \eqref{5.5}, \eqref{dd1}, \eqref{tt}, and \eqref{5.19} that
\begin{align}\label{5.21}
\sup_{0\le t\le T}\big(t\|\sqrt{\rho}\dot{\theta}\|_{L^2}^2\big)
+\int_0^Tt\|\nabla\dot{\theta}\|_{L^2}^2dt\le C.
\end{align}

\textbf{Step 6. Estimate for $\|\nabla^3d(t)\|_{L^2}$ and $\|\nabla^2\theta(t)\|_{L^2}$}.  In view of \eqref{f49} and \eqref{tt}, we get that
\begin{align*}
\|\nabla^3d\|_{L^2}^2
\le C\|\nabla d_t\|_{L^2}^2+C\|\nabla^2d\|_{L^2}^2,
\end{align*}
which together with \eqref{5.19} and \eqref{dd1} leads to
\begin{align}\label{5.22}
\sup_{0\le t\le T}\big(t\|\nabla^3d\|_{L^2}^2\big)\le C.
\end{align}
We obtain from \eqref{a1}$_3$, the standard $L^2$-estimates of elliptic equations, Sobolev's inequality, Gagliardo-Nirenberg inequality, \eqref{tt}, and \eqref{f50} that
\begin{align*}
\|\nabla^2\theta\|_{L^2}^2&\le C\big(\|\sqrt{\rho}\dot{\theta}\|_{L^2}^2+\|\theta\nabla u\|_{L^2}^2+\|\nabla u\|_{L^4}^4
+\|\Delta d+|\nabla d|^2d\|_{L^4}^4\big)\nonumber\\
&\le C\big(\|\sqrt{\rho}\dot{\theta}\|_{L^2}^2+\|\theta\|_{L^\infty}^2\|\nabla u\|_{L^2}^2+\|\nabla u\|_{L^2}\|\nabla u\|_{L^6}^3
+\|\nabla^2d\|_{L^4}^4+\|\nabla d\|_{L^8}^8\big)\nonumber\\
&\le C\|\sqrt{\rho}\dot{\theta}\|_{L^2}^2+C\|\nabla\theta\|_{L^2}\|\nabla^2\theta\|_{L^2}\|\nabla u\|_{L^2}^2
+C\|\nabla^2d\|_{H^1}^2\nonumber\\
&\quad
+C\|\nabla u\|_{L^2}\big(\|\sqrt{\rho}\dot{u}\|_{L^2}^3
+\|\nabla d_t\|_{L^2}^3+\|\nabla\theta\|_{L^2}^3+\|\nabla^2d\|_{L^2}^3\big)
\nonumber\\
&\le C\|\sqrt{\rho}\dot{\theta}\|_{L^2}^2
+\frac12\|\nabla^2\theta\|_{L^2}^2
+C\|\nabla\theta\|_{L^2}^2\|\nabla u\|_{L^2}^4+C\|\nabla^2d\|_{H^1}^2
\nonumber\\ &\quad +C\big(\|\nabla u\|_{L^2}^4+\|\sqrt{\rho}\dot{u}\|_{L^2}^4
+\|\nabla d_t\|_{L^2}^4+\|\nabla\theta\|_{L^2}^4+\|\nabla^2d\|_{L^2}^4\big).
\end{align*}
This along with \eqref{5.21}, \eqref{dd1}, \eqref{5.22}, and \eqref{5.15} implies that
\begin{align*}
\sup_{0\le t\le T}\big(t\|\nabla^2\theta\|_{L^2}^2\big)
&\le C\sup_{0\le t\le T}\big(t\|\sqrt{\rho}\dot{\theta}\|_{L^2}^2\big)
+C\sup_{0\le t\le T}\big(t\|\nabla^2d\|_{H^1}^2\big) \nonumber\\
&\quad
+C\sup_{0\le t\le T}\big[t^\frac12\big(\|\nabla u\|_{L^2}^2
+\|\nabla\theta\|_{L^2}^2+\|\sqrt{\rho}\dot{u}\|_{L^2}^2
+\|\nabla d_t\|_{L^2}^2\big)\big]^2\le C,
\end{align*}
which combined with \eqref{dd1}, \eqref{5.19}, \eqref{5.21}, and \eqref{5.22} yields the desired \eqref{1.12}.
\hfill $\Box$

\end{document}